\documentclass[12pt]{article}
\usepackage[utf8]{inputenc}
\usepackage[T1]{fontenc}
\usepackage{amsmath, amsfonts, amssymb, amsthm, dsfont, mathtools, mathrsfs}
\usepackage[hypertexnames=false]{hyperref}
\hypersetup{
	colorlinks   = true, 
	urlcolor     = blue, 
	linkcolor    = blue, 
	citecolor    = red 
}
\usepackage[nameinlink]{cleveref}
\usepackage{xcolor,colortbl} 
\usepackage{graphicx}
\usepackage{subfig}
\usepackage{caption}
\usepackage{subcaption}
\usepackage{authblk}
\usepackage[thinlines]{easytable}
\usepackage{tikz}
\usetikzlibrary{arrows}
\usepackage{pdflscape} 
\usepackage{nccmath} 
\usepackage{calc}
\usepackage{accents}
\usepackage{multirow}
\usepackage{todonotes}
\usepackage{makecell}

\theoremstyle{plain}
	\newtheorem{theorem}{Theorem}[section]
	\newtheorem{proposition}[theorem]{Proposition}
	\newtheorem{corollary}[theorem]{Corollary}
	
	\newtheorem{assumption}[theorem]{Assumption}
	\newtheorem{lemma}[theorem]{Lemma}
\theoremstyle{definition}
	\newtheorem{remark}[theorem]{Remark}
	
	\newtheorem{definition}[theorem]{Definition}

\theoremstyle{plain}

\usepackage[margin=1in]{geometry}
\providecommand{\keywords}[1]
{
	\small	
	\textbf{Keywords---} #1
}

\usepackage[
	style=authoryear,
	backend=biber,
	dashed=false,
	maxcitenames=2,
	mincitenames=1,
	maxbibnames=99,
	giveninits=true,
	uniquelist=false,
	uniquename=false
]{biblatex}

\DeclareNameAlias{sortname}{family-given}

\DeclareNameFormat{labelname}{%
	\nameparts{#1}%
	\usebibmacro{name:family}
	{\namepartfamily}
	{}
	{\namepartprefix}
	{\namepartsuffix}%
	\usebibmacro{name:andothers}%
}

\AtEveryBibitem{\clearfield{issn}}
\AtEveryBibitem{\clearfield{note}}
\DeclareFieldFormat[article,incollection,inproceedings]{title}{#1}
\DeclareFieldFormat{journaltitle}{\mkbibemph{#1}} 
\DeclareFieldFormat{title}{#1}                    
\renewbibmacro{in:}{}                             

\DeclareCiteCommand{\cite}
	{\usebibmacro{prenote}}
	{%
		\bibhyperref{%
			\printnames{labelname}%
			\setunit{\space}%
 			\printtext[brackets]{%
				\usebibmacro{cite:labeldate+extradate}%
			}%
		}%
	}
	{\multicitedelim}
	{\usebibmacro{postnote}} 
\newcommand{\hquad}{\hspace{0.5em}}

\newcommand{\Romannum}[1]{\textup{\uppercase\expandafter{\romannumeral#1}}}
\newcommand{\romannum}[1]{\textup{\lowercase\expandafter{\romannumeral#1}}}

\newcommand{\set}[1]{\left\{ #1 \right\}}
\newcommand{\dist}{{\operatorname{dist}}}
\newcommand{\op}{{\operatorname{op}}}
\newcommand{\F}{{\operatorname{F}}}
\newcommand{\HS}{{\operatorname{HS}}}
\newcommand{\diag}{{\operatorname{diag}}}
\newcommand{\spec}{{\operatorname{spec}}}
\newcommand{\tr}{{\operatorname{tr}}}
\newcommand{\var}{\mathrm{Var}}
\newcommand{\cov}{\mathrm{Cov}}

\renewcommand{\div}{\mathrm{div}}

\renewcommand{\P}{\mathbb{P}}
\newcommand{\dP}{d\mathbb{P}}
\newcommand{\OP}{O_{\mathbb{P}}}

\newcommand{\LPn}{{\mathbb{L}^2(\mathbb{P}_n)}}
\newcommand{\LP}{{\mathbb{L}^2(\mathbb{P})}}
\newcommand{\Exp}{\mathrm{Exp}}
\newcommand{\PT}{\mathrm{PT}}
\newcommand{\HP}{\mathrm{HP}}
\newcommand{\HPT}{\mathrm{HPT}}
\newcommand{\Jac}{\mathrm{Jac}}

\newcommand{\E}{\mathbb{E}}
\newcommand{\M}{\mathbb{M}}
\newcommand{\N}{\mathbb{N}}

\newcommand{\R}{\mathbb{R}}
\newcommand{\X}{\mathbb{X}}
\newcommand{\e}{\varepsilon}

\newcommand{\bL}{\bar{L}}

\newcommand{\bW}{\bar{W}}

\newcommand{\bgamma}{\bar{\gamma}}
\newcommand{\bmu}{\bar{\mu}}
\newcommand{\blambda}{\bar{\lambda}}

\newcommand{\bvarphi}{\bar{\varphi}}

\newcommand{\cF}{\mathcal{F}}

\newcommand{\cH}{\mathcal{H}}

\newcommand{\cJ}{\mathcal{J}}

\newcommand{\cM}{\mathcal{M}}

\newcommand{\cP}{\mathcal{P}}
\newcommand{\cQ}{\mathcal{Q}}
\newcommand{\cR}{\mathcal{R}}
\newcommand{\cS}{\mathcal{S}}
\newcommand{\cT}{\mathcal{T}}

\newcommand{\cX}{\mathcal{X}}
\newcommand{\cY}{\mathcal{Y}}
\newcommand{\cZ}{\mathcal{Z}}

\newcommand{\hA}{\hat{A}}
\newcommand{\hD}{\hat{D}}

\newcommand{\hH}{\hat{H}}

\newcommand{\hL}{\hat{L}}

\newcommand{\hS}{\hat{S}}
\newcommand{\hU}{\hat{U}}
\newcommand{\hV}{\hat{V}}
\newcommand{\hW}{\hat{W}}
\newcommand{\hu}{\hat{u}}
\newcommand{\hv}{\hat{v}}
\newcommand{\hPi}{\hat{\Pi}}
\newcommand{\hPhi}{\hat{\Phi}}
\newcommand{\hPsi}{\hat{\Psi}}

\newcommand{\hOmega}{\hat{\Omega}}

\newcommand{\hTheta}{\hat{\Theta}}
\newcommand{\halpha}{\hat{\alpha}}
\newcommand{\hbeta}{\hat{\beta}}
\newcommand{\hrho}{\hat{\rho}}

\newcommand{\hlambda}{\hat{\lambda}}

\newcommand{\hpsi}{\hat{\psi}}
\newcommand{\hvarphi}{\hat{\varphi}}
\newcommand{\hxi}{\hat{\xi}}
\newcommand{\heta}{\hat{\eta}}

\newcommand{\hbfH}{\hat{\mathbf{H}}}

\newcommand{\hcH}{\hat{\mathcal{H}}}

\newcommand{\hcP}{\hat{\mathcal{P}}}

\newcommand{\hcS}{\hat{\mathcal{S}}}

\newcommand{\tf}{\tilde{f}}

\newcommand{\tq}{\tilde{q}}

\newcommand{\trho}{\tilde{\rho}}
\newcommand{\tphi}{\tilde{\phi}}

\newcommand{\tfg}{\tilde{\mathfrak{g}}}

\newcommand{\tmbQ}{\tilde{\mathbb{Q}}}
\newcommand{\tbfG}{\tilde{\mathbf{G}}}

\newcommand{\fg}{\mathfrak{g}}
\newcommand{\fh}{\mathfrak{h}}
\newcommand{\fq}{\mathfrak{q}}
\newcommand{\fw}{\boldsymbol{\omega}}
\newcommand{\hfw}{\hat{\boldsymbol{\omega}}}

\newcommand{\msG}{\mathscr{G}}

\newcommand{\mbH}{\mathbb{H}}

\newcommand{\bfA}{\mathbf{A}}

\newcommand{\bfD}{\mathbf{D}}

\newcommand{\bfG}{\mathbf{G}}
\newcommand{\bfH}{\mathbf{H}}
\newcommand{\bfI}{\boldsymbol{\mathrm{I}}}
\newcommand{\bfK}{\mathbf{K}}
\newcommand{\bfL}{\mathbf{L}}
\newcommand{\bfO}{\mathbf{O}}

\newcommand{\bfQ}{\mathbf{Q}}

\newcommand{\bfW}{\mathbf{W}}

\newcommand{\bfPi}{\mathbf{\Pi}}

\newcommand{\bfsigma}{\boldsymbol{\sigma}}
\newcommand{\bfvarsigma}{\boldsymbol{\varsigma}}
\newcommand{\bfkappa}{\boldsymbol{\kappa}}
\newcommand{\bfvarkappa}{\boldsymbol{\varkappa}}

\newcommand{\msF}{\mathscr{F}}

\newcommand{\heat}{\text{heat}}

\newcommand{\unnor}{\text{un}} 
\newcommand{\un}{\text{un}}

\newcommand{\specd}{{\spec_{\downarrow}}}

\newcommand{\Pn}{{\mathbb{P}_n}}
 
\newcommand{\sn}[1]{[\![#1]\!]_n} 
\newcommand{\tmax}{{\mathrm{max}}}
\newcommand{\tmin}{{\mathrm{min}}}

\newcommand{\mycase}[1]{\vspace{.5em}\noindent\underline{{#1}:}\hquad}

\newcommand{\gammaJ}{\gamma_{\raisebox{-2pt}{$\hspace{-0.3mm}\scriptstyle \cJ$}}}

\newcommand{\IF}{{\rm IF}}
\newcommand{\Teps}{F_\e} 
\newcommand{\Tuneps}{{F^\un_\e}}
\newcommand{\LPt}{{\mathbb{L}^2(\mathbb{P}_t)}} 
\newcommand{\citep}[1]{(\cite{#1})}

\usepackage[dvipsnames]{xcolor}
\usepackage[normalem]{ulem}

\newcommand*{\mydot}[1]{\accentset{\scalebox{0.8}{\large\bfseries .}}{#1}}
\usepackage{enumerate}

\title{First-order asymptotic expansions for spectral convergence of normalized graph Laplacians constructed from samples on manifolds}
\author{Eunseong Bae\thanks{Research Institute of Applied Statistics, Sungkyunkwan University.
		Email: \texttt{esbae@skku.edu}}\,}
\author{Wolfgang Polonik\thanks{Department of Statistics, University of California, Davis.
		Email: \texttt{wpolonik@ucdavis.edu}}}
\affil{} 
\date{} 

\begin{document}
\maketitle

\begin{abstract}
	Under i.i.d. uniform sampling from a closed Riemannian manifold embedded in Euclidean space,
	we derive first-order asymptotic expansions for the eigenvalues and eigenvectors of symmetrically normalized graph Laplacians constructed via Gaussian kernels and provide concentration inequalities for the remainder terms.
	As a direct consequence, we establish the corresponding weak convergence results.
	Unlike most related works, our analysis does not impose a multiplicity-one assumption on the eigenvalues, thereby extending the applicability of the framework to general manifolds with potentially nontrivial eigenvalue multiplicities.
\end{abstract}


\tableofcontents

\section{Introduction}

The spectral properties of graph Laplacians play an important role in manifold learning. 
Prominent examples include Laplacian eigenmaps (\cite{belkin2003laplacian}), diffusion maps (\cite{coifman2006diffusion}), spectral clustering (\cite{vonluxburg2007tutorial}), and the heat kernel signature (\cite{sun2009concise}).
These methods leverage the eigenvalues and eigenvectors of graph Laplacians to capture the geometric structure of the data and information about the underlying manifold.

From an asymptotic perspective, when data are sampled from a manifold, the graph Laplacian is known to approximate the Laplace-Beltrami (LB) operator on the manifold. Accordingly, its eigenvalues and eigenvectors converge to the corresponding spectral quantities of the LB operator. Understanding this asymptotic behavior is therefore essential for relating the spectral properties of graph Laplacians to the underlying manifold. This has motivated a substantial body of work on the asymptotic analysis of graph Laplacians and their spectral properties.

Early contributions established consistency for the spectral convergence of graph Laplacians under appropriate sampling and scaling regimes (\cite{belkin2006convergence, luxburg2008consistency}).
Subsequent work has sought to establish explicit convergence rates and derive sharper probabilistic bounds (\cite{singer2017spectral, dunson2021spectral, calder2022improved, calder2022lipschitz, cheng2022eigenconvergence, wahl2024kernelbased}).
Recently, minimax rates for the spectral convergence have been proposed by \cite{chaubet2025minimax} and \cite{trillos2026minimax}.

Despite this substantial progress, the asymptotic distributional theory for the spectral convergence has not been fully studied.
To the best of our knowledge, \cite{li2025central} is the only paper establishing weak convergence of the eigenvalues of unnormalized graph Laplacians; however, this study is restricted to simple eigenvalues (i.e., multiplicity one), and the weak convergence of eigenvectors is not covered there.

The main contribution of this paper is to derive first-order asymptotic expansions for both eigenvalues and eigenvectors of symmetrically normalized graph Laplacians under uniform i.i.d. sampling from a closed (i.e., compact without boundary) manifold, and to provide concentration inequalities for the corresponding remainder terms.
These expansions characterize the leading-order stochastic fluctuations of the spectral quantities.
More importantly, based on these expansions, we establish weak convergence of both the eigenvalues and eigenvectors of normalized graph Laplacians, thereby providing a new result in the asymptotic theory of graph Laplacians.

A distinguishing feature of our analysis is that we do not impose the multiplicity-one assumption on the eigenvalues of the LB operator, which is commonly assumed in related work.
This is particularly important in many common manifold settings where eigenvalue multiplicities naturally arise.
For example, on a circle, every nonzero eigenvalue of the LB operator has multiplicity two.
More generally, higher-dimensional spheres exhibit eigenvalues with even larger multiplicities.

From a technical perspective, the presence of eigenvalue multiplicities leads to substantial differences in the analysis of spectral convergence. For simple eigenvalues, the corresponding eigenvectors are uniquely defined up to sign, and hence their behavior under perturbations can be more readily controlled.
In contrast, for an eigenvalue of nontrivial multiplicity, the associated eigenvectors are not uniquely determined and may change substantially even under small perturbations.

Consequently, a more delicate analysis based on operator perturbation theory is required to characterize the spectral behavior of graph Laplacians. Along these lines, some existing works, such as \cite{wahl2024kernelbased} and \cite{chaubet2025minimax}, adopt operator perturbation-theoretic approaches to study the spectral convergence of graph Laplacians, but do not establish their weak convergence.

The paper is organized as follows.
In Section~\ref{gauss:prelim}, we introduce notation and basic background.
Section~\ref{gauss:main} presents our main results on the first-order expansion of the eigenvalues and eigenvectors. 
Throughout Sections~\ref{smooth}-\ref{graph}, we prove the main results.
Technical proofs and auxiliary results are provided in the appendices.

\section{Preliminaries} \label{gauss:prelim}

We introduce the notation and mathematical preliminaries that will be used throughout the paper.

{\vspace{.5em}\bf Basic notation.} 
We denote by $\R$ and $\N$ the sets of real numbers and natural numbers, respectively.
For $a,b \in \R$, we write $a \wedge b = \min\{a,b\}$.
The Kronecker delta $\delta_{k \ell}$ is defined by $\delta_{k \ell} =1$ if $k=\ell$, and $\delta_{k \ell} =0$ otherwise.
For an index set $\cJ \subset \N$, we write $\{a_k\}_{k \in \cJ} = \{a_k: k \in \cJ\}$.
When $\cJ = \{1,\dots,n\}$ with $n \in \N$, we simply write $\{a_k\}_{k \in [n]}$.

{\vspace{.5em}\bf Matrix Analysis.} 
We let $\|\cdot\|_{\R^m}$ denote the Euclidean 2-norm on $\R^m$, and let $I_m$ denote the $m \times m$ identity matrix. 
For a matrix $A$, we let $\|A\|_\op$, $\|A\|_\F$ and $\tr(A)$ denote its operator norm (w.r.t. 2-norm), Frobenius norm and trace, respectively. Of course, $\tr(A)$ only applies when $A$ is square.

{\vspace{.5em}\bf Hilbert spaces.} 
For a Hilbert space $\mbH$, we denote its inner product by $\langle\cdot,\cdot\rangle_\mbH$ and the induced norm by $\|\cdot\|_\mbH$.
The operator norm and the Hilbert-Schmidt of an operator $\cH: \mbH \rightarrow \mbH$ is denoted by $\|\cH\|_{\op,\mbH}$ and $\|\cH\|_{\HS,\mbH}$, respectively.
For $u,v \in \mbH$, the tensor product $u \otimes_\mbH v: \mbH \rightarrow \mbH$ is defined by
$ (u \otimes_\mbH v) w = \langle v, w \rangle_\mbH u.$
We denote by $\bfI_\mbH:\mbH \to \mbH$ the identity operator on $\mbH$.

{\vspace{.5em}\bf $\mathbf{L^2}$ spaces.} 
Given a probability space $(\X,\P)$, $\LP$ denotes the space of real-valued functions $f$ on $\X$ such that $\int_\X f^2 \, d\P < \infty,$ equipped with the inner product $\langle f, g \rangle_\LP = \int_\X f g \, d\P.$

The empirical measure associated with the observations $X_1, \dots, X_n \in \X$ is denoted by $\P_n$, 
defined for any real-valued measurable function $f$ on $\X$ as
$\int_\X f \, d\P_n = \frac{1}{n} \sum_{i=1}^n f(X_i).$
Correspondingly, $\LPn$ denotes the Euclidean space $\R^n$ endowed with the inner product
$ \langle u, v \rangle_\LPn = \frac{1}{n} \sum_{i=1}^{n} u_i v_i.$

We also define the sampling or evaluation operator $\sn{\cdot}$ as a mapping from real-valued functions on $\X$ to $\LPn$: $\sn{f} = (f(X_i))_{i\in[n]}.$ This satisfies the relation
$ \langle \sn{f}, \sn{g} \rangle_\LPn = \frac{1}{n} \sum_{i=1}^{n} f(X_i) g(X_i),$
which is the empirical counterpart of the population inner product $\langle f, g\rangle_\LP$.
For brevity, we write $\P f = \int_\X f \, d\P$, $\P_n f = \int_\X f \, d\P_n$, $(\P_n - \P)f = \P_n f - \P f.$

We denote by $\mathscr{G}_\P$ the generalized $\P$-Brownian bridge.
This is a centered Gaussian process indexed by functions in $\LP$ with covariance
$ \cov (\mathscr{G}_\P(f), \mathscr{G}_\P (g)) = \P(fg) - (\P f) \, (\P g).$ 
$\mathscr{G}_\P$ arises from the weak limit of the empirical process $\{\sqrt{n} (\P_n - \P) f: f \in \msF \}$ for some appropriate function class $\msF$, i.e., 
$ \sqrt{n} (\P_n-\P) f \stackrel{\mathrm{weakly}}{\rightarrow} \msG_\P f$ for $f \in \msF$, where $\stackrel{\mathrm{weakly}}{\rightarrow}$ denotes weak convergence.

{\vspace{.5em}\bf Manifolds.}
A manifold of dimension $d$ is a topological space that is locally homeomorphic to the Euclidean space $\R^d$.
A Riemannian manifold is a smooth manifold equipped with a smooth, positive-definite metric tensor that defines an inner product on each tangent space.
We denote the Laplace-Beltrami operator on a Riemannian manifold $\M$ by $\Delta_\M  = - \div_\M \nabla_\M$, where $\div_\M$ and $\nabla_\M$ denote the divergence operator and the gradient operator on $\M$, respectively.
Note that we include the negative sign so that $\Delta_\M$ is self-adjoint and positive semi-definite with respect to the $L^2$-inner product induced by the volume measure.	

\section{Main results} \label{gauss:main}

\subsection{Setting}

Throughout the paper, we assume the following:
\begin{assumption} \label{gauss:assump}	
	$\M$ is a $d$-dimensional closed Riemannian sub-manifold isometrically embedded in $\R^p$.
	$X_1, \dots, X_n$ are i.i.d. samples drawn uniformly from $\M$ with respect to its volume measure.
	For simplicity, we assume the volume of $\M$ is 1, so that the uniform sampling probability measure $\P$ coincides with the volume measure.
\end{assumption}

For constructing the graph Laplacian, we use the Gaussian kernel of the form
\begin{equation*}
	\fg_\e (x,y) = \frac{1}{(4\pi \e)^{d/2}} \exp \bigg( -\frac{\|x-y\|^2_{\R^p}}{4\e} \bigg), \quad x,y \in \R^p,
\end{equation*}
where $\e>0$ is the bandwidth.
The Gaussian kernel is a popular choice for constructing graph Laplacians; see, e.g., \cite{belkin2006convergence, cheng2022eigenconvergence, dunson2021spectral, wahl2024kernelbased}.
We adopt this kernel in our analysis since it is highly related to the heat kernel $\fh_\e$ with diffusion time $\e$ associated with $\Delta_\M$, and this relationship enables a more straightforward application of the operator-convergence approach.
See Section~\ref{smooth} for further details.

For constructing graph Laplacians, we define the following quantities: 
\begin{equation*}
	\begin{aligned}
		\hrho_{n,\e} (x) &= \int_\M \fg_\e (x,y) \, d \P_n(y), \quad
		\rho_\e (x) = \int_\M \fg_\e (x,y) \, d \P(y), \\
		\hfw_{n,\e} (x,y) &= \frac{\fg_\e(x,y)}{\sqrt{\hrho_{n,\e}(x)}\sqrt{\hrho_{n,\e}(y)}}, \quad
		\fw_\e (x,y) = \frac{\fg_\e(x,y)}{\sqrt{\rho_\e(x)}\sqrt{\rho_\e(y)}}.
	\end{aligned}
\end{equation*}
The symmetrically normalized Gaussian graph Laplacian $\hL_{n,\e} \in \R^{n \times n}$ is defined as
\begin{equation*}
	\hL_{n,\e} = \frac{I_n - \hW_{n,\e}}{\e},
	\quad\text{where}\hquad
	\hW_{n,\e}  = \bigg( \frac{\hfw_{n,\e} (X_i,X_{i'})}{n} \bigg)_{i,i'\in [n]}.
\end{equation*}
Equivalently, $\hL_{n,\e}$ can be defined in terms of the affinity matrix $\hA_{n,\e} = (\frac{1}{n} \fg_\e (X_i, X_{i'}) )_{i,i' \in [n]}$ and the degree matrix $\hD_{n,\e} = \diag(\hrho_{n,\e}(X_i))_{i \in [n]} = (\text{$i$-th row sum of $\hA_{n,\e}$})_{i \in [n]}$ as
$$ \hL_{n,\e} = \frac{I_n - \hD^{-1/2}_{n,\e} \hA_{n,\e} \hD^{-1/2}_{n,\e}}{\e}.$$
When clear from context, we refer to $\hL_{n,\e}$ simply as the graph Laplacian.

We also consider the smoothed graph Laplacian $\bfL_\e: \LP \to \LP ,$ which is the continuum limit of $\hL_{n,\e}$, defined as
\begin{equation*}
	\bfL_\e = \frac{\bfI_\LP - \bfW_{\e}}{\e},
	\quad\text{where}\hquad 
	\bfW_\e f (x) = \int_\M \fw_\e (x,y) f(y) \dP (y).
\end{equation*}

The eigenpairs of $\hL_{n,\e}$, $\bfL_\e$, $\Delta_\M$  are denoted by 
\begin{equation*}
	\{(\hlambda_{k,n,\e}, \hvarphi_{k,n,\e})\}_{k \in [n]}, \quad
	\{(\lambda_{k,\e}, \phi_{k,\e})\}_{k \in \N}, \quad 
	\{(\lambda_k, \phi_k)\}_{k \in \N},	
\end{equation*}
respectively.
The eigenvalues $\{\hlambda_{k,n,\e}\}_{k\in[n]}$, $\{\lambda_{k,\e}\}_{k\in\N}$, $\{\lambda_k\}_{k\in\N}$ are non-negative and are arranged in non-decreasing order, with repetition allowed accordingly to their finite multiplicity.
The eigenvectors $\{\hvarphi_{k,n,\e}\}_{k \in [n]}$ and
the eigenfunctions $\{\phi_k\}_{k \in \N}$, $\{\phi_{k,\e}\}_{k \in \N}$ are $\LPn$- and $\LP$-orthonormal, respectively, that is, 
$\langle \hvarphi_{k,n,\e}, \hvarphi_{\ell,n,\e} \rangle_\LPn 
=  \langle \phi_{k,\e}, \phi_{\ell,\e} \rangle_\LP 
= \langle \phi_k, \phi_\ell \rangle_\LP
= \delta_{k \ell}$.

To handle the multiplicity of the eigenvalues of $\Delta_\M$, for each $j \in \N$, we define the index set
$\cJ_j = \{ k \in \N: \lambda_k = \lambda_j \}.$
By construction, $|\cJ_j|$ is the multiplicity of $\lambda_j$.
We denote the spectral gap $\gamma_j = \min_{k \notin \cJ_j} |\lambda_j - \lambda_k| = \min \{ \lambda_{j_+} - \lambda_j , \lambda_j - \lambda_{j_-} \}$, where $\lambda_{j_+}$ is the smallest eigenvalue of $\Delta_\M$ larger than $\lambda_j,$ and similarly $\lambda_{j_-}$ is the largest eigenvalue smaller than $\lambda_j$.
Throughout the paper, we fix $j \in \N$ unless stated otherwise.

Our main goal is to establish first-order asymptotic expansions for the spectral convergence of $\hL_{n,\e}$ to $\Delta_\M$, and to decompose the resulting error into two main sources: the stochastic error associated with $\hL_{n,\e}$ and $\bfL_\e$, and the bias associated with $\bfL_\e$ and $\Delta_\M$.

To be more precise, the stochastic error arises from the convergence of $\{(\hlambda_{k,n,\e}, \hvarphi_{k,n,\e}) \}_{k \in \cJ_j}$
to $\{(\lambda_{k,\e}, \phi_{k,\e}) \}_{k \in \cJ_j}$ as $n\to\infty$, with $n$ and $\e$ satisfying an appropriate scaling relation.
We note that $\{ \hvarphi_{k,n,\e} \}_{k \in \cJ_j}$ and $\{\phi_{k,\e}\}_{k \in \cJ_j}$ cannot be directly compared, as the former are vectors in $\R^n$, while the latter are functions in $\LPn$.
The precise formulation of this convergence and the corresponding comparison between the two will be discussed in Section~\ref{gauss:eigvec-secname}.

The bias arises from the convergence of the non-random quantity $\{(\lambda_{k,\e}, \phi_{k,\e}) \}_{k \in \cJ_j}$ to $\{(\lambda_{j}, \phi_{k}) \}_{k \in \cJ_j}.$ This convergence of course only depends on $\e$.

Before proceeding, we note that our results consist of finite sample bounds with constants depending on the underlying manifold $\M$ and the eigenvalues of $\Delta_\M$.
For convenience, we introduce the following notion.

\begin{definition}[$\M$- and $(\M,\lambda_j)$-dependence]
	A positive constant is called $\M$-dependent if it depends on $\M$ through its dimension $d$, volume, diameter, injectivity radius, as well as the uniform bounds on the Riemann curvature tensor and on the second fundamental form, together with their covariant derivatives.
	
	A positive constant is said to be $(\M,\lambda_j)$-dependent if it is $\M$-dependent and additionally depends on $\lambda_j$, $\lambda_{j+}$, $\lambda_{j-}$ and $|\cJ_j|$. 
\end{definition}

\subsection{First-order asymptotic expansions of eigenvalues of graph Laplacians}

We present the first-order asymptotic expansions of the eigenvalues $\{\hlambda_{k,n,\e}\}_{k \in \cJ_j}$ of $\hL_{n,\e}$.
Specifically, considering the multiplicity of $\lambda_j$, we consider sums of differences of eigenvalues over $\cJ_j$ as
\begin{equation}
	\sum_{k \in \cJ_j} (\hlambda_{k,n,\e} - \lambda_{k,\e})
	\quad\text{and}\quad 
	\sum_{k \in \cJ_j} (\hlambda_{k,n,\e} - \lambda_j).
	\label{gauss:eigval-form1}
\end{equation}
Sums of eigenvalues also appear in related literature on Laplacians without assumptions on multiplicities; see, e.g., \cite{chaubet2025minimax}.

The theoretical reason for considering these sums, rather than the vectors of differences $(\hlambda_{k,n,\e} - \lambda_{k,\e})_{k \in \cJ_j}$ and $(\hlambda_{k,n,\e} - \lambda_j)_{k \in \cJ_j}$, is the \emph{clustering} of distinct eigenvalues among $\{\lambda_{k,\e}\}_{k \in \cJ_j}$ of $\bfL_\e$.

More precisely, classical perturbation theory of operators typically considers a single eigenvalue with multiplicity, provided that the eigenvalue is sufficiently separated from the other eigenvalues; that is, the corresponding spectral gap should not be too small.
Under such a spectral gap condition, first-order asymptotic expansions can be formulated for the vector of differences between the perturbed eigenvalues and the associated single eigenvalue, with the dimension of the vector equal to the multiplicity of the single eigenvalue. To apply this theory to our setting, that is, to formulate first-order asymptotic expansions for the vector $(\hlambda_{k,n,\e} - \lambda_{k,\e})_{k \in \cJ_j}$, representing the stochastic error, we need to know the number of distinct eigenvalues among $\{\lambda_{k,\e}\}_{k \in \cJ_j}$, the multiplicities of each distinct eigenvalue, and the spectral gaps between them.
However, this information is generally unavailable.
The only thing we can guarantee is that $\{\lambda_{k,\e}\}_{k \in \cJ_j}$ are close to $\lambda_j$ as $\e \to 0$.
Moreover, in this case, the inner spectral gaps, namely, the spectral gaps between distinct eigenvalues within this collection, may become smaller than the rate of convergence to $\lambda_j$.
Overall, classical perturbation theory may not be directly applicable to our setting.

Recently, \cite{bae2026asymptotic} suggested considering the sum of eigenvalue differences rather than the vector to formulate first-order expansions.
In particular, to control the remainder term in this expansion, one does not require the inner spectral gaps; rather, the outer spectral gap needs to be sufficiently large. In our case, the outer spectral gap of $\{\lambda_{k,\e}\}_{k \in \cJ_j}$ is close to $\gamma_j$, the spectral gap of $\lambda_j$. This suggests that the proposed method is well suited to our setting.

The following theorem shows that the sums \eqref{gauss:eigval-form1} are well approximated by the stochastic term
\begin{equation} 
	(\P_n-\P) \bigg( -\sum_{k\in\cJ_j} \Delta_\M(\phi_k^2) \bigg).
	\label{gauss:eigval-form2}
\end{equation}
We regard this as the first-order leading term in the asymptotic expansion of \eqref{gauss:eigval-form1}.
The proof of the theorem is provided at the end of Section~\ref{graph:eigval}. We also provide a heuristic derivation of \eqref{gauss:eigval-form2} in Appendix~\ref{heuristic} providing insight into why this approximation term emerges.

\begin{theorem} \label{gauss:eigval}
	There exist $(\M,\lambda_j)$-dependent constants $r$, $C$, $C'$, $C''$, $C'''$ such that
	if $\e < r$ and $\frac{n\e^{d/2+2}}{\log(1/\e)} > C'''$, 
	we have with probability at least $1-\tau-C'\exp(-C''n\e^{d/2+2})$,
	\begin{multline}
		\Bigg| \sum_{k \in \cJ_j} (\hlambda_{k,n,\e} - \lambda_{k,\e}) - (\P_n-\P) \bigg( -\sum_{k\in\cJ_j} \Delta_\M(\phi_k^2) \bigg) \Bigg| 
		\\ \leq
		C \Bigg( \e^{5/6} \log(1/\e) \sqrt{\frac{\log(3/\tau)}{n}} + \frac{\log(1/(\e\wedge\tau))}{n \e^{d/2 + 2}} + \frac{\log^2(1/(\e\wedge\tau))}{n^2 \e^{d + 2}} \, \Bigg),
		\label{gauss:eigval-res1}
	\end{multline}
	and by replacing $\lambda_{k,\e}$ with $\lambda_j$ on the left-hand side,
	\begin{multline}
		\Bigg| \sum_{k \in \cJ_j} (\hlambda_{k,n,\e} - \lambda_j) - (\P_n-\P) \bigg( -\sum_{k\in\cJ_j} \Delta_\M(\phi_k^2) \bigg) \Bigg| 
		\\ \leq 
		C \Bigg( \e + \e^{5/6} \log(1/\e) \sqrt{\frac{\log(3/\tau)}{n}}+ \frac{\log(1/(\e\wedge\tau))}{n \e^{d/2 + 2}} + \frac{\log^2(1/(\e\wedge\tau))}{n^2 \e^{d + 2}} \, \Bigg),
		\label{gauss:eigval-res2}
	\end{multline}
\end{theorem}

\begin{remark} \label{gauss:eigval:remark1}
	From the above theorem, by choosing the optimal $\e \asymp n^{-\frac{2}{d+6}}$ (up to log factors), we obtain the error rate for estimating $\lambda_j$ as $O( n^{-\frac{1}{2}} + n^{-\frac{2}{d+6}} ) = O(n^{-\frac{2}{d+6}} )$.
	Here, the term $n^{-\frac{1}{2}}$ comes from the first-order leading term \eqref{gauss:eigval-form2}, while $n^{-\frac{2}{d+6}}$ arises from the bias and stochastic remainder term.
	Since the latter dominates the former, it determines the overall convergence rate.
	
	This rate is slower compared to recently established minimax rates for estimating the eigenvalues of weighted Laplace-Beltrami operators, which are $O(n^{-\frac{2}{d+4}})$ in the general manifold setting \citep{trillos2026minimax}, and $O(n^{-\frac{1}{2}})$ when the underlying manifold is assumed to be known \citep{chaubet2025minimax}.
	
	As seen above, our rate is dominated by the bias and stochastic remainder terms rather than the first-order leading term.
	The gap between our rate and the existing minimax rates suggests that sharper control of these terms may be possible.
	We leave this improvement for future work.
	
	A similar comparison applies to the eigenvector convergence results in Theorem~\ref{gauss:eigvec}.
\end{remark}

Below we provide two important consequences of Theorem~\ref{gauss:eigval}.

Note that the first-order leading term \eqref{gauss:eigval-form2} has zero expectation.
Hence, by applying standard high-probability arguments to Theorem~\ref{gauss:eigval}, we conclude that the expected value of $\sum_{k \in \cJ_j} \hlambda_{k,n,\e}$ is quite close to $\sum_{k \in \cJ} \lambda_{k,\e}$, as stated below.

\begin{corollary} \label{gauss:eigval:expectation}
	There exist $(\M,\lambda_j)$-dependent constants $r$, $C$, $C'$, $C''$, $C'''$ such that
	if $\e <r$ and $\frac{n\e^{d/2+2}}{\log(1/\e)} > C'''$, we have
	\begin{equation*}
		\bigg| \sum_{k\in\cJ_j} (\E \hlambda_{k,n,\e} - \lambda_{k,\e}) \bigg| 
		\leq C \bigg( \frac{\e^{5/6} \log(1/\e)}{\sqrt{n}} + \frac{\log(1/\e)}{n\e^{d/2+2}} \bigg),
	\end{equation*}
	and with probability at least $1-\tau-C'\exp(-C''n\e^{d/2+2})$,
	\begin{multline*}
		\Bigg| \sum_{k \in \cJ_j} (\hlambda_{k,n,\e} - \E\hlambda_{k,n,\e}) - (\P_n-\P) \bigg( -\sum_{k\in\cJ_j} \Delta_\M(\phi_k^2) \bigg) \Bigg| 
		\\ \leq C \Bigg( \e^{5/6} \log(1/\e) \sqrt{\frac{\log(4/\tau)}{n}}+ \frac{\log(1/(\e\wedge\tau))}{n \e^{d/2 + 2}} + \frac{\log^2(1/(\e\wedge\tau))}{n^2 \e^{d + 2}} \, \Bigg),
	\end{multline*}	
\end{corollary}

Now, we turn to the weak convergence of the eigenvalues.
From Theorem~\ref{gauss:eigval}, we see that as $\e \to 0$ and $\frac{n\e^{d+4}}{\log^2(1/\e)} \rightarrow \infty$,
\begin{equation*}
	\sqrt{n} \sum_{k \in \cJ_j} (\hlambda_{k,n,\e} - \lambda_{k,\e}) = \sqrt{n}(\P_n-\P) \bigg( -\sum_{k\in\cJ_j} \Delta_\M(\phi_k^2) \bigg) + \OP \Bigg(\e^{5/6} \log(1/\e) + \sqrt{\frac{\log^2(1/\e)}{n\e^{d+4}}} \Bigg).
\end{equation*}
Hence, by applying  Slutsky's theorem, we obtain the weak convergence results as follows: 

\begin{corollary} \label{gauss:eigval:weak}
	As $\e \rightarrow 0$ and $\frac{n\e^{d+4}}{\log^2(1/\e)} \rightarrow \infty$, we have
	\begin{equation}
		\sqrt{n} \sum_{k \in \cJ_j} (\hlambda_{k,n,\e} - \lambda_{k,\e})
		\stackrel{\mathrm{weakly}}{\rightarrow}
		\msG_\P \bigg( -\sum_{k \in \cJ_j} \Delta_\M (\phi_k^2) \bigg).
		\label{gauss:eigval:weak-res1}
	\end{equation}
	More generally, for distinct eigenvalues $\lambda_{j_1}, \dots, \lambda_{j_m}$,
	\begin{equation*}
		\bigg( \sqrt{n} \sum_{k \in \cJ_j} (\hlambda_{k,n,\e} - \lambda_{k,\e}) \bigg)_{j \in \{j_1,\dots,j_m\}}
		\stackrel{\mathrm{weakly}}{\rightarrow}
		\bigg( \msG_\P \bigg( -\sum_{k \in \cJ_j} \Delta_\M (\phi_k^2) \bigg)
		\bigg)_{j \in \{j_1,\dots,j_m\}}.
	\end{equation*}
	Additionally, by Corollary~\ref{gauss:eigval:expectation}, the above results continue to hold when $\lambda_{k,\e}$ is replaced with $\E \hlambda_{k,n,\e}$.
\end{corollary}

\begin{remark} \label{gauss:eigval:remark2}
	Under the same asymptotic regime and $\sqrt{n}$-scaling, \cite{li2025central} establishes asymptotic normality for simple eigenvalues of the unnormalized graph Laplacian.
	For comparison, under uniform sampling and when $\lambda_j$ is simple, Theorem~1 of \cite{li2025central} provides a mean-zero asymptotic normality with variance
	\begin{align*}
		\int_\M  \Big( \Delta_\M (\phi_j^2(x))  + \lambda_j -\lambda_j \phi^2_j(x) \Big)^2 \dP(x),
	\end{align*}
	which contains the additional term $\lambda_j -\lambda_j \phi^2_j$ in the integrand compared with our variance
	\begin{align*}
		\var \, \msG_\P \big( -\Delta_\M(\phi^2_j) \big) = \int_\M  \Big( \Delta_\M (\phi_j^2(x)) \Big)^2 \dP(x).
	\end{align*}
	The difference may be related to whether the graph Laplacian is unnormalized or normalized.
	We note that the variances cannot be directly compared, since the sign of the cross term involving $\lambda_j - \lambda_j \phi_j^2$ is not determined in general. We refer to the heuristics in Appendix~\ref{heuristic} for further insights into the relations.
\end{remark}

\begin{remark} \label{gauss:eigval:remark3}
	Unfortunately, our results do not guarantee the weak convergence of $\sqrt{n} \sum_{k \in \cJ_j} (\hlambda_{k,n,\e} - \lambda_j)$, 
	that is, when $\lambda_{k,\e}$ in \eqref{gauss:eigval:weak-res1} is replaced by the true $\lambda_j$.
	The bias term requires $\sqrt{n}\e\to0$, while the stochastic remainder requires $\frac{n\e^{d+4}}{\log^2(1/\e)}\to\infty$, and these conditions cannot be satisfied simultaneously. 
	The same issue also appears in \cite{li2025central}.
	
	As discussed in Remark~\ref{gauss:eigval:remark1}, obtaining sharper bounds for the bias and stochastic remainder terms, particularly those achieving the $n^{-1/2}$ rate in \cite{chaubet2025minimax}, may make it possible to extend the weak convergence when replacing $\lambda_{k,\e}$ with $\lambda_j$.
	However, such rates require the manifold to be known, and whether this extension is possible for unknown manifolds remains unclear.
	
	We refer to Corollary~\ref{gauss:eigvec:weak} for a similar issue concerning eigenvectors.
\end{remark}


\subsection{First-order asymptotic expansions of eigenvectors of graph Laplacians} \label{gauss:eigvec-secname}

In this subsection, we present the first-order asymptotic expansions of the eigenvectors $\{\hvarphi_{k,n,\e}\}_{k\in\cJ_j}$ of $\hL_{n,\e}$.
To obtain a theoretically clean expansion, particularly one that is suitable for establishing weak convergence, two issues need to be addressed.

First, the eigenvectors can be viewed as discrete counterparts of the corresponding eigenfunctions $\{\phi_{k,\e}\}_{k \in \cJ_j}$ and $\{\phi_{k}\}_{k \in \cJ_j}$.
Thus, we compare the eigenvectors with the evaluations of these eigenfunctions at the data points, that is, the evaluation vectors $\{\sn{\phi_{k,\e}}\}_{k\in\cJ_j}$ and $\{\sn{\phi_k}\}_{k\in\cJ_j}$.

Second, due to eigenvalue multiplicities, $\{\hvarphi_{k,n,\e}\}_{k \in \cJ_j}$ and $\{ \sn{\phi_{k,\e}} \}_{k\in\cJ_j}$ (and likewise $\{ \sn{\phi_k} \}_{k \in \cJ_j}$) do not necessarily admit a term-by-term correspondence for the same index $k \in \cJ_j$.
Therefore, rather than comparing individual vectors, it is more appropriate to consider quantities characterizing the alignment between the corresponding linear subspaces.

Considering these issues, we investigate the first-order asymptotic expansions 
\begin{equation}
	\sum_{k \in \cJ_j} \sum_{\ell \in \cJ_j} \langle \hvarphi_{k,n,\e}, \sn{\phi_{\ell,\e}} \rangle_\LPn^2 - |\cJ_j|
	\quad\text{and}\quad \sum_{k \in \cJ_j} \sum_{\ell \in \cJ_j} \langle \hvarphi_{k,n,\e}, \sn{\phi_{\ell}} \rangle_\LPn^2 - |\cJ_j|.	
	\label{gauss:eigvec-form1}
\end{equation}
Such quantities also appear in the PCA literature to measure the discrepancy between two linear subspaces spanned by the corresponding two sets of vectors; see, e.g., \cite{fan2019distributedb}. The following theorem shows that both quantities in \eqref{gauss:eigvec-form1} are well approximated by the stochastic term
\begin{equation*}
	(\P_n-\P) \bigg( \sum_{k \in \cJ_j} \phi_k^2 \bigg).
\end{equation*}
The proof is provided at the end of Section~\ref{graph:eigvec}.

\begin{theorem} \label{gauss:eigvec}
	There exist $(\M,\lambda_j)$-dependent constants $r$, $C$, $C'$, $C''$, $C'''$ such that
	if $\e < r$ and $\frac{n\e^{d/2+2}}{\log(1/\e)} > C'''$, 
	we have with probability at least $1-\tau-C'\exp(-C''n\e^{d/2+2})$, 
	\begin{multline}
		\Bigg| \sum_{k \in \cJ_j} \sum_{\ell \in \cJ_j} \langle \hvarphi_{k,n,\e}, \sn{\phi_{\ell,\e}} \rangle^2_\LPn - |\cJ_j|
		- (\P_n-\P) \bigg( \sum_{k \in \cJ_j} \phi_k^2 \bigg) \Bigg|
		\\ \leq C \bigg( \e \sqrt{\frac{\log(3/\tau)}{n}}+ \frac{\log(1/(\e\wedge\tau))}{n \e^{d/2 + 2}} + \frac{\log^2(1/(\e\wedge\tau))}{n^2 \e^{d + 2}} \bigg),
		\label{gauss:eigvec-res1}
	\end{multline}
	and by replacing $\phi_{\ell,\e}$ with $\phi_\ell$ on the left-hand side,
	\begin{multline}
		\Bigg| \sum_{k \in \cJ_j} \sum_{\ell \in \cJ_j} \langle \hvarphi_{k,n,\e}, \sn{\phi_\ell} \rangle^2_\LPn - |\cJ_j|
		- (\P_n-\P) \bigg( \sum_{k \in \cJ_j} \phi_k^2 \bigg) \Bigg|
		\\ \leq C \bigg( \e+ \e \sqrt{\frac{\log(3/\tau)}{n}}+ \frac{\log(1/(\e\wedge\tau))}{n \e^{d/2 + 2}} + \frac{\log^2(1/(\e\wedge\tau))}{n^2 \e^{d + 2}} \bigg).
		\label{gauss:eigvec-res2}
	\end{multline}	
\end{theorem}

Similar to the eigenvalues, we obtain the following two corollaries.
We also refer to Remarks~\ref{gauss:eigval:remark1} and \ref{gauss:eigval:remark3} for related issues concerning the convergence rates of eigenvectors.

\begin{corollary} \label{gauss:eigvec:expectation}
	There exist $(\M,\lambda_j)$-dependent constants $r$, $C$, $C'$, $C''$, $C'''$ such that
	if $\e < r$ and $\frac{n\e^{d/2+2}}{\log(1/\e)} > C'''$, we have
	\begin{equation*}
		\bigg| \sum_{k \in \cJ_j} \sum_{\ell \in \cJ_j} \E \langle \hvarphi_{k,n,\e}, \sn{\phi_{\ell,\e}} \rangle^2_\LPn - |\cJ_j| \bigg| 
		\leq C\bigg( \frac{\e}{\sqrt{n}} + \frac{\log(1/\e)}{n\e^{d/2+2}} \bigg),
	\end{equation*}
	and with probability at least $1-\tau-C'\exp(-C''n\e^{d/2+2})$,
	\begin{multline*}
		\Bigg| \sum_{k \in \cJ_j} \sum_{\ell \in \cJ_j} \Big( \langle \hvarphi_{k,n,\e}, \sn{\phi_{\ell,\e}} \rangle^2_\LPn - \E \langle \hvarphi_{k,n,\e}, \sn{\phi_{\ell,\e}} \rangle^2_\LPn \Big)
		- (\P_n-\P) \bigg( \sum_{k \in \cJ_j} \phi_k^2 \bigg) \Bigg|
		\\ \leq C \Bigg( \e \sqrt{\frac{\log(4/\tau)}{n}}+\Bigg(\sqrt{\frac{\log(1/(\e\wedge\tau))}{n \e^{d/2 + 2}}} + \frac{\log(1/(\e\wedge\tau))}{n \e^{d/2 + 1}}\Bigg)^2 \, \Bigg),
	\end{multline*}	
\end{corollary}

\begin{corollary} \label{gauss:eigvec:weak}
	As $\e \rightarrow 0$ and $\frac{n\e^{d+4}}{\log^2(1/\e)} \to \infty$, we have
	\begin{equation*}
		\sqrt{n}\,\Bigg( \sum_{k \in \cJ_j} \sum_{\ell \in \cJ_j} \langle \hvarphi_{k,n,\e}, \sn{\phi_{\ell,\e}} \rangle^2_\LPn - |\cJ_j| \Bigg)
		\stackrel{\mathrm{weakly}}{\rightarrow}	
		\msG_\P \Bigg( \sum_{k \in \cJ_j} \phi_k^2 \Bigg).
	\end{equation*}
	More generally, for distinct eigenvalues $\lambda_{j_1}, \dots, \lambda_{j_m}$,
	\begin{equation*}
		\Bigg( \sqrt{n}\,\Bigg( \sum_{k \in \cJ_j} \sum_{\ell \in \cJ_j} \langle \hvarphi_{k,n,\e}, \sn{\phi_{\ell,\e}} \rangle^2_\LPn - |\cJ_j| \Bigg) \Bigg)_{j \in \{j_1,\dots,j_m\}}
		\stackrel{\mathrm{weakly}}{\rightarrow}	
		\Bigg( \msG_\P \Bigg( \sum_{k \in \cJ_j} \phi_k^2 \Bigg) \Bigg)_{j \in \{j_1,\dots,j_m\}}.
	\end{equation*}
	Additionally, by Corollary~\ref{gauss:eigvec:expectation}, the above results hold when $|\cJ_j|$ is replace with $\E \langle \hvarphi_{k,n,\e}, \sn{\phi_{\ell,\e}} \rangle_\LPn^2$.
\end{corollary}

\subsection{Discussion}

We discuss possible directions for improving our current results in future research.

{\vspace{.5em}\bf Multiplicity estimation.}
For our results to be practically applicable, the multiplicities of the eigenvalues of the Laplace--Beltrami operator must be known, which is generally unavailable.
However, if a consistent estimator of the multiplicities is available, our results continue to hold with high probability by simply replacing the true multiplicities with their estimates.
Unfortunately, to the best of our knowledge, the estimation of these multiplicities in the manifold setting has not been studied.

Related questions have been considered in the spectral perturbation theory; see \cite{hall2009tierespecting} and Lemma~2.1 of \cite{koltchinskii2020efficient}.
See also Section~6.1.4 of \cite{jolliffe2002principal} for related sequential testing procedures.

{\vspace{.5em}\bf Non-Gaussian kernels.} 
Extending our results beyond the Gaussian kernel is interesting direction.
This extension appears feasible, as convergence results for graph Laplacians with more general kernels have been established.
In particular, simpler kernels like the $0-1$ kernel are of interest; see, e.g., \cite{calder2022improved}.
To apply our arguments, one needs to establish the operator convergence of the corresponding smoothed Laplacian to the smoothed heat graph Laplacian, defined below.
We establish this convergence under the Gaussian kernel, but extending to a broader class of kernels remains an open problem.

{\vspace{.5em}\bf Non-uniform sampling.} 
For simplicity, we assume uniform sampling throughout the paper, although this assumption is unrealistic in practice.
Our results can be extended to non-uniform sampling settings; however, 
in this case, the limiting operator is replaced by a weighted Laplace-Beltrami operator rather than the standard Laplace-Beltrami 
operator, and additional technical details are required.

A more interesting direction would be to consider density-corrected graph Laplacians, as studied in \cite{cheng2022eigenconvergence} and \cite{dunson2021spectral}, whose limiting operator is the standard Laplace-Beltrami operator even under non-uniform sampling.
In particular, this framework captures the geometric information of the underlying manifold without the density bias. We leave this for future work.

\subsection{Outline of proofs} \label{gauss:outline}

The remaining sections are devoted to proving the main results.
Since the proofs are lengthy and involve many steps, we provide a brief summary of the key ideas of each section.

In Section \ref{smooth}, we establish the spectral convergence of $\bfL_\e$ to $\Delta_\M$.
To this end, we introduce the smoothed heat graph Laplacian $\bfL^\heat_\e : \LP \to \LP$, defined by
\begin{equation}
	\bfL^\heat_\e = \frac{\bfI_\LP - \bfH_\e}{\e},
	\quad\text{where}\quad \bfH_\e f (x) = \int_\M \fh_\e (x,y) \frac{f(x)-f(y)}{\e} \dP (y).
	\label{gauss:outline-eq1}
\end{equation}
Since $\bfL^\heat_\e$ shares a similar spectral structure with $\Delta_\M$,
the spectral convergence of $\bfL_\e$ to $\bfL_\e^\heat$ implies that of $\bfL_\e$ to $\Delta_\M$.
The convergence is achieved via relative form boundedness stated below.

Throughout Sections \ref{popul} and \ref{graph}, we study the spectral convergence of $\hL_{n,\e}$ to $\bfL_\e$. 
We split the analysis by introducing the population-normalized graph Laplacian $\bL_{n,\e} \in \R^{n \times n}$, defined by
\begin{equation}
	\bL_{n,\e} = \frac{I_n - \bW_{n,\e}}{\e},
	\quad\text{where}\quad
	\bW_{n,\e} = \bigg( \frac{\fw_\e(X_i,X_{i'})}{n} \bigg)_{i,i' \in [n]}.
	\label{gauss:outline-eq2}
\end{equation}

In Section \ref{popul}, we establish the spectral convergence of $\bL_{n,\e}$ to $\bfL_\e$. 
This is achieved by treating $\bW_{n,\e}$ as a kernel Gram matrix with deterministic kernel $\fw_\e$, 
and $\bfW_\e$ as the corresponding integral operator, and then transferring these results to $\bL_{n,\e}$ and $\bfL_\e$.

In Section \ref{graph}, we study the spectral convergence of $\hL_{n,\e}$ to $\bL_{n,\e}$.
Since $\hL_{n,\e}$ and $\bL_{n,\e}$ differs only in the normalization factor ($\hrho_{n,\e}$ versus $\rho_\e$), the uniform convergence of $\hrho_{n,\e}$ to $\rho_\e$ implies the convergence of $\hL_{n,\e}$ to $\bL_{n,\e}$ in operator norm, from which, the spectral convergence follows.
Technical difficulties stem from the randomness of the eigenvalues and eigenvectors of $\bL_{n,\e}$; this is resolved by replacing them with the deterministic counterparts from $\bfL_\e$ and $\Delta_\M$ under acceptable error bounds.
Finally, combining these results with those from Sections \ref{smooth} and \ref{popul}, we prove our main results, Theorems~\ref{gauss:eigval} and \ref{gauss:eigvec}.

\section{Spectral convergence of smoothed graph Laplacians} \label{smooth} 

In this section, we study the spectral convergence of $\bfL_\e$ to $\Delta_\M$.
Technical issues arise from the fact that $\Delta_\M$ is an unbounded operator on $\LP$.
This unboundedness makes it difficult to apply standard notions of operator norm convergence, which is usually required in operator spectral theory, and complicates the formulation of spectral convergence.

To address this issue, instead of $\Delta_\M$, we use the smoothed heat graph Laplacian $\bfL_\e^\heat$, defined in $\eqref{gauss:outline-eq1}$.
The advantage of using $\bfL_\e^\heat$ instead of $\Delta_\M$ lies in both its spectral properties and its operator-theoretic features.

First, $\bfL_\e^\heat$ shares the same eigenfunctions as $\Delta_\M$.
Specifically, for each $k \in \N$, it holds that $\bfH_\e \phi_k = e^{-\e\lambda_k} \phi_k$ (e.g., see Theorem~10.3 of \cite{grigoryan2009heat}), and consequently,
\begin{equation*}
	\bfL_\e^\heat \phi_k = \frac{1-e^{-\e\lambda_k}}{\e} \phi_k.
\end{equation*}
Moreover, $(1-e^{-\e\lambda_k})/\e$, the eigenvalue of $\bfL_\e^\heat$ associated with $\phi_k$, converges to $\lambda_k$ as $\e \to 0$.
This implies that asymptotically $\bfL_\e^\heat$ preserves the spectral structure of $\Delta_\M$.

Second, unlike $\Delta_\M$, $\bfL_\e^\heat$ is a bounded operator on $\LP$ with operator norm $\|\bfL_\e^\heat\|_{\op,\LP} = 1/\e$.
Although this operator norm diverges as $\e \to 0$, its boundedness for fixed $\e$ allows the application of tools from operator theory when analyzing spectral convergence.

All of this then allows us to establish the spectral convergence of $\bfL_\e$ to $\bfL_\e^\heat$, and to then deduce the spectral convergence to $\Delta_\M$ for sufficiently small $\e$.
For the first step, the convergence is established via the framework of \emph{relative form boundedness}.
Specifically, in Proposition~\ref{smooth:rel}, we establish
\begin{equation*}
	\big| \langle (\bfL_\e- \bfL_\e^\heat) f, f \rangle_\LP \big|
	\leq C \e \log^2(1/\e) \Big( \langle \bfL^\heat_\e f, f \rangle_\LP + \|f\|^2_\LP \Big)
	\quad\text{for all $f \in \LP$.}
\end{equation*}
Heuristically, relative form boundedness means the ratio $\frac{\bfL_\e - \bfL_\e^\heat}{\bfL_\e^\heat}$ is stable even when the magnitude of $\|\bfL_\e^\heat\|_{\op,\LP}$ is large, and moreover, the ratio tends to zero as $\e \to0$ in our case.

As $\bfL_\e^\heat$ is asymptotically unbounded, this notion of relative control is necessary in our case.
Similar relative form boundedness arguments appear in existing literature on unnormalized graph Laplacians; see Theorem~4.1 of \cite{belkin2006convergence} or Theorem~4 of \cite{wahl2024kernelbased}.
We adapt these approaches to our setting.

In the following, we first establish this relative form boundedness, and in the subsequent subsections we use it to derive the spectral convergence of $\bfL_\e$ to $\Delta_\M$.
For the connection between relative form boundedness and spectral convergence, we refer to \cite{bae2026asymptotic}.A brief summary of the relevant argument is provided in Appendix~\ref{gauss:spec}.

\subsection{Relative form boundedness of smoothed graph Laplacian}

To establish the relative form boundedness of $\bfL_\e$ to $\bfL_\e^\heat$, we introduce the unnormalized smoothed graph Laplacian $\bfL_\e^\unnor$, defined by
\begin{equation*}
	\bfL^\unnor_\e f (x) = \int_\M \fg_\e (x,y) \frac{f(x)-f(y)}{\e} \dP (y).
\end{equation*}

The proof proceeds in three steps. First, we establish the relative form boundedness of $\bfL_\e^\unnor$ to $\bfL_\e^\heat$ (Proposition~\ref{smooth:rel1}).
Second, we show the relative form boundedness of $\bfL_\e$ to $\bfL_\e^\unnor$ (Proposition~\ref{smooth:rel2}). 
Finally, by combining these two results, we deduce the relative form boundedness of $\bfL_\e$ to $\bfL_\e^\heat$ (Proposition~\ref{smooth:rel}).

\subsubsection*{Relative form boundedness of $\bfL_\e^\unnor$ to $\bfL_\e^\heat$}

We begin with the relative form boundedness of $\bfL_\e^\unnor$ to $\bfL_\e^\heat$.

\begin{proposition} \label{smooth:rel1}
	There exist $\M$-dependent constants $r$ and $C$ such that $\e < r$ implies
	\begin{equation*}
		\big| \langle (\bfL_\e^\unnor- \bfL_\e^\heat) f, f \rangle_\LP \big|
		\leq C \e \log^2(1/\e) \Big( \langle \bfL^\heat_\e f, f \rangle_\LP + \|f\|^2_\LP \Big)
		\quad\text{for all $f \in \LP$.}
	\end{equation*}	
\end{proposition}

\begin{proof}[Proof of Proposition~\ref{smooth:rel1}]
	We follow arguments similar to those in Theorem~4 of \cite{wahl2024kernelbased}, where the discrete version is proven.
	
	By the definition of $\bfL_\e^\unnor$, we have
	\begin{align*}
		\langle \bfL_\e^\unnor f, f \rangle_\LP
		&= \int_\M \bigg( \int_\M \fg_\e (x,y) \frac{f(x)-f(y)}{\e} \dP(y) \bigg) f(x) \dP(x)
		\\ &= \int_\M \int_\M \fg_\e (x,y) \frac{f(x)-f(y)}{\e} f(x) \dP(y) \dP(x).
	\end{align*}
	Using the identity $(f(x)-f(y)) f(x) = (f(x)-f(y))^2 - (f(y)-f(x))f(y)$, we obtain
	\begin{align*}
		\langle \bfL_\e^\unnor f, f \rangle_\LP
		&= \int_\M \int_\M \fg_\e (x,y) \frac{(f(x)-f(y))^2}{\e} f(x) \dP(y) \dP(x)
		\\&\quad - \int_\M \int_\M \fg_\e (x,y) \frac{f(y)-f(x)}{\e} f(y) \dP(y) \dP(x)
		\\&= \int_\M \int_\M \fg_\e (x,y) \frac{(f(x)-f(y))^2}{\e} f(x) \dP(y) \dP(x) - \langle \bfL_\e^\unnor f, f \rangle_\LP,
	\end{align*}
	and this yields
	\begin{equation}
		2 \langle \bfL_\e^\unnor f, f \rangle_\LP = \int_\M \int_\M \fg_\e (x,y) \frac{(f(x)-f(y))^2}{\e} \dP(x) \dP(y).
		\label{smooth:rel1:eq3}
	\end{equation}
	
	On the other hand, since $\M$ is compact, by stochastic completeness (e.g., see Section 8.4 of \cite{grigoryan2009heat}), it holds that $ \int_\M \fh(x,y) \dP(y) = 1$ for all $x \in \M$.
	Consequently,
	\begin{align*}
		\bfL_\e^\heat f(x) 
		= \frac{f(x) - \int_\M \fh_\e (x,y) f(y) \dP (y)}{\e}
		= \int_\M \fh_\e (x,y) \frac{f(x)-f(y)}{\e} \dP(y).
	\end{align*}
	This implies $\bfL_\e^\heat$ can be obtained by replacing $\fg_\e$ in $\bfL_\e^\unnor$ by $\fh_\e$. Hence, similar to \eqref{smooth:rel1:eq3}, we obtain
	\begin{equation} 
		2 \langle \bfL_\e^\heat f, f \rangle_\LP = \int_\M \int_\M \fh_\e (x,y) \frac{(f(x)-f(y))^2}{\e} \dP(x) \dP(y).
		\label{smooth:rel1:eq4}
	\end{equation}
	
	Now, using \eqref{smooth:rel1:eq3} and \eqref{smooth:rel1:eq4}, we have
	\begin{align}
		2\, \big| \langle (\bfL_\e^\unnor - \bfL_\e^\heat) f, f \rangle \big|
		&\leq \int_\M \int_\M | \fg_\e (x,y) - \fh_\e (x,y) | \frac{(f(x)-f(y))^2}{\e} \dP(x) \dP(y).
		\label{smooth:rel1:eq1}
	\end{align}
	According to Corollary 3 in \cite{wahl2024kernelbased}, there exist $\M$-dependent constants $C_1$ and $r_1$ such that 
	\begin{equation}
		|\fg_\e (x,y) - \fh_\e (x,y)| \leq C_1 \e \log^2(1/\e) (\fh_\e (x,y)+\e),
		\label{smooth:rel1-gauss}
	\end{equation}
	for any $x,y \in \M$ and any $\e < r_1$.
	Thus, we obtain
	\begin{align}
		&\quad \int_\M \int_\M | \fg_\e (x,y) - \fh_\e (x,y)| \frac{(f(x)-f(y))^2}{\e} \dP(x) \dP(y) 
		\notag\\&\leq \int_\M \int_\M  C_1 \e \log^2(1/\e) (\fh_\e (x,y)+\e) \frac{(f(x)-f(y))^2}{\e} \dP(y) \dP(x)
		\notag\\&= C_1 \e \log^2(1/\e) \bigg( \int_\M \int_\M  \fh_\e (x,y) \frac{(f(x)-f(y))^2}{\e} \dP(y) \dP(x)
		+ \int_\M \int_\M (f(x)-f(y))^2 \dP(y) \dP(x) \bigg)
		\notag\\&= C_1 \e \log^2(1/\e) \bigg( 2 \langle \bfL_\e^\heat f, f \rangle_\LP
		+ \int_\M \int_\M (f(x)-f(y))^2 \dP(y) \dP(x) \bigg).
		\label{smooth:rel1:eq2}
	\end{align}
	Since
	$$ \int_\M \int_\M (f(x)-f(y))^2 \dP(y) \dP(x) \leq \int_\M \int_\M 2(f^2(x) + f^2(y)) \dP(x) \dP(y) = 4 \|f\|^2_\LP,$$
	by combining \eqref{smooth:rel1:eq1} and \eqref{smooth:rel1:eq2}, we conclude
	\begin{align*}
		2 | \langle (\bfL_\e^\unnor - \bfL_\e^\heat) f, f \rangle_\LP |
		\leq C_1 \e \log^2(1/\e) \Big( 2 \langle \bfL_\e^\heat f, f \rangle_\LP + 4 \|f\|^2_\LP \Big).
	\end{align*}
	This completes the proof by letting $r = r_1$ and $C = 2C_1$.
\end{proof}

\subsubsection*{Relative form boundedness of $\bfL_\e$ to $\bfL_\e^\unnor$}

To establish the relative form boundedness of $\bfL_\e$ to $\bfL_\e^\unnor$, 
we introduce the following multiplication operator:
\begin{equation*}
	\bfPi_\e f (x) = (\rho_\e (x)-1) f(x).
\end{equation*}
Using $\bfPi_\e$, the difference $\bfL_\e - \bfL^\unnor_\e$ can be approximated by a more tractable form, as stated below:

\begin{lemma} \label{smooth:rel2:aux1}
	If $ \| \bfPi_\e \|_{\op,\LP} < 1/2$, then
	\begin{equation}
		\bigg\| \bfL_\e - \bfL^\unnor_\e + \frac{\bfPi_\e \bfL^\unnor_\e + \bfL^\unnor_\e \bfPi_\e }{2} \bigg\|_{\op,\LP} 
		\leq \frac{17\|\bfL_\e^\unnor\|_{\op,\LP} \|\bfPi_\e\|^2_{\op,\LP}}{4}.
		\label{smooth:rel2:part1-res2}
	\end{equation}
\end{lemma}

\begin{proof}[Proof of Lemma \ref{smooth:rel2:aux1}]
	Let $f \in \LP$ and $x \in \M$.
	From 
	$(\bfI_\LP + \bfPi_\e)^{-1/2} f = f/\sqrt{\rho_\e},$
	we have
	\begin{align}
		(\bfI_\LP + \bfPi_\e)^{-1/2} \bfL^\unnor_\e (\bfI_\LP + \bfPi_\e)^{-1/2} f (x)
		&=(\bfI_\LP + \bfPi_\e)^{-1/2} \bfL^\unnor_\e \bigg( \frac{f}{\sqrt{\rho_\e}} \bigg) (x)
		= \frac{1}{\sqrt{\rho_\e (x)}} \bfL^\unnor_\e \bigg( \frac{f}{\sqrt{\rho_\e}} \bigg) (x).
		\label{smooth:rel2:aux1-eq1}
	\end{align}
	Here, by the definitions of $\bfL^\unnor_\e$, $\bfL_\e$ and $\rho_\e$, it follows that
	\begin{align}
		\bfL^\unnor_\e \bigg( \frac{f}{\sqrt{\rho_\e}} \bigg) (x)
		&= \int_\M \fg_\e (x,y) \frac{ \frac{f(x)}{\sqrt{\rho_\e (x)}} - \frac{f(y)}{\sqrt{\rho_\e (y)}} }{\e} \dP (y)
		\notag\\&= \frac{ \frac{f(x)}{\sqrt{\rho_\e (x)}} \int_\M \fg_\e (x,y) \dP(y)
			- \int_\M \frac{\fg_\e (x,y)}{\sqrt{\rho_\e}(y)} f(y) \dP(y) }{\e}
		\notag\\&= \frac{\sqrt{\rho_\e (x)} f(x) - \sqrt{\rho_\e(x)} \int_\M \frac{\fg_\e (x,y)}{\sqrt{\rho_\e}(x)\sqrt{\rho_\e}(y)} f(y) \dP(y)}{\e}
		\notag\\&= \sqrt{\rho_\e (x)} \bfL_\e f (x).
		\label{smooth:rel2:aux1-eq2}
	\end{align}
Combining \eqref{smooth:rel2:aux1-eq1} and \eqref{smooth:rel2:aux1-eq2}, 
	we obtain
	\begin{equation}
		\bfL_\e	= (\bfI_\LP + \bfPi_\e)^{-1/2} \bfL^\unnor_\e (\bfI_\LP + \bfPi_\e)^{-1/2}.
		\label{smooth:rel2:aux1-eq3}
	\end{equation}
	On the other hand, under the condition $\| \bfPi_\e \|_{\op,\LP} < 1/2$, applying \ref{gauss:sup:taylor} yields
	\begin{multline}
		\bigg\| (\bfI_\LP + \bfPi_\e)^{-1/2} \bfL^\unnor_\e (\bfI_\LP + \bfPi_\e)^{-1/2} - \bfL^\unnor_\e + \frac{\bfPi_\e \bfL^\unnor_\e + \bfL^\unnor_\e \bfPi_\e }{2} \bigg\|_{\op,\LP} 
		\\ \leq \frac{17\|\bfL_\e^\unnor\|_{\op,\LP} \|\bfPi_\e\|^2_{\op,\LP}}{4}.
		\label{smooth:rel2:aux1-eq4}
	\end{multline}
	Combining \eqref{smooth:rel2:aux1-eq3} and \eqref{smooth:rel2:aux1-eq4} completes the proof.
\end{proof}

For (\ref{smooth:rel2:part1-res2}) to be useful, we need to appropriately bound $\|\bfL_\e^\unnor\|_{\op,\LP}$ and $\|\bfPi_\e\|_{\op,\LP}$.
For $\|\bfL_\e^\unnor\|_{\op,\LP}$, we expect a bound analogous to $\|\bfL_\e^\heat\|_{\op,\LP}$,
and the following lemma confirms this:

\begin{lemma} \label{smooth:rel2:aux2}
	There exist $\M$-dependent constants $r$ and $C$ such that $\e < r$ implies
	\begin{equation*}
		\| \bfL_\e^\unnor \|_{\op,\LP} \leq \frac{C}{\e}.
	\end{equation*}
\end{lemma}

\begin{proof}[Proof of Lemma \ref{smooth:rel2:aux2}]
	Let $C_1$ and $r_1$ be as in Proposition~\ref{smooth:rel1} and let $\e < r_1$.
	Using the definition of the operator norm and Proposition~\ref{smooth:rel1}, we have
	\begin{align*}
		\| \bfL_\e^\unnor \|_{\op,\LP}
		&= \sup_{\substack{f:f\in\LP \\ \|f\|_\LP =1}} \langle \bfL_\e^\unnor f, f \rangle_\LP
		\\&\leq \sup_{\substack{f:f\in\LP \\ \|f\|_\LP =1}} \bigg( \langle \bfL_\e^\heat f, f \rangle_\LP
		+ C_1 \e \log^2(1/\e) \Big( \langle \bfL^\heat_\e f, f \rangle_\LP + \|f\|^2_\LP \Big) \bigg)
		\\&\leq \|\bfL_\e^\heat\|_{\op,\LP} + C_1\e \log^2(1/\e) \Big(\|\bfL_\e^\heat\|_{\op,\LP}+1\Big).
	\end{align*}
	Since $0\leq \lambda_k \nearrow \infty$ as $k \to \infty$, it follows that
	$\|\bfL_\e^\heat\|_{\op,\LP} = \sup\limits_{k \in \N} \frac{1-e^{-\e\lambda_k}}{\e} = \frac{1}{\e}. $
	Substituting this into the previous inequality gives
	\begin{equation*}
		\| \bfL_\e^\unnor \|_{\op,\LP}
		\leq \frac{1}{\e} + C_1 \log^2(1/\e) + C_1 \e \log^2(1/\e).
	\end{equation*}
    From this the assertion follows.
\end{proof}

For $\|\bfPi_\e\|_{\op,\LP}$, it follows directly from the definition that
$ \|\bfPi_\e\|_{\op,\LP} = \|\rho_\e -1 \|_\infty. $
So, it suffices to study $\rho_\e$.
The lemma below provides bounds on both $\|\rho_\e -1 \|_\infty$ and $\| \bfL_\e^\unnor \rho_\e \|_\infty$,
the latter of which will be used in subsequent arguments.
The proof is deferred to Appendix~\ref{smooth:KDE-proof}.

\begin{lemma} \label{smooth:KDE}
	There exist $\M$-dependent constants $r$ and $C$ such that $\e < r$ implies
	\begin{equation}
		\|\rho_\e -1 \|_\infty,\, \| \bfL_\e^\unnor \rho_\e \|_\infty \leq C \e. 
		\label{smooth:KDE-res}
	\end{equation}
\end{lemma}

We now state the relative form boundedness of $\bfL_\e$ to $\bfL_\e^\unnor$.

\begin{proposition} \label{smooth:rel2}
	There exist $\M$-dependent constants $r$ and $C$ such that $\e < r$ implies
	\begin{equation*}
		\big| \langle (\bfL_\e- \bfL_\e^\unnor) f, f \rangle \big|
		\leq C \e \Big( \langle \bfL^\unnor_\e f, f \rangle_\LP + \|f\|^2_\LP \Big)
		\quad\text{for all $f \in \LP$.}
	\end{equation*}
\end{proposition}

\begin{proof}[Proof of Proposition~\ref{smooth:rel2}]
	Let $f \in \LP$, and let $r_1$, $C_1$ and $r_2$, $C_2$ be as in Lemmas~\ref{smooth:rel2:aux2} and \ref{smooth:KDE}, respectively.
	Choose $r>0$ smaller than both $r_1$ and $r_2$, and also satisfying $C_2 r <1$.
	Then, for $\e<r$, it holds $\|\bfPi_\e\|_{\op,\LP} < 1/2$ and hence, by Lemma~\ref{smooth:rel2:aux1}, we have
	\begin{align*}
		\big|\langle (\bfL_\e - \bfL^\unnor_\e) f, f \rangle_\LP \big|
		\leq \bigg| \bigg\langle \frac{\bfPi_\e \bfL^\unnor_\e + \bfL^\unnor_\e \bfPi_\e }{2} f, f \bigg\rangle_\LP \bigg|
		+ \frac{17\|\bfL_\e^\unnor\|_{\op,\LP} \|\bfPi_\e\|^2_{\op,\LP} \|f\|^2_\LP}{4}.
	\end{align*}
	By applying Lemmas \ref{smooth:rel2:aux2} and \ref{smooth:KDE}, we further obtain
	\begin{equation} 
		|\langle (\bfL_\e - \bfL^\unnor_\e) f, f \rangle_\LP |
		\leq \bigg|\bigg\langle \frac{\bfPi_\e \bfL^\unnor_\e + \bfL^\unnor_\e \bfPi_\e }{2} f, f \bigg\rangle_\LP \bigg|
		+ \frac{17C_1 C_2^2 \e \|f\|^2_\LP}{4}.
		\label{smooth:rel2-eq1}
	\end{equation}
	To complete the proof, it thus suffices to appropriately bound the term
	\begin{equation*}
		\left\langle \frac{\bfPi_\e \bfL^\unnor_\e + \bfL^\unnor_\e \bfPi_\e }{2} f, f \right\rangle_\LP.
	\end{equation*}
	
	Observe that
	\begin{align*}
		\langle \bfPi_\e \bfL^\unnor_\e f, f \rangle_\LP
		&= \int_\M (\rho_\e (x) -1) \bigg(\int_\M \fg_\e (x,y) \frac{f(x) - f(y)}{\e} \dP(y) \bigg) f(x) \dP(x) 
		\\&= \int_\M (\rho_\e (x) -1) \bigg(\int_\M \fg_\e (x,y) \frac{(f(x) - f(y))}{\e} f(x)  \dP(y) \bigg) \dP(x) 
		\\&= \int_\M (\rho_\e (x) -1) \bigg(\int_\M \fg_\e (x,y) \frac{(f(x) - f(y))}{\e} f(y) \dP(y) \bigg) \dP(x) 
		\\&\quad + \int_\M (\rho_\e (x) -1) \bigg(\int_\M \fg_\e (x,y) \frac{(f(x) - f(y))^2}{\e} \dP(y) \bigg) \dP(x).
	\end{align*}
	By interchanging the order of integration, the first term on the right-hand side becomes
	\begin{align*}
		\int_\M (\rho_\e (x) -1) \left(\int_\M \fg_\e (x,y) \frac{f(x) - f(y)}{\e} f(y) \dP(y) \right)\dP(x)
		\\= \int_\M \left( \int_\M \fg_\e (x,y)  \frac{(\rho_\e (x) -1)(f(x) - f(y))}{\e} \dP(x)\right) f(y)\dP(y).
	\end{align*}
	Using the identity
	\begin{align*}
		(\rho_\e (x) -1) (f(x) - f(y))
		&= (\rho_\e (x) -1) f(x) - (\rho_\e (y) -1) f(y) + (\rho_\e (y) - \rho_\e(x))f(y)
		\\&= \bfPi_\e f (x) - \bfPi_\e f (y) + (\rho_\e (y) - \rho_\e(x))f(y),
	\end{align*}
	we obtain
	\begin{align*}
		&\quad \int_\M \left(\int_\M \fg_\e (x,y) \frac{(\rho_\e (x) -1)(f(x) - f(y))}{\e} \dP(x) \right) f(y) \dP(y)
		\\& = \int_\M \left( \int_\M \fg_\e (x,y) \frac{\bfPi_\e f(x) - \bfPi_\e f(y)}{\e} \dP(x)\right) f(y)\dP(y)
		\\&\quad + \int_\M \left( \int_\M \fg_\e (x,y) \frac{\rho_\e(x) - \rho_\e(y)}{\e} \dP(x)\right) f^2(y)\dP(y)
		\\& = - \langle \bfL_\e^\unnor \bfPi_\e f, f \rangle_\LP
		+ \int_\M \bfL_\e^\unnor \rho_\e (y) f^2(y)\dP(y).
	\end{align*}
	Therefore, we can write
	$\langle (\bfPi_\e \bfL^\unnor_\e + \bfL_\e^\unnor\bfPi_\e) f, f \rangle_\LP = E_1 + E_2,$
	where
	\begin{align*}
		E_1 &= \int_\M (\rho_\e (x) -1) \left(\int_\M \fg_\e (x,y) \frac{(f(x) - f(y))^2}{\e} \dP(y) \right)\dP(x), \\
		E_2 &= \int_\M \bfL_\e^\unnor \rho_\e (y) f^2(y)\dP(y).
	\end{align*}
	
	Using Lemma \ref{smooth:KDE}, $E_1$ can be bounded as
	\begin{align*}
		|E_1|
		&= \left| \int_\M (\rho_\e (x) -1) \left(\int_\M \fg_\e (x,y) \frac{(f(x) - f(y))^2}{\e} \dP(y) \right)\dP(x)  \right|
		\\& \leq \sup_{x \in \M} |\rho_\e(x) -1 | \int_\M \int_\M \fg_\e (x,y) \frac{(f(x) - f(y))^2}{\e} \dP(y) \dP(x)
		\\&= \|\bfPi_\e\|_\infty \langle \bfL^\unnor_\e f, f \rangle_\LP
		\leq C_2 \e \langle \bfL^\unnor_\e f, f \rangle_\LP.
	\end{align*}
	Similarly, $E_2$ is bounded as
	\begin{align*}
		|E_2|
		\leq \|\bfL_\e^\unnor \rho_\e\|_\infty \|f\|_\LP^2
		\leq C_2 \e \|f\|_\LP^2.
	\end{align*}
	Therefore, by combining the bounds for $E_1$ and $E_2$, we obtain
	\begin{equation*}
		\left| \langle (\bfPi_\e \bfL^\unnor_\e + \bfL_\e^\unnor\bfPi_\e) f, f \rangle_\LP \right|
		\leq C_2 \e \left( \langle \bfL^\unnor_\e f, f \rangle_\LP + \|f\|_\LP^2\right).
	\end{equation*}
	Now, plugging the above inequality into \eqref{smooth:rel2-eq1} yields
	\begin{align*}
		\left|\langle (\bfL_\e - \bfL^\unnor_\e) f, f \rangle_\LP \right|
		&\leq \frac{C_2 \e }{2} \left( \langle \bfL^\unnor_\e f, f \rangle_\LP + \|f\|_\LP^2\right)
		+ \frac{17C_1 C_2^2 \e \|f\|^2_\LP}{4}
		\\&\leq \left(\frac{C_2}{2}+ \frac{17C_1 C_2^2}{4} \right) \e \left( \langle \bfL^\unnor_\e f, f \rangle_\LP + \|f\|_\LP^2\right).
	\end{align*}	
	This completes the proof. 
\end{proof}

\subsubsection*{Relative form boundedness of $\bfL_\e$ to $\bfL_\e^\heat$}

Finally, we establish the relative form boundedness of $\bfL_\e$ to $\bfL_\e^\heat$.

\begin{proposition} \label{smooth:rel}
	There exist $\M$-dependent constants $r$ and $C$ such that $\e < r$ implies
	\begin{equation*}
		\big| \langle (\bfL_\e- \bfL_\e^\heat) f, f \rangle_\LP \big|
		\leq C \e \log^2(1/\e) \Big( \langle \bfL^\heat_\e f, f \rangle_\LP + \|f\|^2_\LP \Big)
		\quad\text{for all $f \in \LP$.}
	\end{equation*}
\end{proposition}

\begin{proof}[Proof of Proposition~\ref{smooth:rel}]
	For $f \in \LP$, we separate and bound
	\begin{equation}
		| \langle (\bfL_\e- \bfL_\e^\heat) f, f \rangle | 
		\leq | \langle (\bfL_\e^\unnor- \bfL_\e^\heat) f, f \rangle |
		+ | \langle (\bfL_\e- \bfL_\e^\unnor) f, f \rangle |.
		\label{smooth:rel-eq1}
	\end{equation}
	Let $r_1$, $C_1$ and $r_2$, $C_2$ be as in Lemmas \ref{smooth:rel1} and \ref{smooth:rel2}, respectively.
	For $\e < r = \min\{r_1,r_2\}$, we have
	\begin{equation}
		| \langle (\bfL_\e^\unnor - \bfL_\e^\heat) f, f \rangle |
		\leq C_1 \e \log^2 (1/\e) \Big( \langle \bfL^\heat_\e f, f \rangle_\LP + \|f\|^2_\LP \Big).
		\label{smooth:rel-eq2}
	\end{equation}
	and
	\begin{align}
		| \langle (\bfL_\e - \bfL_\e^\unnor) f, f \rangle |
		&\leq C_2 \e \Big( \langle \bfL^\unnor_\e f, f \rangle_\LP + \|f\|^2_\LP \Big)
		\notag\\&= C_2 \e \Big( \langle (\bfL^\unnor_\e-\bfL^\heat_\e) f, f \rangle_\LP + \langle \bfL^\heat_\e f, f \rangle_\LP + \|f\|^2_\LP \Big).
		\label{smooth:rel-eq3}
	\end{align}
	Using \eqref{smooth:rel-eq1} and $\e <r$, we further bound the last term in \eqref{smooth:rel-eq3} as
	\begin{align}
		&\quad C_2 \e \Big( \langle (\bfL^\unnor_\e-\bfL^\heat_\e) f, f \rangle_\LP + \langle \bfL^\heat_\e f, f \rangle_\LP + \|f\|^2_\LP \Big)
		\notag\\&\leq C_2 \e \Big( C_1 \e \log^2(1/\e) \Big( \langle \bfL^\heat_\e f, f \rangle_\LP + \|f\|^2_\LP \Big) 
		+ \langle \bfL^\heat_\e f, f \rangle_\LP + \|f\|^2_\LP \Big)
		\notag\\&= ( C_1 C_2 \e + C_2 \log^{-2} (1/\e) ) \e \log^2(1/\e) \Big( \langle \bfL^\heat_\e f, f \rangle_\LP + \|f\|^2_\LP \Big)
		\notag\\&\leq ( C_1 C_2 r + C_2 \log^{-2} (1/r) ) \e \log^2(1/\e) \Big( \langle \bfL^\heat_\e f, f \rangle_\LP + \|f\|^2_\LP \Big)
		\label{smooth:rel-eq4}
	\end{align}
    Combining \eqref{smooth:rel-eq1}-\eqref{smooth:rel-eq4}  
    completes the proof.
\end{proof}

\subsection{Eigenvalues of smoothed graph Laplacian}

Based on the relative form boundedness of $\bfL_\e$ to $\bfL^\heat$ stated in Proposition~\ref{smooth:rel},
we study the convergence of the eigenvalues $\{\lambda_{k,\e}\}_{k\in\cJ_j}$ of $\bfL_\e$ to the corresponding eigenvalue $\lambda_j$ of $\Delta_\M$.

We begin by stating and proving some supporting lemmas.

\begin{lemma} \label{smooth:eigval:aux1}
	The following hold:
	\begin{enumerate}[(a)]
		\item There exists an $\M$-dependent constant $C$ such that for all $k \in \N$,
		\begin{equation*}
			\|\phi_k\|_\infty \le C \max\big\{1,\lambda_k^{(d-1)/4}\big\}.
		\end{equation*}
		
		\item There exist ($\M$,$\lambda_j$)-dependent constants $r$ and $C$ such that $\e < r$ implies
		\begin{equation*}
			\max_{k \in \cJ_j} \| \bfL_\e \phi_k - \lambda_j \phi_k \|_\infty \leq C \e.
		\end{equation*}
	\end{enumerate}
\end{lemma}

\begin{proof}[Proof of Lemma~\ref{smooth:eigval:aux1}]
	Part (a) is proved in \cite{donnelly2001bounds}, and part (b) is a direct consequence of Proposition~10 and its proof in \cite{coifman2006diffusion}.
\end{proof}

\begin{lemma} \label{smooth:eigval:aux2}
	If $ \e\lambda_{j+} < 1$, then we have
	\begin{equation*}
		\bigg| \frac{1-e^{-\e\lambda_j}}{\e} - \lambda_j \bigg| \le \lambda_j^2 \e, \quad 
		\min_{k \notin \cJ_j} \bigg| \frac{1-e^{-\e\lambda_j}}{\e} - \frac{1-e^{\e\lambda_k}}{\e}\bigg|
		\geq \frac{\gamma_j}{3}.
	\end{equation*}
\end{lemma}

\begin{proof}[Proof of Lemma~\ref{smooth:eigval:aux2}]
	A straightforward calculation gives
	\begin{equation*}
		\left| \frac{1-e^{-\e\lambda_j}}{\e} - \lambda_j \right| 
		= \lambda_j \left| \frac{1-e^{-\e\lambda_j}}{\e \lambda_j} -1 \right|
		\leq \lambda_j \cdot \e \lambda_j = \e \lambda_j^2,
	\end{equation*}
	provided $\e \lambda_j \leq 1 $, which yields the first inequality. For the second inequality, note that
	\begin{equation*}
		\min_{k \not\in\cJ_j} \left| \frac{1-e^{-\e \lambda_j}}{\e} - \frac{1-e^{-\e \lambda_k}}{\e}\right|
		= \frac{1}{\e} \min_{k \not\in\cJ_j} |e^{-\e\lambda_j}-e^{-\e\lambda_k}|
		= \frac{1}{\e}\min \big\{ e^{-\e\lambda_j} - e^{-\e \lambda_{j_+}},	e^{-\e \lambda_{j_-}} - e^{-\e\lambda_j} \big\}.
	\end{equation*}
	Using the mean value theorem, we have
		$e^{-\e\lambda_j} - e^{-\e\lambda_{j_+}}
		= (-\e \lambda_j + \e \lambda_{j_+}) e^{-c}
		= \e (\lambda_{j_+} - \lambda_j) e^{-c}$  
	for some $c \in (\e \lambda_j, \e \lambda_{j_+})$.
	Hence, $\e \lambda_{j_+} < 1$ implies
	$ \e (\lambda_{j_+} - \lambda_{j}) e^{-1} \leq e^{-\e\lambda_j} - e^{-\e\lambda_{j_+}}. $
	Similarly, $\e \lambda_{j_-} < 1$ implies $\e (\lambda_{j} - \lambda_{j_-}) e^{-1} \leq e^{-\e\lambda_{j_-}} - e^{-\e\lambda_j}.$
	Noting $\min_{k \not\in \cJ_j} |\lambda_k - \lambda_j| = \min \{\lambda_{j+}-\lambda_j, \lambda_j - \lambda_{j-}\},$
	combining the above results completes the proof.
\end{proof}

\begin{lemma} \label{smooth:eigval:aux3}
	There exist $(\M,\lambda_j)$-dependent constants $r$ and $C$ such that 
	$\e < r $ implies
	\begin{equation*}
		\max_{k \in \cJ_j} \|(\bfL_\e -\bfL_\e^\heat) \phi_k \|_\infty \leq C \e.
	\end{equation*}
\end{lemma}

\begin{proof}[Proof of Lemma~\ref{smooth:eigval:aux3}]
	Let $k \in \cJ_j$. Writing
	\begin{align*}
		\bfL_\e^\heat \phi_k 
		= \frac{1-e^{-\e \lambda_j}}{\e} \phi_k
		= \lambda_j \phi_k + \bigg(\frac{1-e^{-\e \lambda_j}}{\e} -\lambda_j \bigg)\phi_k,
	\end{align*}
	we obtain
	\begin{equation*}
		\|(\bfL_\e -\bfL_\e^\heat) \phi_k\|_\infty 
		\leq \|\bfL_\e \phi_k - \lambda_j \phi_k\|_\infty +  \bigg|\frac{1-e^{-\e \lambda_j}}{\e} -\lambda_j \bigg| \cdot \|\phi_k\|_\infty.
	\end{equation*}
	Applying Lemmas~\ref{smooth:eigval:aux1} and \ref{smooth:eigval:aux2} completes the proof.
\end{proof}

We now state the convergence of the eigenvalues.

\begin{theorem} \label{smooth:eigval}
	There exist $(\M,\lambda_j)$-dependent constants $r$ and $C$ such that 
	$\e < r $ implies
	\begin{equation*}
		\sqrt{ \sum_{k \in \cJ_j} | \lambda_{k,\e} - \lambda_j |^2}	\leq C \e.
	\end{equation*}
\end{theorem}

\begin{proof}[Proof of Theorem~\ref{smooth:eigval}]
	Let $r$ and $C$ be as in Propositions~\ref{smooth:rel} and Lemma~\ref{smooth:eigval:aux3}, and let $\e <r$.
	By Propositions~\ref{smooth:rel} and \ref{gauss:spec:eigval1}, we have
	\begin{equation}
		\sqrt{ \sum_{k \in \cJ_j} | \lambda_{k,\e} - \lambda_j |^2}
		\leq \frac{43}{18} \sqrt{\sum_{k \in \cJ_j} \| (\bfL_\e -\bfL_\e^\heat) \phi_k \|_\LP^2 } 
		\label{smooth:eigval-eq1}
	\end{equation}
	provided 
	\begin{equation}
		C \e \log^2(1/\e) \cdot \bigg( \min_{k \notin \cJ_j} 
		\Big| \frac{1-e^{-\e\lambda_j}}{\e} - \frac{1-e^{-\e\lambda_k}}{\e} \Big| \bigg)^{-1} \bigg(2 \cdot \frac{1-e^{-\e\lambda_k}}{\e} + 1\bigg) \le 1/4.
		\label{smooth:eigval-eq2}
	\end{equation}
	By using Lemma~\ref{smooth:eigval:aux2} and adjusting $r$ if necessary, we ensure that \eqref{smooth:eigval-eq2} holds whenever $\e <r$.
	Moreover, applying Lemma~\ref{smooth:eigval:aux3}, we further bound \eqref{smooth:eigval-eq1} as
	\begin{align*}
		\sqrt{ \sum_{k \in \cJ_j} | \lambda_{k,\e} - \lambda_j |^2}
		\leq \frac{43}{18} |\cJ_j| C \e. 
	\end{align*}
	By adjusting $C$ appropriately, we complete the proof.
\end{proof}

\subsection{Eigenfunctions of smoothed graph Laplacians}

In this subsection, we study the convergence of the eigenfunctions $\{\phi_{k,\e}\}_{k \in \cJ_j}$ of $\bfL_\e$ to the corresponding eigenfunctions $\{\phi_k\}_{k \in \cJ_j}$ of $\Delta_\M$.
The key point is that $\phi_{k,\e}$ need not be close to $\phi_k$ for the same index $k$;
rather, each $\phi_{k,\e}$ is expected to be close to an appropriate linear combination of $\{ \phi_\ell\}_{\ell \in \cJ_j}$.

Specifically, we define
\begin{equation*}
	\xi_{k,\e} 
	= \phi_{k,\e} - \sum_{\ell \in \cJ_j} v_{k\ell,\e} \phi_\ell,
	\quad\text{where $v_{k\ell,\e} = \langle \phi_{k,\e}, \phi_\ell \rangle_\mbH$},
\end{equation*}
and investigate its convergence to zero as $\e$ tends to zero.
We also introduce the matrix $V_{j,\e}$ consisting of the coefficients $v_{k\ell,\e}$, defined by
\begin{equation*}
	V_{j,\e} = (v_{k\ell,\e})_{k,\ell \in \cJ_j} = \big( \langle \phi_{k,\e}, \phi_\ell \rangle_\LP \big)_{k,\ell \in \cJ_j},
\end{equation*}
and investigate how $V_{j,\e}$ behaves.
We summarize the results in the following theorem.

\begin{theorem} \label{smooth:eigfun1}
	There exist $(\M,\lambda_j)$-dependent constants $r$ and $C$ such that $\e < r$ implies
	\begin{gather}
		\sqrt{\sum_{k \in \cJ} \|\xi_{k,\e}\|_\LP^2} \leq C \e, \quad
		\| V_{j,\e} V_{j,\e}^\top - I_{|\cJ_j|} \|_\F \leq C \e^2.
		\label{smooth:eigfun1-res1}
	\end{gather}
	In particular, $V_{j,\e}$ satisfies
	\begin{equation}
		\max_{k,\ell \in \cJ_j} |v_{k\ell,\e}|,\, 
		\|V_{j,\e}\|_\F \leq C, 
		\label{smooth:eigfun1-res2}
	\end{equation}
	and for any matrix $A \in \R^{|\cJ_j|\times|\cJ_j|}$,
	\begin{equation} 
		\Big| \tr \big( V_{j,\e} A V_{j,\e}^\top \big) - \tr(A) \Big| \leq C\e^2 \|A\|_\F.
		\label{smooth:eigfun1-res3}
	\end{equation}
\end{theorem}

\begin{proof}[Proof of Theorem~\ref{smooth:eigfun1}]
	Let $r$ and $C$ be as in Propositions~\ref{smooth:rel} and Lemma~\ref{smooth:eigval:aux3}, and let $\e <r$.
	By Propositions~\ref{smooth:rel}, \ref{gauss:spec:eigvec1} and \ref{gauss:spec:eigvec2}, we have
	\begin{gather*}
		\sqrt{\sum_{k\in\cJ_j}\|\xi_{k,\e}\|_\LP^2} 
		\leq \frac{36+7\sqrt{2}}{18} \bigg( \min_{k \notin \cJ_j} 
		\Big| \frac{1-e^{-\e\lambda_j}}{\e} - \frac{1-e^{-\e\lambda_k}}{\e} \Big| \bigg)^{-1} \sqrt{\sum_{k \in \cJ} \| (\bfL_\e - \bfL^\heat) \psi_k \|_\LP^2 }, \\
		\| \hV_{j,\e}^\top \hV_{j,\e} - I_{|\cJ_j|} \|_\F \leq 4 \bigg( \min_{k \notin \cJ_j} 
		\Big| \frac{1-e^{-\e\lambda_j}}{\e} - \frac{1-e^{-\e\lambda_k}}{\e} \Big| \bigg)^{-2} \sum_{k \in \cJ} \|(\bfL_\e-\bfL_\e^\heat)\psi_k\|^2_\mbH,
	\end{gather*}
	provided
	\begin{equation*}
		C \e \log^2(1/\e) \cdot \bigg( \min_{k \notin \cJ_j} 
		\Big| \frac{1-e^{-\e\lambda_j}}{\e} - \frac{1-e^{-\e\lambda_k}}{\e} \Big| \bigg)^{-1} \bigg(2 \cdot \frac{1-e^{-\e\lambda_k}}{\e} + 1\bigg) \le 1/4.		
	\end{equation*}
	Hence, \eqref{smooth:eigfun1-res1} follows by Lemmas~\ref{smooth:eigval:aux2} and \ref{smooth:eigval:aux3}, and adjusting $r$ if necessary.
	
	\eqref{smooth:eigfun1-res2} and \eqref{smooth:eigfun1-res3} can be shown using \eqref{smooth:eigfun1-res1} and 
	\begin{gather*}
		\max_{k,\ell \in \cJ_j} |v_{k\ell,\e}|^2 \leq \max_{k \in \cJ_j} \sum_{\ell \in \cJ_j} |v_{k\ell,\e}|^2
		\leq \|V_{j,\e} V_{j,\e}^\top\|_\op \leq 1+ \|V_{j,\e} V_{j,\e}^\top-I_{|\cJ_j|}\|_\F, 
		\\
		\|V_{j,\e}\|_\F^2 = \tr(V_{j,\e} V_{j,\e}^\top-I_{|\cJ_j|}) + |\cJ_j| 
		\leq \sqrt{|\cJ_j|} \cdot \|V_{j,\e} V_{j,\e}^\top-I_{|\cJ_j|} \|_\F + |\cJ_j|,
		\\
		|\tr(V_{j,\e} A V_{j,\e}^\top) - \tr(A)| 
		\leq |\tr(A V_{j,\e}^\top V_{j,\e} - I_{|\cJ|} )| 
		\leq \|A\|_\F \|V_{j,\e}^\top V_{j,\e} - I_{|\cJ|}\|_\F
		= \|A\|_\F \|V_{j,\e}V_{j,\e}^\top - I_{|\cJ|}\|_\F.
	\end{gather*}
	This completes the proof.
\end{proof}

The above theorem provides the $\LP$-convergence of the eigenfunctions.
This result can be extended to uniform and Lipschitz-type convergence.
Such extensions were studied in \cite{calder2022lipschitz} for unnormalized graph Laplacians with finitely supported kernels.
We adapt these results to our setting, as summarized below.
The proof is postponed to Appendix~\ref{smooth:eigfun2-proof}.

\begin{theorem} \label{smooth:eigfun2}
	There exist $(\M,\lambda_j)$-dependent constants $r$ and $C$ such that $\e < r$ implies
	\begin{equation}
		\max_{k \in \cJ_j} \left\| \xi_{k,\e} \right\|_\infty, \,\,
		\max_{k \in \cJ_j} \left\| \bfL_\e \xi_{k,\e} \right\|_\infty, \,\,
		\max_{k \in \cJ_j} \sup_{x,y \in \M} \frac{|\xi_{k,\e}(x) - \xi_{k,\e} (y)|}{\dist_\M (x,y)+\e^{1/3} \log(1/\e) } \le C \e,
		\label{smooth:eigfun2-res1}
	\end{equation}
	and
	\begin{equation}
		\max_{k \in \cJ} \|\phi_{k,\e}\|_\infty \leq C.
		\label{smooth:eigfun2-res2}
	\end{equation}
\end{theorem}

\subsection{Additional results}

Before finishing this section, we state several results that follow from the analysis above and will be used to prove our main results.

\begin{lemma} \label{smooth:extra1}
	There exist $(\M,\lambda_j)$-dependent constants $r$ and $C$ such that $\e < r$ implies
	\begin{align*}
		\min_{k \in \cJ, \ell \notin \cJ} |\lambda_{k,\e} - \lambda_{\ell,\e}| \ge \frac{\gamma_j}{2}.
	\end{align*}
\end{lemma}

\begin{proof}[Proof of Lemma~\ref{smooth:extra1}]
	The claim follows from Propositions~\ref{smooth:rel} and \ref{gauss:spec:eigval2}, and Lemma~\ref{smooth:eigval:aux2}, after adjusting $r$ if necessary.
\end{proof}

The lemma below will be used to prove Theorem~\ref{gauss:eigval}.
We defer the proof of this lemma to Appendix~\ref{smooth:extra2-proof}.

\begin{lemma} \label{smooth:extra2}
	There exist $(\M,\lambda_j)$-dependent constants $r$ and $C$ such that if $\e < r$,
	we have with probability at least $1-\tau$,
	\begin{equation*}
		\Bigg| (\P_n-\P) \Bigg( \sum_{k \in \cJ_j} \mu_{k,\e} \bfL_\e^\unnor \bigg( \frac{\phi^2_{k,\e} }{\rho_\e} \Bigg)
		- \sum_{k \in \cJ_j} \Delta_\M (\phi^2_k) \Bigg) \Bigg|
		\leq C \e^{5/6} \log(1/\e) \sqrt{\frac{2\log(2/\tau)}{n}},
	\end{equation*}
	where $\mu_{k,\e} = 1 - \e \lambda_{k,\e}$.
\end{lemma}

The following two lemmas will be used to prove Theorem~\ref{gauss:eigvec}.
Roughly speaking, these results imply $\{\phi_{k,\e}\}_{k \in \cJ_j}$ can be replaced by $\{\phi_k\}_{k \in \cJ_j}$ within an acceptable bound.

\begin{lemma} \label{smooth:extra3}
	There exists $(\M,\lambda_j)$-dependent constants $r$ and $C$ such that if $\e< r$,
	we have with probability at least $1-\tau$,
	\begin{equation*}
		\Bigg| (\P_n-\P) \Bigg(\sum_{k\in\cJ_j} \phi_{k,\e}^2 - \sum_{k\in\cJ_j} \phi_k^2 \Bigg) \Bigg| 
		\leq C \e \sqrt{\frac{2\log(2/\tau)}{n}}.
	\end{equation*}
\end{lemma}

\begin{proof}[Proof of Lemma~\ref{smooth:extra3}]
	Once we show 
	\begin{equation}
		\bigg\| \sum_{k\in\cJ_j} \phi_{k,\e}^2 - \sum_{k\in\cJ_j} \phi_k^2 \bigg\|_\infty \leq C \e,
		\label{smooth:extra3-eq1}
	\end{equation}
	the claim follows directly from Hoeffding’s inequality. By letting $B = (\phi_{k,\e})_{k \in \cJ_j}$, $A = (\phi_k)_{k \in \cJ_j}$, we see that
	\begin{equation*}
		\sum_{k\in\cJ_j} \phi_{k,\e}^2 - \sum_{k\in\cJ_j} \phi_k^2 = \tr( B B^\top - A A^\top).
	\end{equation*}
	On the other hand, by letting $E = (\xi_{k,\e})_{k \in \cJ_j}$ and using $\phi_{k,\e} = \sum_{\ell \in \cJ_j} v_{k\ell,\e} \phi_\ell + \xi_{k,\e}$, we have
	$ B = V_{j,\e} A + E$, and this gives
	\begin{equation*}
		B B^\top - A A^\top = V_{j,\e} A A^\top V_{j,\e}^\top - A A^\top + V_{j,\e} A E^\top + E A^\top V_{j,\e}^\top + E E^\top,
	\end{equation*}
	and consequently,
	\begin{equation*}
		|\tr(B B^\top - A A^\top)|
		\leq |\tr(V_{j,\e} A A^\top V_{j,\e}^\top) - \tr(A A^\top)| 
		+ 2 |\tr(V_{j,\e} A E^\top)| + |\tr(E E^\top)|.
	\end{equation*}
	To bound the right-hand side above, we see that
	\begin{gather*}
		|\tr(V_{j,\e} A E^\top)| \leq \|V_{j,\e}\|_\F \|A E^\top\|_\F = \|V_{j,\e}\|_\F \|A\|_2 \|E\|_2, \quad
		\tr(E E^\top) = \|E\|_2^2,
	\end{gather*}
	and using Theorem~\ref{smooth:eigfun1},
	\begin{equation*}
		\big| \tr(V_{j,\e} A A^\top V_{j,\e}^\top) - \tr(A A^\top) \big|
		\leq C_1 \e^2 \| A A^\top \|_\F 
		= C_1 \e^2 \|A\|_2^2,
	\end{equation*}
	provided $\e < r_1$ for some $(\M,\lambda_j)$-dependent constant $C_1$ and $r_1$.
	Combining the above, we obtain
	\begin{equation}
		\bigg| \sum_{k\in\cJ_j} \phi_{k,\e}^2 - \sum_{k\in\cJ_j} \phi_k^2 \bigg| 
		\leq C_1 \e^2 \|A\|_2^2 + 2 \|V_{j,\e}\|_\F \|A\|_2 \|E\|_2 + \|E\|_2^2
		\label{smooth:extra3-eq3}
	\end{equation}
	
	Since 
	\begin{equation*}
		\|A\|_2 = \sqrt{\sum_{k \in \cJ_j} \phi_k^2} \leq \sqrt{|\cJ_j|} \max_{k \in \cJ_j} \|\phi_k\|_\infty, \quad\text{and}\quad
		\|E\|_2 = \sqrt{\sum_{k \in \cJ_j} \xi_{k,\e}^2} \leq \sqrt{|\cJ_j|} \max_{k \in \cJ_j} \|\xi_{k,\e}\|_\infty,
	\end{equation*}	
	we see that \eqref{smooth:extra3-eq1} follows by applying Lemma~\ref{smooth:eigval:aux1} and Theorems~\ref{smooth:eigfun1} and \ref{smooth:eigfun2} to \eqref{smooth:extra3-eq3}.
	This completes the proof.
\end{proof}

\begin{lemma} \label{smooth:extra4}
	There exists $(\M,\lambda_j)$-dependent constants $r$ and $C$ such that if $\e< r$,
	\begin{equation*}
		\bigg| \sum_{k \in \cJ_j} \sum_{\ell \in \cJ_j} \langle \hvarphi_{k,n,\e}, \sn{\phi_{\ell,\e}} \rangle^2_\LPn 
		- \sum_{k \in \cJ_j} \sum_{\ell \in \cJ_j} \langle \hvarphi_{k,n,\e}, \sn{\phi_\ell} \rangle^2_\LPn \bigg|
		\leq C\e.
	\end{equation*}
\end{lemma}

\begin{proof}[Proof of Lemma~\ref{smooth:extra4}]
	By letting $B = (\langle \hvarphi_{k,n,\e}, \sn{\phi_{\ell,\e}} \rangle_\LPn)_{k,\ell \in \cJ_j}$ and $A = (\langle \hvarphi_{k,n,\e}, \sn{\phi_\ell} \rangle_\LPn)_{k,\ell \in \cJ_j}$, we see that
	\begin{equation}
		\sum_{k \in \cJ_j} \sum_{\ell \in \cJ_j} \langle \hvarphi_{k,n,\e}, \sn{\phi_{\ell,\e}} \rangle^2_\LPn 
		- \sum_{k \in \cJ_j} \sum_{\ell \in \cJ_j} \langle \hvarphi_{k,n,\e}, \sn{\phi_\ell} \rangle^2_\LPn
		= \tr(B^\top B -A^\top A).
		\label{smooth:extra4-eq1}
	\end{equation}
	On the other hand, by letting $E = (\langle \hvarphi_{k,n,\e}, \sn{\xi_{\ell,\e}} \rangle_\LPn )_{k,\ell \in \cJ_j}$ and using $\phi_{k,\e} = \sum_{\ell \in \cJ_j} v_{k\ell,\e} \phi_\ell + \xi_{k,\e}$, we have
	$ B = A V^\top_{j,\e} + E$, and this gives
	\begin{equation*}
		B^\top B - A^\top A = V_{j,\e} A^\top A V^\top_{j,\e}  - A^\top A + V_{j,\e} A^\top E + E^\top A V_{j,\e}^T + E^\top E,
	\end{equation*}
	and consequently,
	\begin{align}
		|\tr(B^\top B - A^\top A)|
		&\leq |\tr(V_{j,\e} A^\top A V^\top_{j,\e}) - \tr(A^\top A)|
		+ 2 |\tr(V_{j,\e} A^\top E)| + |\tr(E^\top E)|.
		\notag\\&\leq \|A^T A\|_\F \|V_{j,\e}^\top V_{j,\e} - I_{|\cJ_j|}\|_\F + 2 \|V_{j,\e}\|_\F \|A\|_\F \|E\|_\F +\|E\|_\F^2.
		\label{smooth:extra4-eq2}
	\end{align}
	Noting $\{\hvarphi_{k,n,\e}\}_{k \in \cJ_j}$ are $\LPn$-orthonormal, Bessel's inequality yields
	\begin{equation*}
		\|A\|_\F^2 = \sum_{\ell \in \cJ_j} \sum_{k \in \cJ_j} \langle \hvarphi_{k,n,\e}, \sn{\phi_\ell} \rangle_\LPn^2
		\leq \sum_{\ell \in \cJ_j} \| \sn{\phi_\ell} \|_\LPn^2
		\leq |\cJ_j| \cdot \max_{\ell \in \cJ_j} \|\phi_\ell\|_\infty^2,
	\end{equation*}
	and
	\begin{equation*}
		\|E\|_\F^2 = \tr(E^\top E) = \sum_{k \in \cJ_j} \sum_{\ell \in \cJ_j} \langle \hvarphi_{k,n,\e}, \sn{\xi_{\ell,\e}} \rangle_\LPn^2
		\leq \sum_{\ell \in \cJ_j} \| \sn{\xi_{\ell,\e}} \|_\LPn^2 
		\leq |\cJ_j| \cdot \max_{\ell \in \cJ_j} \|\xi_{\ell,\e}\|_\infty^2.
	\end{equation*}
	Combining \eqref{smooth:extra4-eq1} and \eqref{smooth:extra4-eq2}, and applying Lemma~\ref{smooth:eigval:aux1} and Theorems~\ref{smooth:eigfun1} and \ref{smooth:eigfun2}, we completes the proof.
\end{proof}

\section{First-order asymptotic expansion for spectral convergence of population-normalized graph Laplacians}
\label{popul} 

As mentioned in Section~\ref{gauss:outline}, we consider the population-normalized graph Laplacian $\bL_{n,\e} \in \R^{n \times n}$ defined in \eqref{gauss:outline-eq2}. 
Let $\{ \blambda_{k,n,\e}\}_{k \in [n]}$ denote the eigenvalues of $\bL_{n,\e}$, arranged in non-decreasing order, with corresponding $\LPn$-orthonormal eigenvectors $\{\bvarphi_{k,n,\e}\}_{k \in [n]}.$

The objective of this section is to establish the first-order asymptotic expansions of the eigenvalues and eigenvectors of $\bL_{n,\e}$ to those of $\bfL_\e$.
The key idea is to analyze the spectral convergence of $\bW_{n,\e}$ to $\bfW_\e$.
There are two main reasons for this approach.

First, $\bW_{n,\e}$ and $\bfW_\e$ share the same spectral structure as $\bL_{n,\e}$ and $\bfL_\e$, respectively.
Specifically, for $k =1, \dots, n$, define
\begin{equation*}
	\bmu_{k,n,\e} = 1- \e \blambda_{k,n,\e}, \quad\text{equivalently,}\quad \blambda_{k,n,\e} = \frac{1-\bmu_{k,n,\e}}{\e}.
\end{equation*}
From $\bW_{n,\e} = I_n - \e \bL_{n,\e}$, we see that $\bW_{n,\e} \bvarphi_{k,n,\e} = \bmu_{k,n,\e} \bvarphi_{k,n,\e}.$
Thus, $\bmu_{k,n,\e}$ is the eigenvalue of $\bW_{n,\e}$ corresponding to the eigenvector $\bvarphi_{k,n,\e}$.
Similarly, we see that $\mu_{k,\e}$, defined as
\begin{equation*}
	\mu_{k,\e} = 1- \e \lambda_{k,\e}, \quad\text{equivalently,}\quad \lambda_{k,\e} = \frac{1-\mu_{k,\e}}{\e},
\end{equation*}
is the eigenvalue of $\bfW_\e$ corresponding to the eigenfunction $\phi_{k,\e}$.
Therefore, the spectral convergence of $\bW_{n,\e}$ and $\bfW_\e$ directly implies that of $\bL_{n,\e}$ and $\bfL_\e$.

Second, the spectral convergence $\bW_{n,\e}$ is more conveniently analyzed within the framework of kernel Gram matrices.
We can view $\bW_{n,\e}$ as the kernel Gram matrix associated with $\fw_\e$, and $\bfW_\e$ is the corresponding kernel integral operator. In particular, $\fw_\e$ satisfies Mercer's conditions, and this enables us to analyze the spectral quantities of $\bW_{n,\e}$ and $\bfW_\e$ via the reproducing kernel Hilbert space associated with $\fw_\e$.

The spectral convergence of kernel Gram matrices to their associated kernel integral operators is a well-studied topic in statistics and machine learning.
In particular, we rely on \cite{bae2026asymptotic}, which provides first-order asymptotic expansions that are suitable for our purpose.
We summarize the relevant results and present their adapted forms in Appendix~\ref{gauss:mercer}.

\subsection{Eigenvalues of population-normalized graph Laplacians}

To simplify the proofs, in this and the next section, generic constants (e.g., $r_k$, $C_k$) are assumed to be $(\M,\lambda_j)$-dependent unless specified otherwise.

Here we derive the first-order asymptotic expansion for the sum of the differences of eigenvalues corresponding to $\cJ_j$: 
\begin{equation}
	\sum_{k \in \cJ_j} ( \blambda_{k,n,\e} - \lambda_{k,\e} ).
	\label{gauss:popul:eigen-goal}
\end{equation}

To this end, we introduce the following quantities:
\begin{equation*}
	\varrho_\e = \sup_{x \in \M} \fw_\e (x,x), \quad
	\gamma_{j,\e} = \min_{k \in \cJ_j, \ell \notin \cJ_j} |\lambda_{k,\e} - \lambda_{\ell,\e}|, \quad
	\lambda_{\tmax,j,\e} = \max_{k \in \cJ_j} \lambda_{k,\e}, \quad
	\lambda_{\tmin,j,\e} = \min_{k \in \cJ_j} \lambda_{k,\e}.
\end{equation*}
By Theorem~\ref{smooth:eigval} and Lemma~\ref{smooth:extra1}, we can control $\gamma_{j,\e}$, $\lambda_{\tmax,j,\e}$, and $\lambda_{\tmin,j,\e}$. 
For $\varrho_\e$, we formulate a lemma below.

\begin{lemma} \label{popul:eigval:aux}
	There exist $\M$-dependent constant $r$, $C_1$ and $C_2$ such that $\e < r$ implies
	\begin{equation*}
		\frac{C_1}{\e^{d/2}} \leq \varrho_\e \leq \frac{C_2}{\e^{d/2}}.
	\end{equation*}
\end{lemma}

\begin{proof}[Proof of Lemma~\ref{popul:eigval:aux}]
	By definition, we see that
	\begin{equation*}
		\fw_\e (x,x) = \frac{\fg_\e (x,x)}{\rho_\e (x)} = \frac{1}{(4\pi\e)^{d/2} \rho_\e (x)}.
	\end{equation*}
	By Lemma~\ref{smooth:KDE}, $\rho_\e$ satisfies
	$\|\rho_\e -1\|_\infty \leq C \e$, provided $\e <r$ for some $\M$-dependent constant $r$ and $C$.
	Therefore, the claim follows by adjusting $r$ appropriately.
\end{proof}

Now, we present the first-order asymptotic expansion of \eqref{gauss:popul:eigen-goal}.

\begin{theorem} \label{popul:eigval}
	There exist $(\M,\lambda_j)$-dependent constants $r$, $C$, $C'$, $C''$, $C'''$ such that 
	if $\e < r$ and	$ \frac{n\e^{d/2+2}}{\log(1/\e)} > C'''$,
	we have with probability at least $1- \tau - C'\exp ( -C'' n\e^{d/2+2} )$, 
	\begin{equation*}
		\Bigg| \sum_{k \in \cJ_j} (\blambda_{k,n,\e} - \lambda_{k,\e})
		- (\P_n-\P) \Bigg( -\sum_{k \in \cJ_j} \frac{\mu_{k,\e}\phi_{k,\e}^2}{\e} \Bigg) \Bigg|
		\leq C \Bigg( \sqrt{\frac{\log(2/\tau)}{n \e^{d/2 + 2}}} +  \frac{\log(2/\tau)}{n \e^{d/2 + 1}} \Bigg)^2.
	\end{equation*}	
\end{theorem}

\begin{proof}[Proof of Theorem \ref{popul:eigval}]
	Our strategy is to apply spectral convergence results to $\bW_{n,\e}$ and simplify the resulting expressions by replacing them with more tractable forms.
	The same argument applies to other proofs in this section; hence, we refer to this proof for the details and omit the repeated derivations.
	
	Noting $\mu_{k,\e} = 1 - \e \lambda_{k,\e}$, we see that
	\begin{equation*}
		\min_{k\in\cJ_j, \ell\notin\cJ_j} |\mu_{k,\e} - \mu_{\ell,\e}| = \e \gamma_{j,\e}, \quad
		\max_{k \in \cJ_j} \mu_{k,\e} = 1 - \e \lambda_{\tmin,j,\e}, \quad
		\min_{k \in \cJ_j} \mu_{k,\e} = 1 - \e \lambda_{\tmax,j,\e}.
	\end{equation*}
	Equipped with these, applying Proposition~\ref{gauss:mercer:eigval1} to $\bW_{n,\e}$ yields the following:		
	
	If $n$ satisfies
	\begin{equation}
		n \ge \frac{38 \varrho_\e\,(1-\e \lambda_{\tmin,j,\e}) }{\e^2\gamma_{j,\e}^2},
		\label{popul:eigval-eq1}
	\end{equation}
	then with probability at least 
	\begin{equation}
		1-\tau-\frac{2\varrho_\e}{\e\gamma_{j,\e}} \exp \bigg( -\frac{n\e^2\gamma_{j,\e}^2}{70 \varrho_\e\, (1-\e\lambda_{\tmin,j,\e}) } \bigg),
		\label{popul:eigval-eq2}
	\end{equation}
	we have
	\begin{multline}
		\bigg| \sum_{k\in\cJ_j} (\bmu_{k,n,\e} - \mu_{k,\e})
		- (\Pn-\P) \Big( \sum_{k \in \cJ} \mu_{k,\e} \phi_{k,\e}^2 \Big) \bigg| 
		\\ \leq \bigg(\frac{5}{2\e\gamma_{j,\e}} + \frac{4\e \lambda_{\tmax,j,\e} - 4\e \lambda_{\tmin,j,\e})}{\e^2\gamma_{j,\e}^2} \bigg) 
		\bigg( \sqrt{\frac{2\varrho_\e \sum_{k\in\cJ}(1-\e\lambda_{k,\e})\log(2/\tau)}{n}} +  \frac{\varrho_\e \log(2/\tau)}{n} \bigg)^2.
		\label{popul:eigval-eq3}
	\end{multline}
	
	Let $r$ satisfy Theorem~\ref{smooth:eigval} and Lemmas~\ref{smooth:extra1} and \ref{popul:eigval:aux}, so that
	\begin{equation*}
		\frac{C_1}{\e^{d/2}} \leq \varrho_\e \leq \frac{C_2}{\e^{d/2}}, \quad
		\gamma_{j,\e} \ge \frac{\gamma_j}{2}, \quad
		|\lambda_{\tmax,j,\e}-\lambda_j| \le C_3 \e, \quad
		|\lambda_{\tmin,j,\e}-\lambda_j| \le C_3 \e.
	\end{equation*}	
	Now, we see that \eqref{popul:eigval-eq1} is equivalent to $n \ge C'''_1 \e^{-(d/2+2)}$, hence it holds when 
	$ \frac{n\e^{d/2+2}}{\log(1/\e)} > C'''$. 
	For \eqref{popul:eigval-eq2}, we can bound
	\begin{align*}
		\frac{2\varrho_\e}{\e\gamma_{j,\e}} \exp \bigg( -\frac{n\e^2\gamma_{j,\e}^2}{70 \varrho_\e\, (1-\e\lambda_{\tmin,j,\e}) } \bigg)
		&\leq C' \e^{-d/2-1} \exp(- C''_{1} n\e^{d/2+2} )
		\\&=  C' \exp \bigg(- C''_1 \bigg( 1- \frac{(d/2+1)\log(1/\e)}{C''_1 n\e^{d/2+2}} \bigg) n\e^{d/2+2} \bigg)
		\\&\leq C' \exp (- C'' n\e^{d/2+2} ),
	\end{align*}
	provided $C'''$ is sufficiently large. Lastly, using $\bmu_{k,n,\e} = 1 -\e \blambda_{k,n,\e}$, \eqref{popul:eigval-eq3} is equivalent to
	\begin{equation*}
		\bigg| -\e \sum_{k\in\cJ_j} (\blambda_{k,n,\e} - \lambda_{k,\e})
		- (\Pn-\P) \Big( \sum_{k \in \cJ} \mu_{k,\e} \phi_{k,\e}^2 \Big) \bigg| 
		\leq \frac{C}{\e} \bigg( \sqrt{\frac{\log(2/\tau)}{n\e^{d/2}}} +  \frac{\log(2/\tau)}{n\e^{d/2}} \bigg)^2.
	\end{equation*}
	Therefore, dividing both sides by $\e$ completes the proof.
\end{proof}

\subsection{Eigenvectors of population-normalized graph Laplacians}

We present the first-order asymptotic expansion of the eigenvectors of $\bL_{n,\e}$ as follows.

\begin{theorem} \label{popul:eigvec}
	There exist $(\M,\lambda_j)$-dependent constants $r$, $C$, $C'$, $C''$, $C'''$ such that 
	if $\e < r$ and	$ \frac{n\e^{d/2+2}}{\log(1/\e)} > C'''$,
	we have with probability at least $1- \tau - C'\exp \left( -C''n\e^{d/2+2} \right)$,
	\begin{multline*}
		\Bigg| \sum_{k \in \cJ_j} \sum_{\ell \in \cJ_j} \langle \bvarphi_{k,n,\e}, \sn{\phi_{\ell,\e}} \rangle^2_\LPn - |\cJ_j|
		- (\P_n-\P)\bigg(\sum_{k \in \cJ_j}\phi_{k,\e}^2\bigg) \Bigg|
		 \leq C \Bigg( \sqrt{\frac{\log(2/\tau)}{n \e^{d/2 + 2}}} +  \frac{\log(2/\tau)}{n \e^{d/2 + 1}} \Bigg)^2.
	\end{multline*}	
\end{theorem}

\begin{proof}[Proof of Theorem~\ref{popul:eigvec}]
	Applying Proposition~\ref{gauss:mercer:eigvec1} to $\bW_{n,\e}$ yields the following:
	
	If $n \ge \frac{38 \varrho_\e (1-\e\lambda_{\tmin,j,\e}) }{\e^2\gamma_{j,\e}^2}$, we have with probability at least 
	$1-\tau-\frac{2\varrho_\e}{\e\gamma_{j,\e}} \exp \Big( -\frac{n\e^2\gamma_{j,\e}^2}{70 \varrho_\e\,(1-\e\lambda_{\tmin,j\e})} \Big) $,
	\begin{multline*}
		\Bigg| \sum_{k \in \cJ} \sum_{\ell \in \cJ} \langle \bvarphi_{k,n,\e}, \sn{\phi_{\ell,\e}} \rangle_\LP^2 - |\cJ_j| - (\P_n-\P) \Big( \sum_{k \in \cJ} \phi_{k,\e}^2 \Big) \Bigg|
		\\ \leq
		\frac{11}{2\e^2 \gamma_{j,\e}^2} 
		\bigg( \sqrt{\frac{2 \varrho_\e \sum_{k\in\cJ}(1-\e\lambda_{k,\e})\log(2/\tau)}{n}} + \frac{\varrho_\e \log(2/\tau)}{n} \bigg)^2.
	\end{multline*}
	By following the same arguments as those in the proof of Theorem~\ref{popul:eigval}, we complete the proof.
	
\end{proof}

\subsection{Additional results}

We state several results related to the spectral quantities of $\bL_{n,\e}$, which will be used in the following section.

\begin{lemma} \label{popul:extra1}
	There exist $(\M,\lambda_j)$-dependent constants $r$, $C$, $C'$, $C''$, $C'''$ such that if $\e < r$ and $ \frac{n\e^{d/2+2}}{\log(1/\e)} > C'''$,
	we have with probability at least $1-C'\exp \left( -C''n\e^{d/2+2} \right)$,
	\begin{gather}
		\bgamma_{j,n,\e} := \min_{k \in \cJ_j, \ell \not\in \cJ_j} |\blambda_{k,n,\e}-\blambda_{\ell,n,\e}| \ge \frac{\gamma_j}{4}, 
		\label{popul:extra1-res1} \\
		\max_{k \in \cJ_j} |\blambda_{k,n,\e} - \lambda_j| \le 1.
		\label{popul:extra1-res2}
	\end{gather}
\end{lemma}

\begin{proof}[Proof of lemma~\ref{popul:extra1}]
	Applying Proposition~\ref{gauss:mercer:eigval2} to $\bW_{n,\e}$ yields the following:
	
	If $n \ge \frac{38 \varrho_\e \,(1-\e\lambda_{\tmin,j,\e})}{\e \gamma_{j,\e}^2}$, we have with probability at least
	$1-\tau-\frac{2\varrho_\e}{\e\gamma_{j,\e}} \exp \Big( -\frac{n\e^2 \gamma_{j,\e}^2}{70 \varrho_\e \,(1-\e\lambda_{\tmin,j,\e})} \Big) $,
	\begin{gather*}
		\min_{k \in \cJ_j, \ell \not\in \cJ_j} |\e\blambda_{k,n,\e}-\e\blambda_{\ell,n,\e}| \ge \frac{\e \gamma_{j,\e}}{2}, \\
		\max_{k \in \cJ_j} |\e\blambda_{k,n,\e} - \e\lambda_{k,\e}| 
		\leq \sqrt{\frac{4(1-\e\lambda_{\tmin,j,\e}) \log(2\varrho_\e/(\e\gamma_{j,\e}\tau))}{n}} + \frac{2\varrho_\e\log(2\varrho_\e/(\e\gamma_{j,\e} \tau))}{3n}.
	\end{gather*}
	By following the same arguments as those in the proof of Theorem~\ref{popul:eigval}, we obtain \eqref{popul:extra1-res1} and 
	\begin{equation*}
		\max_{k \in \cJ_j} |\blambda_{k,n,\e} - \lambda_{k,\e}| \le C \bigg( \sqrt{\frac{\log(1/(\e\wedge\tau)}{n\e^{d/2+2}}} + \frac{\log(1/(\e\wedge\tau)}{n\e^{d/2+2}} \bigg).
	\end{equation*}
	with probability at least $1-\tau-C'\exp \left( -C''n\e^{d/2+2} \right)$.
	Now, using the inversion of the high probability argument and $\max_{k \in \cJ}|\lambda_{k,\e}-\lambda_j|\leq C_1\e$ from Theorem~\ref{smooth:eigval}, we can obtain \eqref{popul:extra1-res2}. This completes the proof.
\end{proof}

\begin{lemma} \label{popul:extra2}
	There exist $(\M,\lambda_j)$-dependent constants $r$, $C$, $C'$, $C''$, $C'''$ such that if $\e < r$ and $ \frac{n\e^{d/2+2}}{\log(1/\e)} > C'''$,
	we have with probability at least $1-C'\exp \left( -C'' n\e^{d/2+2} \right)$,
	\begin{multline*}
		\bigg| \sum_{k \in \cJ_j} \bmu_{k,n,\e} \langle A \bvarphi_{k,n,\e}, \bvarphi_{k,n,\e} \rangle_\LPn
		- \sum_{k \in \cJ_j} \mu_{k,\e} \langle A \sn{\phi_{k,\e}}, \sn{\phi_{k,\e}} \rangle_\LPn \bigg|
		\\ \leq C \|A\|_{\op,\LPn} \|\bW_{n,\e}\|_{\op,\LPn} \Bigg( \sqrt{\frac{\log(2/\tau)}{n \e^{d/2 + 2}}} +  \frac{\log(2/\tau)}{n \e^{d/2 + 1}} \Bigg),
	\end{multline*}
	and
	\begin{multline*}
		\bigg| \sum_{k \in \cJ_j} \langle A \bvarphi_{k,n,\e}, \bvarphi_{k,n,\e} \rangle_\LPn
		- \sum_{k \in \cJ_j} \langle A \sn{\phi_{k,\e}}, \sn{\phi_{k,\e}} \rangle_\LPn \bigg|
		\\ \leq C \|A\|_{\op,\LPn} \bigg( 1+ \sqrt{\|\bW_{n,\e}\|_{\op,\LPn}} + \|\bW_{n,\e}\|_{\op,\LPn} \bigg) 
		\Bigg( \sqrt{\frac{\log(4/\tau)}{n \e^{d/2 + 2}}} +  \frac{\log(4/\tau)}{n \e^{d/2 + 1}} \Bigg),
	\end{multline*}
	for any symmetric matrix $A \in \R^{n \times n}$.
\end{lemma}

\begin{proof}[Proof of Lemma~\ref{popul:extra2}]
	Applying Proposition~\ref{gauss:mercer:eigvec2} to $\bW_{n,\e}$ yields the following:
	
	If $n \ge \frac{38 \varrho_\e (1-\e\lambda_{\tmin,j,\e})}{\e^2\gamma_{j,\e}^2}$, with probability at least 
	$$1-\tau-\frac{2\varrho_\e}{\e\gamma_{j,\e}} \exp \Big( -\frac{n\e^2\gamma_{j,\e}^2}{70 \varrho_\e\,(1-\e\lambda_{\tmin,j,\e})} \Big) - 2 \exp \Bigg(-\frac{n\e \gamma_{j,\e}}{2(\sqrt{n \varrho_\e \sum_{k\in\cJ}(1-\e\lambda_{k,\e})}+ \varrho_\e)} \Bigg),$$
	we have 
	\begin{multline*}
		\bigg| \sum_{k \in \cJ} \bmu_{k,n,\e} \langle A \bvarphi_{k,n,\e}, \bvarphi_{k,n,\e} \rangle_\LPn
		- \sum_{k \in \cJ} \mu_{k,\e} \langle A \sn{\phi_{k,\e}}, \sn{\phi_{k,\e}} \rangle_\LPn \bigg|
		\\ \leq C_{|\cJ|} \|A\|_\op \| \bW_n \|_{\op,\LPn} \bigg( \sqrt{ \frac{\varrho_\e \sum_{k \in \cJ} (1-\e\lambda_{k,\e})}{n\e^2\gamma_{j,\e}^2}} +  \frac{\varrho_\e}{n\e \gamma_{j,\e}} \bigg) \log(2/\tau),
	\end{multline*}
	where $C_{|\cJ|}$ is a constant depending only on $|\cJ|$.
	By following the same arguments as those in the proof of Theorem~\ref{popul:eigval}, we obtain the first inequality.
	
	The second inequality can be obtained in a similar way by using Proposition~\ref{gauss:mercer:eigvec3} and additionally considering the uniform control of $\{\phi_{k,\e}\}_{k \in \cJ}$ stated in Theorem~\ref{smooth:eigfun2}.
	This completes the proof.
\end{proof}

\section{First-order asymptotic expansion for spectral convergence of graph Laplacians} \label{graph}

In this section, we prove our main results regarding the first-order asymptotic expansions of the eigenvalues and eigenvectors of $\hL_{n,\e}$.

Throughout Sections~\ref{smooth} and \ref{popul}, we study the spectral convergence of $\bL_{n,\e}$ to $\bfL_\e$ and $\bfL_\e$ to $\Delta_\M$. 
Hence, by further analyzing the spectral properties of $\hL_{n,\e}$ and $\bL_{n,\e}$ and combining these results with the previous analyses, we obtain the desired first-order asymptotic expansions.

To this end, we first establish the operator-norm closeness between $\hL_{n,\e}$ and $\bL_{n,\e}$.
This operator-level closeness leads to a spectral correspondence between the two matrices, which serves as a key step toward proving our main results. We rely on the spectral convergence results regarding symmetric matrices established in \cite{bae2026asymptotic}, which are well suited to our setting, and summarize their relevant forms in Appendix~\ref{gauss:spec}.

\subsection{Perturbation bounds of graph Laplacians}

To formulate the operator-norm closeness between $\hL_{n,\e}$ and $\bL_{n,\e}$, we introduce
\begin{align*}
	\hPi_{n,\e} = \diag\bigg( \frac{\hrho_{n,\e}(X_i)- \rho_\e(X_i)}{\rho_\e(X_i)} \bigg)_{i\in[n]} \quad\text{and}\quad
	\hTheta_{n,\e} = \frac{\hPi_{n,\e} \bW_{n,\e} + \bW_{n,\e} \hPi_{n,\e} }{2\e}.
\end{align*}
The matrix $\hPi_{n,\e}$ characterizes the perturbation in $\hL_{n,\e} - \bL_{n,\e}$, while $\hTheta_{n,\e}$ serves as its first-order asymptotic expansion and plays an important role in the eigenvalues analysis of $\hL_{n,\e}$ in the subsequent subsection.

We formalize this as a lemma below.

\begin{lemma} \label{graph:perturb:aux1}
	It holds that
	\begin{equation}
		\| \bW_{n,\e} \|_{\op,\LPn} \leq 1 + \|\hPi_{n,\e}\|_{\op,\LPn}.
		\label{graph:perturb:aux1-res1}
	\end{equation}
	Moreover, if $\|\hPi_{n,\e}\|_{\op,\LPn} < 1/2$, then
	\begin{gather}
		\|\hTheta_{n,\e} \|_{\op,\LPn} \leq \frac{3\|\hPi_{n,\e}\|_{\op,\LPn}}{2\e}, 
		\label{graph:perturb:aux1-res2} \\
		\|\hL_{n,\e}-\bL_{n,\e} \|_{\op,\LPn} \leq \frac{15\|\hPi_{n,\e}\|_{\op,\LPn}}{4\e}, 
		\label{graph:perturb:aux1-res3} \\
		\big\| \hL_{n,\e}-\bL_{n,\e} - \hTheta_{n,\e} \big\|_{\op,\LPn} 
		\leq \frac{51 \|\hPi_{n,\e}\|^2_{\op,\LPn} }{8\e}.
		\label{graph:perturb:aux1-res4}
	\end{gather}
\end{lemma}

\begin{proof}[Proof of Lemma \ref{graph:perturb:aux1}]
	We express $\hW_{n,\e}$ as
	\begin{align}
		\hW_{n,\e} 
		&= \diag\left( \sqrt{\frac{\rho_\e (X_i)}{\hrho_{n,\e}(X_i)}} \, \right)_{i\in[n]} 
		\left( \frac{\fg_\e (X_i, X_{i'})}{\sqrt{\rho_\e (X_i)} \sqrt{\rho_\e (X_{i'})} } \right)_{i,i' \in [n]}
		\diag\left( \sqrt{\frac{\rho_\e (X_i)}{\hrho_{n,\e}(X_i)}} \, \right)_{i\in[n]} 
		\notag\\&= (I_n + \hPi_{n,\e})^{-1/2} \bW_{n,\e} (I_n + \hPi_{n,\e})^{-1/2}.
		\label{graph:perturb:aux1-eq1}
	\end{align}
	Noting $\|\hW_{n,\e}\|_{\op,\LPn} \leq 1$, we obtain \eqref{graph:perturb:aux1-res1} from
	\begin{align*}
		\| \bW_{n,\e} \|_{\op,\LPn}
		&= \| (I_n + \hPi_{n,\e})^{1/2} \hW_{n,\e} (I_n + \hPi_{n,\e})^{1/2} \|_{\op,\LPn}
		\\&\leq \|I_n+\hPi_{n,\e}\|_{\op,\LPn}
		\leq 1+\|\hPi_{n,\e}\|_{\op,\LPn}.
	\end{align*}
	
	\eqref{graph:perturb:aux1-res2} follows from the definition of $\hTheta_{n,\e}$, \eqref{graph:perturb:aux1-res1}, and the condition $\|\hPi_{n,\e}\|_{\op,\LPn} < 1/2$.
	
	$\eqref{graph:perturb:aux1-res3}$ can be shown in the same way as proving $\eqref{graph:perturb:aux1-res4}$ below.
	
	Under the condition $\|\hPi_{n,\e}\|_{\op,\LPn} < 1/2$, applying Lemma~\ref{gauss:sup:taylor} to \eqref{graph:perturb:aux1-eq1} yields
	\begin{align*}
		\big\| \hL_{n,\e}-\bL_{n,\e} - \hTheta_{n,\e} \big\|_{\op,\LPn} 
		&= \frac{1}{\e} \bigg\| \hW_{n,\e} - \bW_{n,\e} + \frac{\bW_{n,\e} \hPi_{n,\e}+\hPi_{n,\e}\bW_{n,\e} }{2} \bigg\|_{\op,\LPn}
		\\&\leq \frac{17 \|\bW_{n,\e}\|_{\op,\LPn} \|\hPi_{n,\e}\|_{\op,\LPn}^2 }{4\e}.
	\end{align*}
	Thus, applying \eqref{graph:perturb:aux1-res1} yields \eqref{graph:perturb:aux1-res4}, completing the proof.
\end{proof}

The next lemma establishes a concentration bound for the operator norm of $\hPi_{n,\e}$.

\begin{lemma} \label{graph:perturb:aux2}
	There exist $\M$-dependent constants $r$, $C$, $C''$, $C'''$ such that if $\e < r$,
	we have with probability at least $1-\tau$,
	\begin{equation*}
		\| \hPi_{n,\e} \|_{\op,\LPn}
		\leq C \bigg( \sqrt{\frac{\log(1/(\e\wedge\tau))}{n\e^{d/2}}}
		+\frac{\log(1/(\e\wedge\tau))}{n\e^{d/2}} \bigg),
	\end{equation*}
	and for $t>0$ with $t \wedge t^2 \ge \frac{C''' \log(1/\e)}{n\e^{d/2}}$,
	\begin{equation*}
		\P \left (\|\hPi_{n,\e}\|_{\op,\LPn} \geq t \right) 
		\leq \exp \big(- C'' n\e^{d/2} (t \wedge t^2) \big).
	\end{equation*}
\end{lemma}

\begin{proof}[Proof of Lemma~\ref{graph:perturb:aux2}]
	The second claim follows directly from the first and a standard inversion of the high-probability bound.
	
	To show the first, note that
	\begin{equation*}
		\|\hPi_{n,\e}\|_{\op,\LPn} 
		= \max_{i=1,\dots,n} \left| \frac{\hrho_{n,\e}(X_i)- \rho_\e(X_i)}{\rho_\e(X_i)} \right|.
	\end{equation*}
	For the denominator, by Lemma~\ref{smooth:KDE}, we can choose a $\M$-dependent constant $r_1$ such that $\e<r_1$ implies $\sup_{x \in \M} \frac{1}{\rho_\e(x)} \leq 2.$
	For the numerator, note that $\hrho_{n,\e}$ is the KDE corresponding to the smoothed KDE $\rho_\e$. Hence, by the uniform convergence of KDE on $\M$ (e.g., Theorem 3.6 of \cite{bae2026kernel}), 
	we can choose $\M$-dependent constants $r_2$ and $C_2$ such that if $\e < r_2$, then with probability at least $1-\tau$, it holds that
	\begin{equation*}
		\sup_{x \in \M} |\hrho_{n,\e}(x)- \rho_\e(x)| 
		\leq C \left( \sqrt{\frac{\log(1/(\e\wedge\tau))}{n\e^{d/2}}}+\frac{\log(1/(\e\wedge\tau))}{n\e^{d/2}} \right).
	\end{equation*}
	Therefore, the proof is complete by adjusting $r$ and $C$ appropriately.
\end{proof}

With the above two lemmas, we quantify perturbation of $\hL_{n,\e}- \bL_{n,\e}$ and its first-order approximation $\hTheta_{n,\e}$.

\begin{proposition} \label{graph:perturb}
	There exist $\M$-dependent constants $r$, $C$, $C''$, $C'''$ such that 
	if $\e < r$ and $\frac{n\e^{d/2}}{\log(1/\e)} > C'''$, we have with probability at least $1- \tau - \exp \big(-C''n\e^{d/2} \big)$,
	\begin{gather}
		\|\bW_{n,\e} \|_{\op,\LPn} \leq \frac{3}{2},
		\label{graph:perturb-res1} \\
		\|\hTheta_{n,\e}\|_{\op,\LPn},\, \| \hL_{n,\e} - \bL_{n,\e} \|_{\op,\LPn}
		\leq C \left( \sqrt{\frac{\log (1/(\e\wedge\tau))}{n\e^{d/2+2}}} + \frac{\log (1/(\e\wedge\tau))}{n\e^{d/2+1}} \right),
		\label{graph:perturb-res2} \\
		\| \hL_{n,\e} - \bL_{n,\e} - \hTheta_{n,\e}\|_{\op,\LPn}
		\leq C \left( \sqrt{\frac{\log (1/(\e\wedge\tau))}{n\e^{d/2+2}}} + \frac{\log (1/(\e\wedge\tau))}{n\e^{d/2+1}} \right)^2.
		\label{graph:perturb-res3}
	\end{gather}
	Moreover, for $t>0$ with $(t \wedge t^2) \ge \frac{C'''_1\log(1/\e)}{n\e^{d/2+2}}$,
	\begin{equation}
		\P\big( \| \hL_{n,\e} - \bL_{n,\e} \|_{\op,\LPn} \|_{\op,\LPn} \ge t \big) \leq \exp( -C''_1 n\e^{d/2+2} (t \wedge t^2) ) + \exp(-C'' n\e^{d/2}),
		\label{graph:perturb-res4}
	\end{equation}
	provided $\e <r$, where $C''_1$ is $\M$-dependent.
\end{proposition}

\begin{proof}[Proof of Proposition \ref{graph:perturb}]
	Let $r$ satisfy Lemma~\ref{graph:perturb:aux2}.
	By adjusting $C'''$ appropriately and plugging $t=1/2$ into Lemma~\ref{graph:perturb:aux2}, we obtain
	\begin{equation*}
		\P ( \|\hPi_{n,\e}\|_{\op,\LPn} \geq 1/2 ) \leq \exp(-C''_2 n\e^{d/2}).
	\end{equation*}	
	Therefore, we obtain \eqref{graph:perturb-res1}-\eqref{graph:perturb-res3} by applying Lemmas~\ref{graph:perturb:aux1} and \ref{graph:perturb:aux2}.
	
	\eqref{graph:perturb-res4} follows from \eqref{graph:perturb-res2} and a standard inversion argument.
	This completes the proof.
\end{proof}

\subsection{Eigenvalues of graph Laplacians} \label{graph:eigval}

The goal of this subsection is to prove Theorem~\ref{gauss:eigval}.
To achieve this, we first provide an initial first-order expansion of $\sum_{k \in \cJ_j} (\hlambda_{k,n,\e}-\blambda_{k,n,\e})$ as below.

\begin{lemma} \label{graph:eigval:aux1}
	There exist $(\M,\lambda_j)$-dependent constant $r$, $C$, $C'$, $C''$, $C'''$ such that 
	if $\e < r$ and	$ \frac{n\e^{d/2+2}}{\log(1/\e)} > C'''$,
	we have with probability at least $1- \tau - C'\exp \left( -C''n\e^{d/2+2} \right)$,
	\begin{equation*}
		\left| \sum_{k \in \cJ_j} (\hlambda_{k,n,\e}-\blambda_{k,n,\e})
		- \sum_{k \in \cJ_j} \langle \hTheta_{n,\e} \bvarphi_{k,n,\e}, \bvarphi_{k,n,\e} \rangle_\LPn \right|
		\leq C \left( \sqrt{\frac{\log(1/(\e\wedge\tau))}{n\e^{d/2+2}}} + \frac{\log(1/(\e\wedge\tau))}{n\e^{d/2+1}} \right)^2
	\end{equation*}
\end{lemma}

\begin{proof}[Proof of Lemma \ref{graph:eigval:aux1}]
	Let $r$ and $C'''$ satisfy Lemma~\ref{popul:extra1} and Proposition~\ref{graph:perturb}.	
	We decompose and bound
	\begin{multline*}
		\bigg| \sum_{k\in\cJ_j} \hlambda_{k,n,\e} - \sum_{k\in\cJ_j} \blambda_{k,n,\e}
		- \sum_{k \in \cJ_j} \langle \hTheta_{n,\e} \bvarphi_{k,n,\e}, \bvarphi_{k,n,\e} \rangle_\LPn \bigg|
		\\ \leq \bigg| \sum_{k\in\cJ_j} \hlambda_{k,n,\e} - \sum_{k\in\cJ_j} \blambda_{k,n,\e}
		- \sum_{k \in \cJ_j} \langle (\hL_{n,\e}-\bL_{n,\e}) \bvarphi_{k,n,\e}, \bvarphi_{k,n,\e} \rangle_\LPn \bigg| 
		\\ + \bigg| \sum_{k \in \cJ_j} \langle (\hL_{n,\e}-\bL_{n,\e} - \hTheta_{n,\e}) \bvarphi_{k,n,\e}, \bvarphi_{k,n,\e} \rangle_\LPn \bigg|
		=: E_1 + E_2.
	\end{multline*}
	
	Applying Proposition~\ref{gauss:spec:eigval3} to $\hL_{n,\e}$ yields
	\begin{equation*}
		E_1	\leq \bigg( \frac{5}{2\bgamma_{j,n,\e}} + \frac{4 \blambda_{\max,j,\e} - 4 \blambda_{\min,j,\e} }{\bgamma_{j,n,\e}^2} \bigg) 
		\|\hL_{n,\e}-\bL_{n,\e}\|_{\op,\LPn}^2,
	\end{equation*}
	provided that
	\begin{equation*}
		\frac{\|\hL_{n,\e}-\bL_{n,\e}\|_{\op,\LPn}}{\bgamma_{j,n,\e}} \leq \frac{1}{4}.
	\end{equation*}
	Using Lemma~\ref{popul:extra1}, we have with probability at least $1-C'_1 \exp(-C''_1 n\e^{d/2+2})$,
	\begin{equation*}
		\frac{5}{2\bgamma_{j,n,\e}} + \frac{4 \blambda_{\max,j,\e} - 4 \blambda_{\min,j,\e} }{\bgamma_{j,n,\e}^2} \le C_1, \quad
		\frac{1}{\bgamma_{j,n,\e}} \leq C_1.
	\end{equation*}
	Hence, applying Proposition~\ref{graph:perturb} gives with probability at least $1-\tau-C'_2 \exp(-C''_2 n\e^{d/2+2})$,
	\begin{equation}
		E_1 \leq C_2 \left( \sqrt{\frac{\log(1/(\e\wedge\tau))}{n\e^{d/2+2}}} + \frac{\log(1/(\e\wedge\tau))}{n\e^{d/2+1}} \right)^2.
		\label{graph:eigval:aux1-eq1}
	\end{equation}
	
	Next, using that $\{\bvarphi_{k,n,\e}\}_{k\in\cJ_j}$ are $\LPn$-orthonormal and applying Proposition~\ref{graph:perturb}, it follows that with probability at least $1- \tau - C_3'\exp ( -C''_3 n\e^{d/2} )$,
	\begin{align}
		E_2 \leq |\cJ_j| \cdot \|\hL_{n,\e} - \bL_{n,\e} - \hTheta_{n,\e} \|_{\op,\LPn}
		\leq C_3 \left( \sqrt{\frac{\log(1/(\e\wedge\tau))}{n\e^{d/2+2}}} + \frac{\log(1/(\e\wedge\tau))}{n\e^{d/2+1}} \right)^2.
		\label{graph:eigval:aux1-eq2}
	\end{align}
	
	Therefore, the claim follows by \eqref{graph:eigval:aux1-eq1} and \eqref{graph:eigval:aux1-eq2}.
\end{proof}

Now, our goal is to approximate the term $\sum_{k \in \cJ_j} \langle \hTheta_{n,\e} \bvarphi_{k,n,\e}, \bvarphi_{k,n,\e} \rangle_\LPn$ by a more tractable expression.
To this end, we replace the random quantities $\{\bvarphi_{k,n,\e}\}_{k \in \cJ_j}$ and $\hTheta_{n,\e}$ with deterministic quantities,
and this is carried out in the following two lemmas.


\begin{lemma} \label{graph:eigval:aux2}
	There exist $(\M,\lambda_j)$-dependent constant $r$, $C$, $C'$, $C''$, $C'''$ such that 
	if $\e < r$ and	$ \frac{n\e^{d/2+2}}{\log(1/\e)} > C'''$,
	we have with probability at least $1- \tau - C'\exp \left( -C'' n\e^{d/2+2} \right)$,	
	\begin{multline*}
		\bigg| \sum_{k \in \cJ_j} \langle \hTheta_{n,\e} \bvarphi_{k,n,\e}, \bvarphi_{k,n,\e} \rangle_\LPn
		- \sum_{k \in \cJ_j} \frac{\mu_{k,\e}}{\e} \langle \hPi_{n,\e} \sn{\phi_{k,\e}}, \sn{\phi_{k,\e}} \rangle_\LPn \bigg|
		\\ \leq C \Bigg( \sqrt{\frac{\log(1/(\e\wedge\tau))}{n\e^{d/2+2}}} + \frac{\log(1/(\e\wedge\tau))}{n\e^{d/2+1}} \Bigg)^2.
	\end{multline*}
\end{lemma}

\begin{proof}[Proof of Lemma~\ref{graph:eigval:aux2}]
	Let $r$ and $C'''$ satisfy Lemmas~\ref{popul:extra2} and \ref{graph:perturb:aux2} and Proposition~\ref{graph:perturb}.
	From $\bW_{n,\e} \bvarphi_{k,n,\e} = \bmu_{k,n,\e} \bvarphi_{k,n,\e}$, we have
	\begin{equation*}
		\langle \hTheta_{n,\e} \bvarphi_{k,n,\e}, \bvarphi_{k,n,\e} \rangle_\LPn
		= \frac{1}{\e} \langle \hPi_{n,\e} \bvarphi_{k,n,\e}, \bW_{n,\e} \bvarphi_{k,n,\e} \rangle_\LPn
		= \frac{\bmu_{k,n,\e}}{\e} \langle \hPi_{n,\e} \bvarphi_{k,n,\e}, \bvarphi_{k,n,\e} \rangle_\LPn.
	\end{equation*}
	On the other hand, applying Lemma~\ref{popul:extra2} gives with probability at least $1-\tau-C_1'\exp(-C_1''n\e^{d/2+2})$.
	\begin{multline}
		\bigg| \sum_{k \in \cJ_j} \bmu_{k,n,\e} \langle \hPi_{n,\e} \bvarphi_{k,n,\e}, \bvarphi_{k,n,\e} \rangle_\LPn
		- \sum_{k \in \cJ_j} \mu_{k,\e} \langle \hPi_{n,\e} \sn{\phi_{k,\e}}, \sn{\phi_{k,\e}} \rangle_\LPn \bigg|
		\\ \leq C \|\hPi_{n,\e}\|_{\op,\LPn} \| \bW_{n,\e} \|_{\op,\LPn} \Bigg( \sqrt{\frac{\log(2/\tau)}{n \e^{d/2 + 2}}} +  \frac{\log(2/\tau)}{n \e^{d/2 + 1}} \Bigg).
		\label{graph:eigval:aux2-eq2}
	\end{multline}	
	Therefore, we complete the proof by diving both sides of \eqref{graph:eigval:aux2-eq2} by $\e$ and applying Lemma~\ref{graph:perturb:aux2} and Proposition~\ref{graph:perturb}.
\end{proof}

\begin{lemma} \label{graph:eigval:aux3}
	There exist $(\M,\lambda_j)$-dependent constant $r$ and $C$ such that 
	if $\e < r$, we have with probability at least $1- \tau$,
	\begin{multline*}
		\bigg| \sum_{k \in \cJ_j} \frac{\mu_{k,\e}}{\e} \langle \hPi_{n,\e} \sn{\phi_{k,\e}}, \sn{\phi_{k,\e}} \rangle_\LPn
		- (\P_n -\P) \bigg( \sum_{k \in \cJ_j} \frac{\mu_{k,\e}\phi_{k,\e}^2}{\e} - \sum_{k \in \cJ_j} \mu_{k,\e} \bfL_\e^\unnor \Big( \frac{\phi^2_{k,\e}}{\rho_\e} \Big) \bigg) \bigg|
		\\ \leq C \bigg( \sqrt{\frac{\log(1/(\e\wedge\tau))}{n\e^{d/2+2}}} + \frac{\log(1/(\e\wedge\tau))}{n\e^{d/2+1}} \bigg)^2.
	\end{multline*}
\end{lemma}

\begin{proof}[Proof of Lemma~\ref{graph:eigval:aux3}]
	We show the high-probability bound for
	\begin{equation*}
		\frac{\mu_{k,\e}}{\e} \langle \hPi_{n,\e} \sn{\phi_{k,\e}}, \sn{\phi_{k,\e}} \rangle_\LPn
		- (\P_n - \P) \bigg( \frac{\mu_{k,\e}\phi_{k,\e}^2}{\e} - \mu_{k,\e} \bfL_\e^\unnor \Big( \frac{\phi^2_{k,\e}}{\rho_\e} \Big) \bigg),
	\end{equation*}
	and the result follows by taking the sum over $k \in \cJ_j$.
	We split the proof into several steps
	
	\mycase{1. Hoeffding decomposition}
	Using the definitions of $\hPi_{n,\e}$ and $\hrho_{n,\e}$, we expand
	\begin{align*}
		\langle \hPi_{n,\e} \sn{\phi_{k,\e}}, \sn{\phi_{k,\e}} \rangle_\LPn
		&= \frac{1}{n} \sum_{i=1}^n \frac{\hrho_{n,\e}(X_i) - \rho_\e (X_i)}{\rho(X_i)}  \phi_{k,\e}^2 (X_i)
		\\&= \frac{1}{n^2} \sum_{i=1}^n \sum_{i'=1}^n \frac{\fg_\e(X_i,X_{i'}) - \rho_\e (X_i)}{\rho(X_i)}  \phi_{k,\e}^2 (X_i).
	\end{align*}
	By letting
	\begin{equation}
		f(x,y) = \frac{\mu_{k,\e}}{2\e} \bigg(\frac{\fg_\e (x,y) - \rho_\e(x)}{\rho_\e(x)}\bigg) \phi_{k,\e}^2 (x)
		+ \frac{\mu_{k,\e}}{2\e} \bigg(\frac{\fg_\e (y,x) - \rho_\e(y) }{\rho_\e(y)}\bigg) \phi_{k,\e}^2 (y),
		\label{graph:eigval:aux3-eq2}
	\end{equation}
	it follows that
	\begin{equation}
		\frac{\mu_{k,\e}}{\e} \langle \hPi_{n,\e} \sn{\phi_{k,\e}}, \sn{\phi_{k,\e}} \rangle_\LPn 
		= \frac{1}{n^2} \sum_{i=1}^n \sum_{i'=1}^n f(X_i,X_{i'}),
		\label{graph:eigval:aux3-eq1}
	\end{equation}
	which is a symmetric $V$-statistic.
	For i.i.d. copies $X$ and $Y$ of $X_1$, let
	\begin{equation*}
		\theta  = \E f(X, Y), \quad 
		f_1 (x) = \E f(x, Y) - \theta, \quad 
		f_2 (x,y) = f(x, y) - f_1 (x) - f_1 (y) - \theta.
	\end{equation*}
	Then, we have the following Hoeffding decomposition
	\begin{equation}
		\frac{1}{n^2} \sum_{i=1}^n \sum_{i'=1}^n f(X_i,X_{i'})
		= \underbrace{\theta + \frac{2}{n} \sum_{i=1}^n f_1(X_i)}_{E_1}
		+ \underbrace{\frac{1}{n^2} \sum_{i=1}^n f_2(X_i, X_i)}_{E_2}
		+ \underbrace{\frac{1}{n^2} \sum_{1\leq i \ne i' \le n} f_2(X_i, X_{i'})}_{E_3}.
		\label{graph:eigval:aux3-eq3}
	\end{equation}
	Below, we analyze each term.
	
	\mycase{2. Evaluating $E_1$}
	Since $\int_\M \fg_\e (x,y) \dP(y) = \rho_\e (x)$ and $\int_\M \phi_{k,\e}^2 (y)\dP(y) = 1$, we have
	\begin{equation*}
		\E f(x,Y) = \frac{\mu_{k,\e}}{2\e} \int_\M \frac{\fg_\e(y,x) - \rho_\e(y)}{\rho_\e(y)} \phi_{k,\e}^2 (y) \dP(y)
		= \frac{\mu_{k,\e}}{2\e} \int_\M \fg_\e(y,x) \frac{\phi_{k,\e}^2(y)}{\rho_\e (y)} \dP(y) - \frac{\mu_{k,\e}}{2\e}.
	\end{equation*}
	This yields
	\begin{equation}
		\theta = 0 \quad\text{and}\quad
		f_1 (x) = \frac{\mu_{k,\e}}{2\e}\int_\M \fg_\e(y,x) \frac{\phi_{k,\e}^2(y)}{\rho_\e (y)} \dP(y) - \frac{\mu_{k,\e}}{2\e}.
		\label{graph:eigval:aux3-eq4}
	\end{equation}
	From this, we have
	\begin{align*}
		- \mu_{k,\e} \bfL_\e^\unnor \Big( \frac{\phi^2_{k,\e}}{\rho_\e} \Big) (x)
		&= \frac{\mu_{k,\e}}{\e} \int_\M \fg_\e (x,y) \left( \frac{\phi^2_{k,\e}(y)}{\rho_\e (y)} - \frac{\phi^2_{k,\e}(x)}{\rho_\e (x)} \right) \dP (y)
		\\&= \frac{\mu_{k,\e}}{\e} \int_\M \fg_\e(y,x) \frac{\phi_{k,\e}^2(y)}{\rho_\e (y)} \dP(y) - \frac{\mu_{k,\e}\phi^2_{k,\e} (x)}{\e}
		= 2f_1(x) - \frac{\mu_{k,\e}}{\e} (\phi_{k,\e}^2(x) -1 ),
	\end{align*}
	and it follows that
	\begin{align}
		E_1 &= \frac{2}{n} \sum_{i=1}^n f_1 (X_i)
		= \frac{1}{n} \sum_{i=1}^n \bigg( \frac{\mu_{k,\e} (\phi^2_{k,\e} (X_i)-1)}{\e}  -\mu_{k,\e} \bfL_\e^\unnor \Big( \frac{\phi^2_{k,\e}}{\rho_\e} \Big) (X_i) \bigg)
		\notag\\&= (\P_n -\P) \bigg(\frac{\mu_{k,\e}\phi_{k,\e}^2}{\e} - \mu_{k,\e} \bfL_\e^\unnor \Big( \frac{\phi^2_{k,\e}}{\rho_\e} \Big) \bigg),
		\label{graph:eigval:aux3-eq5}
	\end{align}
	where the last equality uses and $\P\, \bfL_\e^\unnor \Big( \frac{\phi^2_{k,\e}}{\rho_\e} \Big) =0 $ and $\P \phi_{k,\e}^2 = 1$.
	
	\mycase{3. Bounding $E_2$}
	From $\fg_\e(x,y) \leq (4\pi\e)^{-d/2}$, $\mu_{k,\e} = 1- \e \lambda_{k,\e}$ and $\|\rho_\e -1\|_\infty \leq C_1 \e$ from Lemma~\ref{smooth:KDE}, applying Theorems~\ref{smooth:eigval} and \ref{smooth:eigfun2} gives for $\e < r_2$,
	\begin{equation}
		\|f\|_\infty, \|f_2\|_\infty \leq \frac{C_2}{\e^{d/2+1}}.
		\label{graph:eigval:aux3-eq6}
	\end{equation}
	Hence, it follows that
	\begin{equation*}
		|E_2| \leq \frac{1}{n} \|f_2\|_\infty \leq \frac{C_3}{n\e^{d/2+1}}.
	\end{equation*}
	
	\mycase{4. Bounding $E_3$}
	Note that $E_2$ is a $U$-statistic with kernel $f_2$, which is $\P$-canonical, that is, $\E f_2(x,Y) = \E f_1 (X,y) = 0$ for all $x, y \in \M$.
	By applying the concentration inequality for a $U$-statistic with canonical kernel stated in Lemma~\ref{gauss:sup:ustat}, we have with probability at least $1-\tau$,
	\begin{equation}
		|E_3| \leq C_4 \bigg( \frac{ \sqrt{\E f_2^2(X,Y)} \log(1/\tau)}{n}+ \frac{\|f_2\|_\infty \log^{3/2}(1/\tau)}{n^{3/2}} \bigg).
		\label{graph:eigval:aux3-eq7}
	\end{equation}	
	where $C_4$ is a universal constant.
	We can bound $\|f_2\|_\infty$ using \eqref{graph:eigval:aux3-eq6}, and it remains to bound $\E f_2^2(X,Y)$.
	
	From \eqref{graph:eigval:aux3-eq2} and \eqref{graph:eigval:aux3-eq4}, we evaluate $f_2$ as
	\begin{align*}
		f_2 (x,y) 
		&= \frac{\mu_{k,\e}}{2\e} \bigg(\frac{\fg_\e (x,y) - \rho_\e(x)}{\rho_\e(x)}\bigg) \phi_{k,\e}^2 (x)
		+ \frac{\mu_{k,\e}}{2\e} \bigg(\frac{\fg_\e (x,y) - \rho_\e(y) }{\rho_\e(y)}\bigg) \phi_{k,\e}^2 (y)
		\\&\quad - \frac{\mu_{k,\e}}{2\e}\int_\M \fg_\e(x,z) \frac{\phi_{k,\e}^2(z)}{\rho_\e (z)} \dP(z) + \frac{\mu_{k,\e}}{2\e} 
		- \frac{\mu_{k,\e}}{2\e}\int_\M \fg_\e(y,z) \frac{\phi_{k,\e}^2(z)}{\rho_\e (z)} \dP(z) + \frac{\mu_{k,\e}}{2\e}
		\\&= \frac{\mu_{k,\e}}{2\e} \bigg( \fg_\e(x,y) \frac{\phi_{k,\e}^2(y)}{\rho_\e(y)} - \int_\M \fg_\e(x,z) \frac{\phi_{k,\e}^2(z)}{\rho_\e(z)} \dP(z) + \phi_{k,\e}^2(x) - 1
		\\&\quad\quad\quad + \fg_\e(y,x) \frac{\phi_{k,\e}^2(x)}{\rho_\e(x)} - \int_\M \fg_\e(y,z) \frac{\phi_{k,\e}^2(z)}{\rho_\e(z)} \dP(z) + \phi_{k,\e}^2(y) - 1 \bigg),
	\end{align*}
	and this gives
	\begin{align*}
		f_2^2 (x,y) &\leq \frac{\mu_{k,\e}^2}{\e^2} \Bigg\{ \bigg( \fg_\e(x,y) \frac{\phi_{k,\e}^2(y)}{\rho_\e(y)} - \int_\M \fg_\e(x,z) \frac{\phi_{k,\e}^2(z)}{\rho_\e(z)} \dP(z)\bigg)^2 + (\phi_{k,\e}^2(x) - 1)^2
		\\&\quad  + \bigg( \fg_\e(y,x) \frac{\phi_{k,\e}^2(x)}{\rho_\e(x)} - \int_\M \fg_\e(y,z) \frac{\phi_{k,\e}^2(z)}{\rho_\e(z)} \dP(z) \bigg)^2 + (\phi_{k,\e}^2(y) - 1)^2 \bigg\}.
	\end{align*}
	Since each term above is centered, we obtain
	\begin{align*}
		\E f_2^2 (X,Y) &\leq \frac{2 \mu_{k,\e}^2}{\e^2} \bigg( \int_\M \int_\M \fg^2_\e(x,y) \frac{\phi_{k,\e}^4(y)}{\rho^2_\e(y)} \dP(x) \dP(y) + \int_\M \phi_{k,\e}^4(x) \dP(x) \bigg)
		\\&\leq \frac{2 \mu_{k,\e}^2}{\e^2} \bigg( \|\fg_\e\|_\infty \int_\M \int_\M \fg_\e(x,y) \frac{\phi_{k,\e}^4(y)}{\rho_\e^2(y)} \dP(x) \dP(y) + \|\phi_{k,\e}\|_\infty^2 \int_\M \phi_{k,\e}^2 (x) \dP(x) \bigg)
		\\&= \frac{2 \mu_{k,\e}^2}{\e^2} \bigg( \|\fg_\e\|_\infty \int_\M \frac{\phi_{k,\e}^3(y)}{\rho_\e(y)} \dP(y) + \|\phi_{k,\e}\|_\infty^2 \bigg).
	\end{align*}
	By the same argument used to derive \eqref{graph:eigval:aux3-eq6}, it follows that
	\begin{equation*}
		\E f_2^2 (X,Y) \leq \frac{C_5}{\e^{d/2+2}},
	\end{equation*}	
	and we can replace \eqref{graph:eigval:aux3-eq7} with
	\begin{equation}
		|E_3| \leq C_6 \bigg( \frac{\log(1/\tau)}{n \e^{d/4+1}} + \frac{\log^{3/2}(1/\tau)}{n^{3/2}\e^{d/2+1}} \bigg).
		\label{graph:eigval:aux3-eq8}
	\end{equation}
	
	\mycase{5. Concluding the proof}
	By combining \eqref{graph:eigval:aux3-eq1} \eqref{graph:eigval:aux3-eq3}, \eqref{graph:eigval:aux3-eq5} and \eqref{graph:eigval:aux3-eq8}, if $\e < r_2$, we have with probability at least $1-\tau$,
	\begin{multline*}
		\bigg| \frac{\mu_{k,\e}}{\e} \langle \hPi_{n,\e} \sn{\phi_{k,\e}}, \sn{\phi_{k,\e}} \rangle_\LPn
		- (\P_n - \P) \bigg( \frac{\mu_{k,\e}\phi_{k,\e}^2}{\e} - \mu_{k,\e} \bfL_\e^\unnor \Big( \frac{\phi^2_{k,\e}}{\rho_\e} \Big) \bigg) \bigg|
		\\ \leq C_7 \bigg( \frac{1}{n\e^{d/2+1}} + \frac{\log(1/\tau)}{n \e^{d/4+1}} + \frac{\log^{3/2}(1/\tau)}{n^{3/2}\e^{d/2+1}} \bigg).
	\end{multline*}
	By letting $r=1 \wedge r_2$, it can be easily shows that the right-hand side of the above inequality is dominated by
	\begin{equation*}
		C_8 \bigg( \sqrt{\frac{\log(1/(\e\wedge\tau))}{n\e^{d/2+2}}} + \frac{\log(1/(\e\wedge\tau))}{n\e^{d/2+1}} \bigg)^2.
	\end{equation*}
	By taking the sum over $k \in \cJ_j$, we conclude the proof.
\end{proof}

By simply combining Lemmas~\eqref{graph:eigval:aux1}-\eqref{graph:eigval:aux3}, we have the following result.

\begin{proposition} \label{graph:eigval:aux4}
	There exist $(\M,\lambda_j)$-dependent constant $r$, $C$, $C'$, $C''$, $C'''$ such that 
	if $\e < r$ and	$ \frac{n\e^{d/2+2}}{\log(1/\e)} > C'''$,
	we have with probability at least $1- \tau - C'\exp ( -C'' n\e^{d/2+2} )$,
	\begin{multline*}
		\bigg| \sum_{k \in \cJ_j} (\hlambda_{k,n,\e} - \blambda_{k,n,\e})
		- (\P_n -\P) \bigg( \sum_{k \in \cJ_j} \frac{\mu_{k,\e}\phi_{k,\e}^2}{\e} - \sum_{k \in \cJ_j} \mu_{k,\e} \bfL_\e^\unnor \Big( \frac{\phi^2_{k,\e}}{\rho_\e} \Big) \bigg) \bigg|
		\\ \leq C \bigg( \sqrt{\frac{\log(1/(\e\wedge\tau))}{n\e^{d/2+2}}} + \frac{\log(1/(\e\wedge\tau))}{n\e^{d/2+1}} \bigg)^2.
	\end{multline*}
\end{proposition}

Finally, we are ready to prove our main result for the eigenvalues.

\begin{proof}[Proof of Theorem \ref{gauss:eigval}]
	\eqref{gauss:eigval-res2} is obtained from \eqref{gauss:eigval-res1} and $\sum_{k \in \cJ_j} |\lambda_{k,\e}-\lambda_j| \leq C_1 \e$ in Theorem~\ref{smooth:eigval}.
	
	To show \eqref{gauss:eigval-res1}, the expression inside the absolute value on the left-hand side of \eqref{gauss:eigval-res1} can be decomposed as
	\begin{align*}
		\sum_{k \in \cJ_j} (\hlambda_{k,n,\e} - \lambda_{k,\e}) - (\P_n-\P) \bigg( -\sum_{k\in\cJ_j} \Delta_\M(\phi_k^2) \bigg)
		= E_1 + E_2 + E_3,
	\end{align*}
	where
	\begin{align*}
		E_1 &= -(\P_n - \P) \bigg( \sum_{k \in \cJ_j} \mu_{k,\e} \bfL_\e^\unnor \Big(\frac{\phi_{k,\e}^2}{\rho_\e}\Big)	- \sum_{k \in \cJ_j} \Delta_\M (\phi^2_k) \bigg), \\
		E_2 &= \sum_{k \in \cJ_j} (\blambda_{k,n,\e} - \lambda_{k,\e})
		- (\P_n-\P) \bigg( -\sum_{k \in \cJ_j} \frac{\mu_{k,\e} \phi_{k,\e}^2}{\e} \bigg), \\
		E_3 &= \sum_{k \in \cJ_j} (\hlambda_{k,n,\e}-\blambda_{k,n,\e}) - (\P_n - \P) \bigg( \sum_{k \in \cJ_j} \frac{\mu_{k,\e} \phi_{k,\e}^2}{\e} -\sum_{k \in \cJ_j} \mu_{k,\e} \bfL_\e^\unnor \Big(\frac{\phi_{k,\e}^2}{\rho_\e}\Big) \bigg).
	\end{align*}
	We can bound $E_1$, $E_2$, $E_3$ using Lemma~\ref{smooth:extra2}, Theorem~\ref{popul:eigval}, Proposition~\ref{graph:eigval:aux4}, respectively.
	Note that the sum of the bounds is dominated by the right-hand side of \eqref{gauss:eigval-res2}, considering the elementary inequality  $(a+b)^2 \leq 2a^2 + 2b^2$.
	This yields \eqref{gauss:eigval-res1}, completing the proof.
\end{proof}

\subsection{Eigenvectors of graph Laplacians} \label{graph:eigvec}

The goal of this subsection is to prove Theorem~\ref{gauss:eigvec}.
To this end, we need the following result.

\begin{proposition} \label{graph:eigvec:aux1}
	There exist $(\M,\lambda_j)$-dependent constants $r$, $C$, $C'$, $C''$, $C'''$ such that if $\e < r$ and $ \frac{n\e^{d/2+2}}{\log(1/\e)} > C'''$,
	we have with probability at least $1-C'\exp \left( -C''n\e^{d/2+2} \right)$,
	\begin{multline*}
		\bigg| \sum_{k \in \cJ} \sum_{\ell \in \cJ} \langle \hvarphi_{k,n,\e}, \sn{\phi_{\ell,\e}} \rangle_\LPn^2
		- \sum_{k \in \cJ} \sum_{\ell \in \cJ} \langle \bvarphi_{k,n,\e}, \sn{\phi_{\ell,\e}} \rangle_\LPn^2 \bigg|
		\\ \leq C \bigg( \sqrt{\frac{\log(1/(\e\wedge\tau))}{n\e^{d/2+2}}} + \frac{\log(1/(\e\wedge\tau))}{n\e^{d/2+1}} \bigg)^2.
	\end{multline*}
\end{proposition}

\begin{proof}[Proof of Proposition~\ref{graph:eigvec:aux1}]
	Let $r$ and $C'''$ be as in Theorem~\ref{smooth:eigfun2}, Lemmas~\ref{popul:extra1} and \ref{popul:extra2}, and Proposition~\ref{graph:perturb}, and let $\e <r$.
	
	Applying Proposition~\ref{gauss:spec:eigvec3} to $\hL_{n,\e}$ with $u_\ell = \sn{\phi_{\ell,\e}}$ yields
	\begin{multline}
		\bigg| \sum_{k \in \cJ_j} \sum_{\ell \in \cJ_j} \langle \hvarphi_{k,n,\e}, \sn{\phi_{\ell,\e}} \rangle_\LPn^2
		- \sum_{k \in \cJ_j} \sum_{\ell \in \cJ_j} \langle \bvarphi_{k,n,\e}, \sn{\phi_{\ell,\e}} \rangle_\LPn^2 \bigg|
		\\ \leq \bigg| \sum_{\ell \in \cJ_j} \langle \hS_{j,n,\e} \sn{\phi_{\ell,\e}}, \sn{\phi_{\ell,\e}} \rangle_\mbH \bigg|
		+ \bigg( \sum_{\ell \in \cJ_j} \|\sn{\phi_{\ell,\e}}\|_{\LPn}^2 \bigg) (8|\cJ_j|+32) \bigg(\frac{\|\hL_{n,\e}-\bL_{n,\e}\|_{\op,\LPn}}{\bgamma_{j,n,\e}}\bigg)^2,
		\label{graph:eigvec:aux1-eq2}
	\end{multline}
	where $\hS_{j,n,\e}$ satisfies
	\begin{equation}
		\| \hS_{j,n,\e} \|_{\op,\LPn} \leq 2 \sqrt{D} \bigg(\frac{\|\hL_{n,\e}-\bL_{n,\e}\|_{\op,\LPn}}{\bgamma_{j,n,\e}}\bigg), \quad  \langle \hS_{j,n,\e} \bvarphi_{\ell,n,\e}, \bvarphi_{\ell,n,\e} \rangle_\LPn = 0 \quad\text{for all $\ell \in \cJ_j$.}
		\label{graph:eigvec:aux1-eq1}
	\end{equation}
	By Lemma~\ref{popul:extra2} and \eqref{graph:eigvec:aux1-eq1}, with probability at least $1-\tau-C_1'\exp(-C_1'' n\e^{d/2+2})$,
	the first term on the right-hand side of \eqref{graph:eigvec:aux1-eq2} is bounded as
	\begin{multline*}
		\bigg| \sum_{\ell \in \cJ_j} \langle \hS_{j,n,\e} \sn{\phi_{\ell,\e}}, \sn{\phi_{\ell,\e}} \rangle_\mbH \bigg|
		= \bigg| \sum_{\ell \in \cJ_j} \langle \hS_{j,n,\e} \bvarphi_{\ell,n,\e}, \bvarphi_{\ell,n,\e} \rangle_\LPn
		- \sum_{\ell \in \cJ_j} \langle \hS_{j,n,\e} \sn{\phi_{\ell,\e}}, \sn{\phi_{\ell,\e}} \rangle_\LPn \bigg|
		\\\leq C_1 \|\hS_{j,n,\e} \|_{\op,\LPn} \bigg(1+\sqrt{\|\bW_{n,\e}\|_{\op,\LPn}}+{\|\bW_{n,\e}\|_{\op,\LPn}} \bigg) \bigg( \sqrt{\frac{\log(4/\tau)}{n\e^{d/2+2}}} + \frac{\log(4/\tau)}{n\e^{d/2+1}} \bigg)
		\\\leq 2 C_1 \sqrt{D} \bigg(1+\sqrt{\|\bW_{n,\e}\|_{\op,\LPn}}+{\|\bW_{n,\e}\|_{\op,\LPn}} \bigg)
		\bigg(\frac{\|\hL_{n,\e}-\bL_{n,\e}\|_{\op,\LPn}}{\bgamma_{j,n,\e}}\bigg) \\ \times \bigg( \sqrt{\frac{\log(4/\tau)}{n\e^{d/2+2}}} + \frac{\log(4/\tau)}{n\e^{d/2+1}} \bigg).
	\end{multline*}	
	Now, all the bounds can be controlled using Theorem~\ref{smooth:eigfun2}, Lemma~\ref{popul:extra1} and \ref{graph:perturb}.
	This completes the proof.
\end{proof}

Now, we prove our second main result.

\begin{proof}[Proof of Theorem~\ref{gauss:eigvec}]
	\eqref{gauss:eigvec-res2} is obtained from \eqref{gauss:eigvec-res1} and Lemma~\ref{smooth:extra4}.
	
	To show \eqref{gauss:eigvec-res1}, the expression inside the absolute value on the left-hand side of \eqref{gauss:eigvec-res1} can be decomposed as
	\begin{equation*}
		\sum_{k \in \cJ_j} \sum_{\ell \in \cJ_j} \langle \hvarphi_{k,n,\e}, \sn{\phi_{\ell,\e}} \rangle^2_\LPn - |\cJ_j|
		- (\P_n-\P) \bigg( \sum_{k \in \cJ_j} \phi_k^2 \bigg)
		= E_1 + E_2 + E_3,
	\end{equation*}
	where
	\begin{align*}
		E_1 &= (\P_n-\P) \bigg( \sum_{k \in \cJ_j} \phi_{k,\e}^2 - \sum_{k \in \cJ_j} \phi_k^2 \bigg), \\
		E_2 &= \sum_{k \in \cJ_j} \sum_{\ell \in \cJ_j} \langle \bvarphi_{k,n,\e}, \sn{\phi_{\ell,\e}} \rangle^2_\LPn 
		- |\cJ_j| - (\P_n-\P) \bigg( \sum_{k \in \cJ_j} \phi_{k,\e}^2 \bigg), \\
		E_3 &= \sum_{k \in \cJ_j} \sum_{\ell \in \cJ_j} \langle \hvarphi_{k,n,\e}, \sn{\phi_{\ell,\e}} \rangle^2_\LPn - \sum_{k \in \cJ_j} \sum_{\ell \in \cJ_j} \langle \bvarphi_{k,n,\e}, \sn{\phi_{\ell,\e}} \rangle^2_\LPn.
	\end{align*}
	We can bound $E_1$, $E_2$, $E_3$ using Lemma~\ref{smooth:extra3}, Theorem~\ref{popul:eigvec}, Proposition~\ref{graph:eigvec:aux1}, respectively.
	This yields \eqref{gauss:eigvec-res1}, completing the proof.
\end{proof}

\section*{Acknowledgement}
This material is based upon work supported by the National Science Foundation under Grant Number DMS-2515934.
We also note that this work was conducted while Eunseong Bae was a Ph.D. student at the University of California, Davis.

\appendix


\section{Heuristic derivation of first-order terms of the eigenvalues of graph Laplacians} \label{heuristic}

This section provides heuristic arguments for why the first-order terms of the eigenvalues of graph Laplacians given in \eqref{gauss:eigval-form2} emerge, and we compare them to the first-order terms in the unnormalized case studied in \cite{li2025central}. This provides insights into how the differences between the two cases arise. While we only consider eigenvalues, we expect a similar heuristic applies to the case of eigenvectors. Our approach relies on deriving influence functions of statistical functionals via finding G\^ateaux derivatives.
For further details on this approach, we refer the reader to Chapter 2.1 of \cite{hampel986robust}.

For the sake of exposition, we assume throughout this section that the index $k \in \N$ is fixed and all eigenvalues under consideration have trivial multiplicity. 
We also omit some calculations or details to keep the presentation concise.

Consider the eigenvalues of the graph Laplacians as statistical functionals, i.e., $\lambda_{k,\e} = \lambda_{k,\e} [\P] = \Teps[\P]$, and $\hlambda_{k,n,\e} = \Teps[\P_n]$.
The \emph{influence function}  $\IF_\e$ of this functional then is  the G\^ateaux derivative of $\Teps[\P]$: 
With $\delta_x$ denoting the Dirac measure in $x$, let $\P_t = t \delta_x + (1-t)\P = \P + t(\delta_x - \P), t \in [0,1]$ (note that $\P_t$ also depends on $x$).
Then $$\IF_\e(x) := \frac{d}{dt} \Teps [\P_t] \Big|_{t = 0^+}. $$
Von Mises calculus then suggests that
$$\hlambda_{k,n,\e} - \lambda_{k,\e} = \Teps[\P_n] - \Teps[\P] 
= (\P_n - \P) (\IF_\e) \hquad+\hquad \text{remainder}.$$
Our goal now is to show that for the symmetrically normalized graph Laplacian, we have
$$\IF_\e \stackrel{\e\to0}{\approx} \IF := - \Delta_\M (\phi_k^2),$$
so that $(\P_n - \P) (\IF)$ coincides with the first-order term stated in \eqref{gauss:eigval-form2}. 
The analogue in the unnormalized case (under uniform sampling) is: 
$$ \IF^\un_\e \stackrel{\e\to0}{\approx} \IF^\un := -\Delta_\M (\phi_k^2) + \lambda_k \big(\phi^2_k - 1\big)
= -\lambda_k \phi_k^2 + 2\|\nabla_\M \phi_k\|_{T\M}^2\ - \lambda_k$$
(cf. \cite{li2025central}).
Comparing the two influence functions $\IF$ and $\IF^\un$, we see that they differ by the additive term $\lambda_k (\phi^2_k - 1)$, and our heuristic will show how this term is canceled out by degree normalization.

\subsection{Influence functions in case of unnormalized graph Laplacians}

In the remaining parts of Appendix \ref{heuristic}, we use the short-hand notation $\mydot{\bfO} = \frac{d}{dt} \bfO[\P_t] \Big|_{t = 0^+}$ for an operator $\bfO$.

We start with the unnormalized case, and view the eigenvalue as a statistical functional $\F_\e^\un [\P]$, i.e.,
$$ \lambda^\un_{k,\e} 
= \Tuneps[\P] 
= \langle \bfL_\e^\un[\P] \phi^\un_{k,\e}[\P], \phi^\un_{k,\e} [\P] \rangle_{\LP},$$
where $(\lambda^\un_{k,\e},\phi^\un_{k,\e}) = (\lambda^\un_{k,\e}[\P],\phi^\un_{k,\e}[\P])$ is the $k$-th eigenpair of $\bfL^\un_\e = \e^{-1}( \bfD_\e [\P] - \bfA_\e [\P])$ with
$$ \bfD_\e [\P] f (z) = \int_\M \fg_\e (y,z) \dP (y) \cdot f(z), \quad\text{and}\quad
\bfA_\e [\P] (z) = \int_\M \fg_\e (y,z) f(y) \dP (y).$$
The corresponding influence function is given by
$$\IF^\un_\e = \mydot F^{\rm un}_\e  
= \frac{d}{dt}\Big\langle \bfL^\un_\e[\P_t]\phi^\un_{k,\e} [\P_t], \phi^\un_{k,\e} [\P_t] \Big\rangle_\LPt \, \bigg|_{t=0^+}. $$

We will show below that
\begin{align}
	\hspace*{3.5cm}\IF_\e^\un
	= \Big\langle  \mydot\bfL^\un_\e \phi^\un_{k,\e} [\P], \phi^\un_{k,\e} [\P] \Big\rangle_\LP
    = \Big\langle  \mydot\bfL^\un_\e \phi^\un_{k,\e}, \phi^\un_{k,\e} \Big\rangle_\LP,
    \label{heuristic:eq1}
\end{align}
meaning that the eigenfunction perturbation does not affect the eigenvalue perturbation.
Also note that we can use the dot product in $\LP$, rather than $\LPt$.

Assuming \eqref{heuristic:eq1} to hold, we now need to derive $\mydot\bfL_\e \phi^\un_{k,\e}$. 
%
%
From
$$ \mydot\bfD_\e f (z) = \fg_\e(x,z) f(z) - (\bfD_\e f)(z), \hquad\text{and}\hquad
\mydot\bfA_\e f (z) 
= \fg_\e(x,z) f(x) - (\bfA_\e f)(z), $$
we have
$$
	\e\; \Big\langle \mydot\bfL^\un_\e \phi^\un_{k,\e}, \phi^\un_{k,\e} \Big\rangle_\LP
	= \int_\M \fg_\e(x,z) \big(\phi^\un_{k,\e} (z) - \phi^\un_{k,\e} (x) \big) \phi^\un_{k,\e}(z) d\P(z) - \e \lambda^\un_{k,\e}.
$$
Using well-known properties of the Gaussian kernel with variance 2,  the Taylor expansion under the integral sign yields
\begin{align*}
	\int_\M \fg_\e(x,z) \big(\phi^\un_{k,\e} (z) - \phi^\un_{k,\e} (x) \big) \phi^\un_{k,\e}(z) d\P(z)
	\approx - \e \phi_{k,\e}^\un (x) \Delta_\M \phi_{k,\e}^\un (x)
	+ 2 \e \|\nabla_\M \phi_{k,\e}^\un (x) \|_{T_x \M}^2,
\end{align*}
so that
\begin{align*}
	\Big\langle \mydot \bfL^\un_\e \phi^\un_{k,\e}, \phi^\un_{k,\e} \Big\rangle_\LP
	\stackrel{\e\to0}{\approx} - \phi_{k,\e}^\un \Delta_\M \phi_{k,\e}^\un
	+2 \|\nabla_\M \phi_{k,\e}^\un \|_{T \M}^2 - \lambda_{k,\e}^\un.
\end{align*}

Finally, since $\bfL^\un_\e \to \Delta_\M$ as $\e \to 0$, we have $\lambda^\un_{k,\e} \to \lambda_k$, $\phi^\un_{k,\e} \to \phi_k$ and $\phi^\un_{k,\e} \Delta_\M \phi^\un_{k,\e} \to \phi_k\Delta_\M \phi_k = \lambda_k \phi_k^2$. 
Therefore, the influence function for the unnormalized graph Laplacian for small $\e$ becomes
\begin{align*}
	\IF^\un_\e 
	&= \Big\langle \mydot \bfL^\un_\e\phi^\un_{k,\e}, \phi^\un_{k,\e} \Big\rangle_\LP
	\stackrel{\e\to0}{\approx} -\lambda_k \phi_k^2 + 2\|\nabla_\M \phi_k\|_{T\M}^2\ - \lambda_k,
\end{align*}
which is the influence function in \cite{li2025central} (up to an unimportant sign change), if specified to the uniform distribution.

\subsection{Influence functions in case of symmetrically normalized graph Laplacians}

Here our statistical functional of interest is given by  
$$\lambda_{k,\e} = \Teps[\P] = \langle \bfL_\e [\P] \phi_{k,\e} [\P], \phi_{k,\e} [\P] \rangle_\LP,$$
where now $(\lambda_{k,\e},\phi_{k,\e})$ are $k$-th eigenpairs of the {\em symmetrically normalized} graph Laplacian $\bfL_\e = \bfD_\e^{-1/2} (\bfD_\e - \bfA_\e) \bfD^{-1/2}_\e$.
Similar to above, 
it can be shown that the corresponding influence function satisfies
\begin{equation}
	\IF_\e
	= \Big\langle \mydot\bfL_\e \phi_{k,\e} [\P], \phi_{k,\e} [\P] \Big\rangle_\LP
    = \Big\langle \mydot\bfL_\e \phi_{k,\e}, \phi_{k,\e} \Big\rangle_\LP.
	\label{heuristic:eq2}
\end{equation}

By definition of the 
multiplication operator $\bfD_\e[\P_t] f(z)$ 
we have $\mydot{\bfD}_\e f(z) := \big(\fg_\e(y,z) - \rho_\e(z)\big) f(z)$.
Also note that $\big(\bfD_\e[\P_t]\big)^{-\alpha} f(z) = \Big(\displaystyle{\int_\M} \fg_\e(y,z)d\P_t(y)\Big)^{-\alpha} f(z)$ for $\alpha > 0$, and
\begin{equation}
	\mydot{\bfD}^{-1/2}_\e f(z)
	= -\frac{1}{2} \bfD_\e^{-3/2} \mydot{\bfD}_\e f(z).
	\label{heuristic:eq3}
\end{equation}
We thus can argue as follows. Using the product rule on $\bfL_\e[\P_t] = \bfD_\e^{-1/2}[\P_t] \big(\bfD_\e[\P_t]
- \bfA_\e[\P_t] \big) \bfD_\e^{-1/2}[\P_t] $ in the right-hand side of \eqref{heuristic:eq2}, we obtain a sum of three dot-products
\begin{align*}
	\IF_\e &= \Big\langle \bfD_\e^{-1/2} \big( \mydot{\bfD}_\e - \mydot{\bfA}_\e \big) \bfD_\e^{-1/2} \phi_{k,\e}, \phi_{k,\e} \Big\rangle_\LP
	+ \Big\langle \mydot{\bfD}_\e^{-1/2} (\bfD_\e - \bfA_\e) \bfD_\e^{-1/2} \phi_{k,\e},\phi_{k,\e}\Big\rangle_\LP
	\\&\quad + \Big\langle \bfD_\e^{-1/2} (\bfD_\e - \bfA_\e) \mydot{\bfD}_\e^{-1/2}\phi_{k,\e},\phi_{k,\e}\Big\rangle_\LP
	=: \Romannum{1} + \Romannum{2} + \Romannum{3}.
\end{align*}
Using $\phi _{k,\e} = \bfD_\e^{1/2}\phi^\un_{k,\e}$ and \eqref{heuristic:eq3}, it can be shown that
$$ \Romannum{1} = \IF^\un_\e, \hquad\text{and}\hquad
\Romannum{2} = \Romannum{3} = -\frac{1}{2} \lambda_{k,\e} \Big\langle \mydot{\bfD}_\e \phi^\un_{k,\e}, \phi^\un_{k,\e} \Big\rangle_\LP,
$$
and consequently,
$$ \IF_\e = \IF^\un_\e - \lambda_{k,\e} \big\langle \mydot{\bfD}_\e \phi_{k,\e}^\un, \phi_{k,\e}^\un \big\rangle_\LP. $$
As we have derived $\IF_\e^\un$ above already, it remains to consider the second term. For small $\e$, we have
\begin{align*}
	\big\langle \mydot{\bfD}_\e \phi_{k,\e}^\un, \phi_{k,\e}^\un \big\rangle_\LP
	= \int_\M \mydot{\bfD}_\e (\phi_{k,\e}^\un (z))^2 d \P(z)
	\stackrel{\e \to 0}{\approx} (\phi_{k,\e}^\un)^2 - \int_\M (\phi_{k,\e}^2 (z))^2 d \P(z)
	\stackrel{\e \to 0}{\approx} \phi_k^2 - 1.
\end{align*}

Therefore, we conclude
$$ \IF_\e = \IF^\un_\e - \lambda_{k,\e} \big\langle \mydot{\bfD}_\e \phi_{k,\e}^\un, \phi_{k,\e}^\un \big\rangle_\LP
\stackrel{\e\to0}{\approx} \big(- \Delta_\M \phi_k^2 + \lambda_k (\phi^2_k - 1)\big) - \lambda_k (\phi^2_k - 1)
= - \Delta_\M \phi_k^2. $$
We clearly see how the 'degree correction' term comes in and modifies the influence function from the unnormalized case.

\subsection{Negligibility of eigenfunction perturbation}

Here we will show \eqref{heuristic:eq1} and \eqref{heuristic:eq2}:
\begin{align*}
	\frac{d}{dt} \Big\langle \bfL_\e[\P_t] \phi_{k,\e} [\P_t], \phi_{k,\e} [\P_t] \Big\rangle_\LPt \bigg|_{t=0^+}
	= \big\langle \mydot\bfL_\e \phi_{k,\e}, \phi_{k,\e} \big\rangle_\LP,
\end{align*}
where our proof works for both Laplacians, symmetrically normalized and unnormalized.
To simplify notation, we use the shorthand $\bfL_{\e,t}$, $\phi_{\e,t}$, $\phi_\e$ and $\dot{\phi}_\e$ to denote $\bfL_\e[\P_t]$, $\phi_{k,\e} [\P_t]$, $\phi_{k,\e}$ and $\frac{d}{dt} \phi_{\e,t} \big|_{t=0^+}$, respectively.

First observe that
\begin{align*}
	\frac{d}{dt} \big\langle \bfL_{\e,t} \phi_{\e,t}, \phi_{\e,t}\big\rangle_\LPt
	&= \frac{d}{dt} \Big\langle \bfL_{\e,t}\phi_{\e,t}, \phi_{\e,t}\Big\rangle_{L_2(\P)}
	+ \int_\M \big(\bfL_{\e,t}\phi_{\e,t} ) \phi_{\e,t}\,d(\delta_x - \P)
	\\&\quad + t \cdot \frac{d}{dt} \int_\M \big(\bfL_{\e,t}\phi_{\e,t} ) \phi_{\e,t}\,d(\delta_x - \P)
	=: \mydot{I}_1 (t) + \mydot{I}_2 (t) + \mydot{I}_3 (t).
\end{align*}
Obviously, $\mydot{I}_2(0)= \int_\M (\bfL_\e \phi_\e) \phi_e d(\delta_x - \P) $ and $\mydot{I}_3(0)= 0$.
To evaluate $\mydot{I}_1(0)$, note that for operators $\{H_t\}$ and functions $\{h_t\}$, the functional $\Phi(H_t,h_t) = H_t(h_t)$ satisfies $\frac{d}{dt}\Phi(H_t,h_t) \Big|_{t=0^+} = \Phi(\mydot{H},h_0) + \Phi(H_0,\mydot{h})$. Thus,
$$
\mydot I_{1}(0) 
= \big\langle \mydot{\bfL}_\e \phi_\e,\phi_\e\big\rangle_{L_2(\P)} 
+ \big\langle \bfL_\e\mydot{\phi}_\e,\phi_\e\big\rangle_{L_2(\P)}
+ \big\langle \bfL_\e\phi_\e,\mydot{\phi}_\e\big\rangle_{L_2(\P)}
= \big\langle \mydot{\bfL}_\e \phi_\e,\phi_\e\big\rangle_\LP 
+ 2 \big\langle \bfL_\e\mydot{\phi}_\e,\phi_\e\big\rangle_\LP,
$$
and consequently,
\begin{align*}
	\frac{d}{dt} \big\langle \bfL_{\e,t} \phi_{\e,t}, \phi_{\e,t}\big\rangle_\LPt \bigg|_{t=0^+}
	&= \mydot{I}_1(0) + \mydot{I}_2(0) + \dot{I}_3(0)
	\\&= \big\langle \mydot{\bfL}_{\e} \phi_\e, \phi_\e\big\rangle_\LP
	+ 2\big\langle \bfL_\e\mydot{\phi}_\e,\phi_\e\big\rangle_\LP
	+ \int_\M \big(\bfL_\e\phi_\e ) \phi_\e\, d(\delta_x - \P).	
\end{align*} 
Now, we need to show that the sum of the last two terms above vanishes. To this end, we use $\|\phi_{\e,t}\|_\LPt^2 = 1$ for all $t$. This gives
\begin{align*}
	0 = \frac{d}{dt} \|\phi_{\e,t}\|_\LPt^2 \Big|_{t=0^+}
	&= \frac{d}{dt} \Big[\int_\M \phi^2_{\e,t}d\P + t \int_\M \phi^2_{\e,t}d(\delta_x - \P)\Big] \,\bigg|_{t=0^+}
	= \int_\M 2 \phi_\e \mydot{\phi}_\e d\P  +  \int_\M \phi^2_\e d(\delta_x - \P).
\end{align*}
Multiplying by $\lambda_\e$, we have
\begin{align*}
	0 
	= \int_\M 2 \lambda_\e \phi_\e \mydot{\phi}_\e d\P
    + \int_\M \lambda_\e \phi^2_\e d(\delta_x - \P)
	= 2 \big\langle \bfL_\e\phi_\e, \mydot{\phi}_\e \big\rangle_\LP
    + \int_\M (\bfL_\e\phi_\e) \phi_\e d(\delta_x - \P),
\end{align*}
and this is what we wanted to show.

\section{Technical proofs} \label{gauss:tech}

\subsection{Proof of Lemma~\ref{smooth:KDE}} \label{smooth:KDE-proof}

We first present a key technical result for evaluating integrals over the manifold $\M$.
This result will be used to prove subsequent technical results.

For $x \in \M$, we let $T_x \M$ denote the tangent space of $\M$ at $x$.
Under the identification $T_x \M \cong \R^d$, we regard $\Exp_x$, the exponential map at $x$, as a map from $\R^d$ to $\M \subset \R^D$.
Using $\Exp_x$, we can transform integrals over $\M$ into integrals over $\R^d$ by a change of variables.
Moreover, the exponential decay of the Gaussian kernel enables the transformed integral to be concentrated in an open ball in $\R^d$ with diverging radius.

We formulate this as a lemma below.
For the detailed proof, see Remark~B.2 and the proof of Theorem~3.5 in \cite{bae2026kernel}.

\begin{lemma} \label{gauss:tech:integral}
	There exist $\M$-dependent constants $r$ and $C$ such that if $\e <r$, then any bounded measurable function $f:\M \to \R$ satisfies
	\begin{gather*}
		\sup_{x \in \M} \bigg| \int_\M \fg_\e(x,y) f(y) \dP(y) - \int_{B_{\e}} \tfg (z) f(\Exp_x (\sqrt{\e} z)) dz \bigg| \leq C \e \|f\|_\infty, \\
		\sup_{x \in \M} \bigg| \int_\M \fg_\e(x,y) f(y) \dP(y) - \int_{B_{\e}} \tfg (z) f(\Exp_x (\sqrt{\e} z)) (1+\e R(z)) dz \bigg| \leq C \e^2 \|f\|_\infty,
	\end{gather*}
	where $B_{\e}$ and $\tfg(z)$ are defined by
	\begin{align*}
		B_{\e} = \set{z \in \R^d : \|z\|_{\R^d} \le \nu_\e := \sqrt{(4d+8)\log(1/\e)}},
		\quad \tfg(z) = \frac{1}{(4 \pi)^{d/2}} \exp\Big(-\frac{\|z\|^2_{\R^d}}{4}\Big),
	\end{align*}
	and $R(z)$ satisfies $|R(z)| \le C \|z\|_{\R^d}^2.$
\end{lemma}

Note that $\nu_\e \to \infty$ but $\sqrt{\e} \nu_\e \to 0$ as $\e \to 0$.
The function $\tfg$ is the PDF of the $d$-dimensional normal distribution with covariance matrix $2 I_d$. Hence, $\tfg$ is symmetric and its moments admit explicit forms, e.g., 
\begin{equation*}
	\int_{\R^d} \tfg(z) dz = 1, \quad \int_{\R^d} \tfg(z) \, (z^\top A z) dz = 2 \, \tr(A) \quad\text{for symmetric matrix $A \in \R^{d \times d}$}.
\end{equation*}

Now, we prove Lemma~\ref{smooth:KDE} as follows.

\begin{proof}[Proof of Lemma~\ref{smooth:KDE}]
	Note that we can view $\rho_\e$ as a smoothed KDE for the uniform probability measure on $\M$.
	According to Section~6.1 of \cite{bae2026kernel}, the smoothed KDE and its derivatives satisfy uniform convergence. 
	Specifically, there exist $\M$-dependent constants $r_1$ and $C_1$ such that $\e < r_1$ implies
	\begin{equation}
		\|\rho_\e -1\|_\infty, 
		\|\nabla_\M \rho_\e\|_\infty, 
		\|\nabla_\M^2 \rho_\e\|_\infty \le C_1 \e,
		\label{smooth:rel2:KDE-eq1}
	\end{equation}
	where $\nabla_\M^2$ denotes the Hessian operator.
	Hence, the first inequality in \eqref{smooth:KDE-res} follows directly from \eqref{smooth:rel2:KDE-eq1}.
	
	For the second inequality in \eqref{smooth:KDE-res}, by plugging $f (\cdot) = (\rho_\e(x) - \rho_\e(\cdot))/\e$ into Lemma~\ref{gauss:tech:integral} and adjusting $r$ and $C$ appropriately, we have for $\e <r$,
	\begin{equation}
		\sup_{x \in \M} \bigg| \bfL_\e^\unnor \rho_\e (x) - \int_{B_\e} \tfg(z) \frac{\rho_\e (x) - \rho_\e(\Exp_x(\sqrt{\e}z))}{\e} (1+\e R(z)) dz \bigg| \le C_2 \e \|\rho_\e\|_\infty,
		\label{smooth:rel2:KDE-eq2}
	\end{equation}
	for some $\M$-constant $C_2$.
	Let $ \trho_{x,\e} (z) = (\rho_\e \circ \Exp_x )(z). $
	The Taylor expansion of $\trho_{x,\e}$ around $z=0$ gives
	\begin{equation*}
		\trho_{x,\e} (\sqrt{\e}z) -\trho_{x,\e} (0) 
		= \sqrt{\e} z^\top \, \nabla_{\R^d} \trho_{x,\e} (0) + \frac{\e}{2} z^\top \bigg[ \nabla^2_{\R^d} \trho_{x,\e} (z^*) \bigg] z.
	\end{equation*}
	Using \eqref{smooth:rel2:KDE-eq1}, it can be easily shown that $\|\nabla^2_{\R^d} \trho_{x,\e} (z^*)\|_\F$ is uniformly bounded over $x \in \M$ by $C_3 \e$ for some $\M$-constant $C_3$.
	Hence, together with the symmetry of $\tfg$, we obtain 
	\begin{equation}
		\sup_{x \in \M} \Bigg| \int_{B_\e} \tfg(z) \frac{\rho_\e (x) - \rho_\e(\Exp_x(\sqrt{\e}z))}{\e} (1+\e R(z)) dz \Bigg| \leq C_4 \e,
		\label{smooth:rel2:KDE-eq3}
	\end{equation}
	for some $\M$-constant $C_4$.
	Therefore, the second inequality in \eqref{smooth:KDE-res} follows by \eqref{smooth:rel2:KDE-eq2} and \eqref{smooth:rel2:KDE-eq3}.
	This completes the proof.
\end{proof}

\subsection{Proof of Theorem~\ref{smooth:eigfun2}} \label{smooth:eigfun2-proof}

As noted earlier, we follow the proof techniques used for Theorem~2.6 in \cite{calder2022lipschitz} and adapt it to our setting.
The key idea is to use the Lipschitz regularity estimate stated in Proposition~\ref{Lipschitz:main}, namely,
\begin{equation*}
	|f(x)-f(y)| \leq C (\|f\|_\infty+\|\bfL_\e f\|_\infty) \big( \dist_\M(x,y) + \e^{1/3} \log(1/\e) \big).
\end{equation*}
This proposition is a modified version of Theorem~2.11 of \cite{calder2022lipschitz}, where kernels with finite support are used to construct graph Laplacians, whereas we use Gaussian kernels with full support.

\subsubsection*{Lipschitz regularity: bounded-away case}

We first establish Lipschitz regularity for the difference $|f(x)-f(y)|$ when $\dist_\M (x,y)$ is bounded-away by $\sqrt{\e} \nu_\e$, stated in Proposition~\ref{Lipschitz:away}.
We note that this proposition is a modification of Theorem~8.1 in \cite{calder2022lipschitz}, adapted to our setting.
For the proof, we introduce some lemmas below.

Using the notions in Lemma~\ref{gauss:tech:integral}, we define the truncated Gaussian integral operator by
\begin{equation*}
	\tbfG_\e f (x) = \int_{B_\e} \tq_\e (z) f(\Exp_x(\sqrt{\e}z)) dz,
\end{equation*}
where $\tq_\e (z) = \tfg (z)/A_\e$ and $A_\e = \int_{B_\e} \tfg (z) dz$. Note that $\tq_\e$ can be viewed as a probability measure on $\R^d$ with support $B_\e$.

The lemma below provides how $\bfL_\e$ can be replaced in terms of $\tbfG_\e$.

\begin{lemma} \label{Lipschitz:approx}
	There exist $\M$-dependent constant $r$ and $C$ such that if $\e < r$,
	then any bounded measurable function $f: \M \to \R$ satisfies
	\begin{equation*}
		\Bigg\| \, \bfL_\e f - \frac{f-\tbfG_\e f}{\e} \, \Bigg\|_\infty 
		\leq C  \|f\|_\infty.
	\end{equation*}
\end{lemma}

\begin{proof}[Proof of Lemma~\ref{Lipschitz:approx}]
	Let $r$ as in Lemmas~\ref{smooth:KDE} and \ref{gauss:tech:integral}, and let $\bfG_\e$ be the (untruncated) Gaussian integral operator defined by
	\begin{equation*}
		\bfG_\e f (x) = \int_\M \fg_\e(x,y) f(y) \dP(y).
	\end{equation*}
	
	For a bounded measurable function $f:\M\to\R$, we first consider the decomposition
	\begin{align}
		\bfL_\e f - \frac{f - \tbfG_\e f}{\e} 
		= \frac{\bfW_\e f - \tbfG_\e f}{\e}
		= \frac{\bfW_\e f - \bfG_\e f}{\e} + \frac{\bfG_\e f - \tbfG_\e f}{\e}.
		\label{Lipschitz:approx-eq1}
	\end{align}
	
	The first term in \eqref{Lipschitz:approx-eq1} is expressed as
	\begin{align*}
		\bfW_\e f (x) - \bfG_\e f(x)
		= \int_\M \fg_\e(x,y) \bigg( \frac{1}{\sqrt{\rho_\e(x)}\sqrt{\rho_\e(y)}} -1 \bigg) f(y) \dP(y).
	\end{align*}
	Using $\|\rho_\e - 1\|_\infty < C_1 \e$ stated in Lemma~\ref{smooth:KDE} and adjusting $r$ if needed, it can be easily shown that 
	\begin{align*}
		\sup_{x,y \in \M} \bigg| \frac{1}{\sqrt{\rho_\e(x)}\sqrt{\rho_\e(y)}} -1 \bigg| \leq C_2 \e.
	\end{align*}
	for some $\M$-dependent constant $C_2$.
	Thus, again using $\|\rho_\e - 1\|_\infty < C_1 \e$, we obtain
	\begin{equation}
		\big\| \bfW_\e f - \bfG_\e f \big\|_\infty
		\leq C_2 \e \cdot \|f\|_\infty \e \|\rho_\e\|_\infty
		\leq C_3 \|f\|_\infty \e.
		\label{Lipschitz:approx-eq2}
	\end{equation}
	for some $\M$-dependent constant $C_3$.
	
	For the second term in \eqref{Lipschitz:approx-eq1}, we decompose
	\begin{equation}
		\bfG_\e f (x) - \tbfG_\e f(x)
		= \bfG_\e f (x) - \int_{B_\e} \tfg(z) f (\Exp_x(\sqrt{\e}z)) dz 
		+ \bigg(1 - \frac{1}{A_\e}\bigg) \int_{B_\e} \tfg(z) f (\Exp_x(\sqrt{\e}z)) dz.
		\label{Lipschitz:approx-eq3}
	\end{equation}
	By Lemma~\ref{gauss:tech:integral}, we have
	\begin{equation}
		\sup_{x \in \M} \bigg|\bfG_\e f (x) - \int_{B_\e} \tfg(z) f (\Exp_x(\sqrt{\e}z)) dz \bigg| 
		\leq C_4 \e \|f\|_\infty.
		\label{Lipschitz:approx-eq4}
	\end{equation}
	for some $\M$-dependent constant $C_4$.
	A simple calculation gives
	\begin{equation}
		|1-A_\e| \leq \frac{1}{(4\pi)^{d/2}} \exp \Big( -\frac{1}{4} \cdot (4d+8) \log(1/\e) \Big)
		= \frac{\e^{d+2}}{(4\pi)^{d/2}},
		\label{Lipschitz:approx-eq7}
	\end{equation}
	and consequently,
	\begin{equation}
		\sup_{x \in \M} \bigg| \bigg(1 - \frac{1}{A_\e}\bigg) \int_{B_\e} \tfg(z) f (\Exp_x(\sqrt{\e}z)) dz \bigg|
		\leq |1-A_\e| \cdot \|f\|_\infty
		\leq \frac{\e^{d+2} \|f\|_\infty}{(4\pi)^{d/2}}.
		\label{Lipschitz:approx-eq5}
	\end{equation}
	Combining \eqref{Lipschitz:approx-eq3}-\eqref{Lipschitz:approx-eq5} yields
	\begin{equation}
		\|\bfG_\e f - \tbfG_\e f \|_\infty \leq C_5 \e \|f\|_\infty,
		\label{Lipschitz:approx-eq6}
	\end{equation}
	for some $\M$-dependent constant $C_5$.
	
	Finally, combining \eqref{Lipschitz:approx-eq1}, \eqref{Lipschitz:approx-eq2} and \eqref{Lipschitz:approx-eq6}, the claim follows.
\end{proof}

The following lemma gives integral estimates for $\tq_\varepsilon$.

\begin{lemma} \label{Lipschitz:normalapprox}
	There exist $\M$-dependent constants $r$ and $C$ such that if $\e<r$ and $v \in \R^d$ is a unit vector, then
	\begin{equation*}
		\bigg| \int_{B_\e} \tq_\e(z) \langle z, v \rangle_{\R^d}^2 dz - 2 \bigg|, \hquad
		\bigg| \int_{B_\e} \tq_\e(z) \|z\|_{\R^d}^2 dz - 2d \bigg| \leq C \e^{d+2} \log(1/\e).
	\end{equation*}
\end{lemma}

\begin{proof}[Proof of \ref{Lipschitz:normalapprox}]
	We prove only the first inequality, and the second can be shown similarly.
	
	For the first inequality, we separate and bound
	\begin{equation}
		\bigg| \int_{B_\e} \tq_\e(z) \langle z, v \rangle_{\R^d}^2 dz - 2 \bigg|
		\leq 2 \bigg|\frac{1}{A_\e}-1\bigg|
		+ \frac{1}{A_\e} \bigg| \int_{B_\e} \tfg (z) \langle z, v \rangle_{\R^d}^2 dz - 2 \bigg| =: E_1 + E_2.
		\label{Lipschitz:normalapprox-eq1}
	\end{equation}
	Using \eqref{Lipschitz:approx-eq7}, we can control $A_\e$, and consequently, we obtain
	\begin{equation}
		E_1 \leq C_1 \e^{d+2},
		\label{Lipschitz:normalapprox-eq2}
	\end{equation}
	provided $\e<r_1$ for some $\M$-dependent constants $r$ and $C_1$.
	To bound $E_2$, noting $\tfg$ is the PDF of $d$-dim normal distribution with covariance matrix $2I_d$, we see that	
	\begin{equation*}
		\int_{\R^d} \tfg(z) \langle z, v \rangle_{\R^d}^2 dz = 2.
	\end{equation*}
	This gives
	\begin{equation*}
		E_2 
		= \frac{1}{A_\e} \int_{\R^d \backslash B_\e} \tfg (z) \langle z, v \rangle_{\R^d}^2 dz
		\leq \frac{1}{A_\e} \int_{\{z \in \R^d: \|z\|_{\R^d} > \nu_\e\}} \tfg (z) \|z\|_{\R^d}^2 dz.
	\end{equation*}
	By the tail properties of Gaussian kernels (see, e.g., Lemma~B.3 in \cite{bae2026kernel}), we bound
	\begin{equation*}
		\int_{\R^d \backslash B_\e} \tfg (z) \|z\|_{\R^d}^2 dz
		\leq C_2 \nu_\e^2 \exp\bigg( -\frac{1}{4} \nu_\e^2 \bigg)
		= C_2(4d+8) \e^{d+2} \log(1/\e),
	\end{equation*}
	for some constant $C_2$ depending only on $d$.
	Hence, again using \eqref{Lipschitz:approx-eq7}, we control $A_\e$ and obtain
	\begin{equation}
		E_2 \leq C_3 \e^{d+2} \log(1/\e),
		\label{Lipschitz:normalapprox-eq3}
	\end{equation}
	provided $\e < r_3$ for some $\M$-dependent constants $r_3$ and $C_3$.
	
	Therefore, combining \eqref{Lipschitz:normalapprox-eq1}-\eqref{Lipschitz:normalapprox-eq3} completes the proof.
\end{proof}

For $x,y \in \M$, we denote by 
\begin{equation*}
	\PT_{x}^{y}: T_x \M \to T_y \M
\end{equation*}
the parallel transport along the minimal geodesic connecting $x$ to $y$.
This operator allows us to move tangent vectors from one point on $\M$ to another while preserving their lengths and angles.
When $x$ and $y$ are sufficiently close (in particular, within the injectivity radius), $\PT_{x}^{y}$ is well-defined and is an isometry; i.e., for any $v,w \in T_x \M$,
$ \langle \PT_{x}^{y} (v), \PT_{x}^{y} (w) \rangle_{T_y \M} =  \langle v, w \rangle_{T_x \M}. $
In particular, the isometry property implies
$ \Jac(\PT_{x}^{y}) = 1, $
where $\Jac$ denotes the absolute value of the determinant of the corresponding Jacobian matrix.

The map $\HP_{x|y}: T_x \M \to T_x \M$ denotes the Householder projection with respect to $\Exp_x^{-1}(y)$.
It maps a tangent vector in $T_x \M$ to its reflection across the hyperplane perpendicular to $\Exp_x^{-1}(y)$.
As a formula, for $v \in T_x \M$,
\begin{equation*}
	\HP_{x|y} (v) 
	= v - \frac{2\langle v, \Exp_x^{-1}(y)\rangle_{T_x \M} }{\langle \Exp_x^{-1}(y), \Exp_x^{-1}(y)\rangle_{T_x \M}} \Exp_x^{-1} (y)
	= v - \frac{2\langle v, \Exp_x^{-1}(y)\rangle_{T_x \M} }{\dist^2_\M(x,y)} \Exp_x^{-1} (y).
\end{equation*}
Similar to parallel transport, $\HP_{x|y}$ is well-defined when $x$ and $y$ are sufficiently close and it is an isometry. 

Lastly, we define $\HPT_{x}^{y}: T_x \M \to T_y \M$ as the composition of $\HP_{x|y}$ and $\PT_{x}^{y}$, i.e., for $v \in T_x \M$,
\begin{equation*}
	\HPT_{x}^{y} (v) = (\PT_{x}^{y} \circ \HP_{x|y}) (v).
\end{equation*}
This combined map first reflects $v \in T_x \M$ with respect to $\Exp_x^{-1}(y)$ and then transports it to $T_y \M$.
It will be a key tool for comparing $\tbfG_\e f(x)$ and $\tbfG_\e f(y)$ for different points $x,y \in \M$.

The following lemma concerns the Levi-Civit\`a quadrilateral, 
which provides an estimate for the distance between two points perturbed along tangent vectors.

\begin{lemma} \label{Lipschitz:Levi}
	There exist $\M$-dependent constants $r$ and $C$ such that
	for $s>0$, $x,y \in \M$, $v \in T_x \M$, $w \in T_y \M$ satisfying
	$s, \dist_\M (x,y), \|sv\|_{T_x\M}, \|sw\|_{T_y\M} < r, $
	we have
	\begin{align}
		\dist_\M^2 ( \Exp_x (sv), \Exp_y (sw) ) - \dist_\M^2 (x,y)
		&= F_1 + E_1, 
		\label{Lipschitz:Levi-res1} \\
		\dist_\M^2 ( \Exp_x (sv), \Exp_y (s\HPT_x^y(v)) ) - \dist_\M^2 (x,y)
		&= F_2 + E_2,
		\label{Lipschitz:Levi-res2}
	\end{align}
	where
	\begin{align*}
		F_1 &= 2 s \langle \Exp_x^{-1}(y), \PT_{y}^{x} (w) -v \rangle_{T_x \M} + s^2 \| \PT_{y}^{x} (w) -v \|_{T_x \M}^2, \\
		F_2 &= -\frac{4 s \langle \Exp_x^{-1}(y), v \rangle_{T_x \M}}{\dist_\M(x,y)} + \frac{4 s^2 \langle \Exp_x^{-1}(y), v \rangle^2_{T_x \M}}{\dist^2_\M(x,y)}, \\
		|E_1|&\leq C s^2 \big( \dist_\M^2(x,y) + s^2(\|v\|_{T_x \M}^2 + \|w\|_{T_y \M}^2) \big), \\
		|E_2|&\leq C s^2 \big( \dist_\M^2(x,y) + s^2 \|v\|_{T_x \M}^2 \big).
	\end{align*}
	Additionally, if $\dist_\M(x, y) \ge \frac{3 s \langle \Exp_x^{-1}(y), v \rangle_{T_x \M}}{\dist_\M(x,y)}$, 
	we have
	\begin{equation}
		\dist_\M ( \Exp_x (sv), \Exp_y (s\HPT_x^y(v)) ) - \dist_\M (x,y)
		= -\frac{2s \langle \Exp_x^{-1}(y), v \rangle_{T_x \M}}{\dist_\M(x,y)} + E_3,
		\label{Lipschitz:Levi-res3}
	\end{equation}
	where $|E_3| \leq Cs^2 \Big( \dist_\M(x,y) + \frac{s^2 \|v\|^2_{T_x \M}}{\dist_\M(x,y)} \Big)$.
	In particular, if $\dist_\M(x, y) \ge 3 \|sv\|_{T_x \M}$, then \eqref{Lipschitz:Levi-res3} holds with
	$|E_3| \leq Cs^2 \dist_\M(x,y)$.
\end{lemma}

\begin{proof}[Proof of Lemma~\ref{Lipschitz:Levi}]
	\eqref{Lipschitz:Levi-res1} is a slight modification of Lemma~5.1 of \cite{calder2022lipschitz}; hence we omit the proof.
	
	\eqref{Lipschitz:Levi-res2} follows directly from \eqref{Lipschitz:Levi-res1} by replacing $w$ with $\HPT_x^y(v)$ and using $\PT_y^x \circ \HPT_x^y = \HP_{x|y}$.
	
	To show \eqref{Lipschitz:Levi-res3}, from the identity
	$ \sqrt{c+h} = \sqrt{c} + \frac{h}{\sqrt{c+h}+\sqrt{c}}, $
	substituting $c = \dist^2_\M(x, y) + F_2(x,y,s,v)$, $h = E_2(x,y,s,v)$ yields
	\begin{equation}
		\dist_\M ( \Exp_{x} (sv), \Exp_{y} (s \HPT_{x}^{y} (v))
		= \sqrt{F_3} + E_3,
		\label{Lipschitz:Levi-eq1}
	\end{equation}
	where
	$F_3 = \dist^2_\M(x,y) + F_2$ and $E_3 = \frac{E_2}{\sqrt{F_3 + E_2}+ \sqrt{F_3}}.$
	We see that
	\begin{align*}
		\dist^2_\M(x, y) + F_2(x,y,s,v)
		&= \dist^2_\M(x, y) -\frac{4 s \langle \Exp_x^{-1}(y), v \rangle_{T_x \M}}{\dist_\M(x,y)} + \frac{4 s^2 \langle \Exp_x^{-1}(y), v \rangle^2_{T_x \M}}{\dist^2_\M(x,y)}
		\\&= \bigg( \dist_\M(x, y) -\frac{2 s \langle \Exp_x^{-1}(y), v \rangle_{T_x \M}}{\dist_\M(x,y)} \bigg)^2,
	\end{align*}
	Since $\frac{2 s \langle \Exp_x^{-1}(y), v \rangle_{T_x \M}}{\dist_\M(x,y)}
	\le \frac{2}{3}\dist_\M(x,y),$ 
	we obtain
	\begin{equation}
		\sqrt{F_3} = \sqrt{\dist^2_\M(x, y) + F_2} = \dist_\M(x, y)-\frac{2 s \langle \Exp_x^{-1}(y), v \rangle_{T_x \M}}{\dist_\M(x,y)} \ge \frac{1}{3} \dist_\M (x,y),
		\label{Lipschitz:Levi-eq2}
	\end{equation}
	and
	\begin{equation}
		|E_3| \leq \frac{|E_2|}{\sqrt{F_3}} \leq 3Cs^2 \bigg( \dist_\M(x,y) + \frac{s^2 \|v\|^2_{T_x \M}}{\dist_\M(x,y)} \bigg).
		\label{Lipschitz:Levi-eq3}
	\end{equation}
	Therefore, \eqref{Lipschitz:Levi-res3} follows by combining \eqref{Lipschitz:Levi-eq1}-\eqref{Lipschitz:Levi-eq3} and adjusting $C$ appropriately.
\end{proof}


Equipped with these lemmas, we establish the Lipschitz regularity in the bounded-away case as follows.

\begin{proposition} \label{Lipschitz:away}
	There exist $\M$-dependent constants $r$ and $C$ such that if $\e < r$ and $x_0,y_0 \in \M$ with $ 3 \sqrt{\e} \nu_\e < \dist_\M (x_0,y_0) < r$, then any bounded measurable function $f:\M \to \R$ satisfies
	\begin{multline*}
		|f(x_0) - f(y_0) | \leq \sup \{ \, |f(x)-f(y)| : \dist_\M(x,y) \leq 3 \sqrt{\e}\nu_\e \, \}
		\\+ C \left( \|f\|_\infty + \|\bfL_\e f \|_\infty \right) \dist_\M(x_0,y_0).
	\end{multline*}
\end{proposition}

\begin{proof}[Proof of Proposition~\ref{Lipschitz:away}]
	We proceed in several steps.
	
	\mycase{1. Initialization}
	Let $r$ be as in both Lemmas~\ref{Lipschitz:approx}-\ref{Lipschitz:Levi}, and let $\e <r$ and $\N_0 = \N \cup \{0\}$.
	
	We define the random variables $\{\cX_m\}_{m\in\N_0}$ taking values in $\M$ via the recursive formula:
	\begin{align*}
		\cX_0 = x_0, \quad
		\cX_{m+1} &= \begin{cases}
			\Exp_{\cX_m} (\sqrt{\e} \, \PT_{x_0}^{\cX_m} (\cZ_{m+1}) ), &\text{if } \dist_\M (\cX_m, x_0) < r, \\
			\cX_m, &\text{otherwise.} 
		\end{cases}
	\end{align*}
	Here, $\{\cZ_m\}_{m\in\N}$ are i.i.d. random tangent vectors in $T_{x_0} \M$ distributed according to the probability density $\tq_\e$.
	
	Similarly, we also define the random variables $\{\cY_m\}_{m\in\N_0}$ taking values in $\M$ via the recursive formula:
	\begin{align*}
		\cY_0 = y_0, \quad
		\cY_{m+1} &= \begin{cases}
			\Exp_{\cY_m} (\sqrt{\e} (\HPT_{\cX_m}^{\cY_m} \circ \PT_{x_0}^{\cX_m}) (\cZ_{m+1})), 
			&\text{if } \dist_\M (\cX_m, x_0) < r \\
			& \text{and } 3 \sqrt{\e}\nu_\e < \dist_\M (\cX_m,\cY_m) < r, \\
			\cY_m, &\text{otherwise.}
		\end{cases}
	\end{align*}
	
	Lastly, we define the random variable $\cS$ by
	\begin{equation*}
		\cS = \min \{ m \in \N:\dist_\M (\cX_m, x_0) \geq r 
		\text{ or } \dist_\M (\cX_m,\cY_m) \leq 3 \sqrt{\e}\nu_\e 
		\text{ or } \dist_\M (\cX_m,\cY_m) \geq r \}.
	\end{equation*}
	By construction, $\cS$ is clearly a stopping time with respect to the filtration $\{\cF_m\}_{m \in \N_0}$, where $\cF_0$ is the trivial $\sigma$-algebra and $\cF_m$ is the $\sigma$-algebra generated by $\cZ_1, \dots, \cZ_m$ for $m\ge1$.
	Moreover, it satisfies $\tmbQ(\cS<\infty) = 1$, where $\tmbQ=\tmbQ_{x_0,\e}$ is the probability measure induced by $\{\cZ_m\}_{m \in \N}$.
	
	\mycase{2. Construction of a submartingale}
	Define the sequence of random variables $\{\cM_m\}_{m \in \N_0}$ by
	\begin{equation*}
		\cM_m = |f(\cX_{\cS \wedge m}) - f(\cY_{\cS \wedge m})| + 2 \|f-\tbfG_\e f\|_\infty (\cS \wedge m).
	\end{equation*}
	We claim that $\{\cM_m\}_{m \in \N_0}$ is a submartingale with respect to  $\{\cF_m\}_{m \in \N_0}$.
	
	On the event $\{\cS \leq m\}$, it is clear that
	$ \E [\cM_{m+1} \mid \cF_m] = \cM_m.$ 
	
	On the event $\{\cS > m\}$, note that
	\begin{align*}
		\E [ f(\cX_{m+1}) \mid \cF_m ]
		= \E [ f(\Exp_{\cX_m} (\sqrt{\e} \, \PT_{x_0}^{\cX_m} (\cZ_{m+1}) )) ]
		= \int_{B_\e} \tq_\e(z) f(\Exp_{\cX_m} (\sqrt{\e} \, \PT_{x_0}^{\cX_m} (z) )) dz.
	\end{align*}
	Since $\PT_{x_0}^{\cX_m}$ is an isometry, a change of variables yields
	\begin{align*}
		\int_{B_\e} \tq_\e (z) f(\Exp_{\cX_m} (\sqrt{\e} \, \PT_{x_0}^{\cX_m} (z) )) dz
		= \int_{B_\e} \tq_\e (z) f(\Exp_{\cX_m} (\sqrt{\e} \, z )) dz
		= \tbfG_\e f (\cX_m),
	\end{align*}
	which implies
	\begin{equation}
		\E [ f(\cX_{m+1}) \,\big|\, \cF_m ] = \tbfG_\e f (\cX_m).
		\label{smooth:lipschitz:away-eq1}
	\end{equation}
	Similarly, one obtains $\E [ f(\cY_{m+1}) \,\big|\, \cF_m ] = \tbfG_\e f (\cY_m)$.
	Using these identities, on the event $\{\cS > m\}$, we have
	\begin{align*}
		\E [ \, |f(\cX_{m+1}) - f(\cY_{m+1})| \mid \cF_m ]
		&\geq \left| \E [ f(\cX_{m+1}) - f(\cY_{m+1}) \mid \cF_m ] \right| 
		\\&= \left| \tbfG_\e f (\cX_m) - \tbfG_\e f (\cY_m) \right| 
		\\&\geq |f(\cX_{m}) - f(\cY_{m})| -\left| \tbfG_\e f (\cX_m) - f(\cX_m) \right| -\left| \tbfG_\e f (\cY_m) - f(\cY_m) \right|
		\\&\geq |f(\cX_{m}) - f(\cY_{m})| -2 \| \tbfG_\e f - f \|_\infty,
	\end{align*}
	and hence,
	\begin{align*}
		\E [ \cM_{m+1} | \cF_m ] - \cM_m 
		= \E [ \, |f(\cX_{m+1}) - f(\cY_{m+1})| \mid \cF_m ] - |f(\cX_{m}) - f(\cY_{m})| + 2 \| \tbfG_\e f - f \|_\infty
		\ge 0.
	\end{align*}
	Therefore, by combining the two cases, the claim follows.
	
	\mycase{3. Initial bound via the submartingale}
	Applying Doob’s optional stopping theorem to the submartingale $\{\cM_m\}_{m \in \N_0}$ yields
	\begin{align*}
		|f(x_0)-f(y_0)| =  \E \cM_0 \leq \E \cM_\cS 
		= \E|f(\cX_{\cS}) - f(\cY_{\cS})| + 2 \|\tbfG_\e f - f\|_\infty \E \cS.
	\end{align*}
	By Lemma~\ref{Lipschitz:approx}, we can choose an $\M$-dependent constant $C_1$ such that
	\begin{align*}
		2 \|\tbfG_\e f - f\|_\infty
		\leq C_1 \e \big( \|f\|_\infty + \|\bfL_\e^\unnor f\|_\infty \big).
	\end{align*}
	Next, we bound $E|f(\cX_{\cS}) - f(\cY_{\cS})|$ as
	\begin{multline*}
		\E|f(\cX_{\cS}) - f(\cY_{\cS})|
		\leq \int_{\{\dist_\M(\cX_\cS, \cY_\cS) \leq 3 \sqrt{\e}\nu_\e \,\}} |f(\cX_\cS)-f(\cY_\cS)| d \tmbQ
		+ 2 \|f\|_\infty \tmbQ(\dist_\M(\cX_\cS,\cY_\cS)>3 \sqrt{\e}\nu_\e\,).
	\end{multline*}
	Each term on the right-hand side can be bounded as
	\begin{align*}
		\int_{\{\dist_\M(\cX_\cS, \cY_\cS) \leq 3 \sqrt{\e}\nu_\e \}} |f(\cX_\cS)-f(\cY_\cS)| d \tmbQ
		&\leq \sup \{ |f(x)-f(y)|: x,y \in \M, \, \dist_\M(x,y) \leq 3 \sqrt{\e}\nu_\e \}, \\
		\tmbQ(\dist_\M(\cX_\cS,\cY_\cS)>3 \sqrt{\e}\nu_\e \,)
		&\leq \tmbQ(\dist_\M (\cX_\cS, x_0) \geq r) + \tmbQ(\dist_\M (\cX_\cS, \cY_\cS) \geq r),
	\end{align*}
	where the second inequality follows from the definition of $\cS$.
	
	Combining these bounds yields
	\begin{multline}
		|f(x_0)-f(y_0)| \leq \sup \{ |f(x)-f(y)|: x,y \in \M, \, \dist_\M(x,y) \leq 3 \sqrt{\e}\nu_\e \}
		\\ + 2 \|f\|_\infty \big( \tmbQ(\dist_\M (\cX_\cS, x_0) \geq r) + \tmbQ(\dist_\M (\cX_\cS, \cY_\cS) \geq r) \big)
		+ C_1 \big( \|f\|_\infty + \|\bfL_\e^\unnor f\|_\infty \big) \e \E \cS.
		\label{smooth:lipschitz:away-eq2}
	\end{multline}
	Hence, to complete the proof, it suffices to bound the following three quantities:
	\begin{equation*}
		\tmbQ(\dist_\M (\cX_\cS, x_0) \geq r), \quad \tmbQ(\dist_\M (\cX_\cS, \cY_\cS) \geq r), \quad \e \E \cS.
	\end{equation*}
	For this purpose, we introduce two supermartingales.
	
	\mycase{4. Construction of the first supermartingale}
	We show that the sequence of random variables $\{\cM'_m\}_{m\in\N_0}$, defined by
	\begin{equation*}
		\cM'_m = \dist^2_\M (\cX_{\cS \wedge m}, x_0) - 3d \e (\cS \wedge m),
	\end{equation*}
	forms a super-martingale with respect to $\{\cF_m\}_{m\in\N_0}$.
	
	On the event $\{\cS \leq m\}$, it is clear that
	$ \E [\cM'_{m+1} | \cF_m] = \cM'_m. $ 
	
	On the event $\{\cS > m\}$, similar to \eqref{smooth:lipschitz:away-eq1}, it can be shown
	\begin{align*}
		\E [ \dist^2_\M (\cX_{m+1}, x_0) \mid \cF_m ]
		&= \int_{B_\e} \tq_\e (z) \dist_\M^2 (\Exp_{\cX_m} (\sqrt{\e} z),x_0) dz,
	\end{align*}
	and this gives
	\begin{align*}
		\E [ \dist^2_\M (\cX_{m+1}, x_0) - \dist^2_\M (\cX_{m}, x_0) \mid \cF_m ]
		= \int_{B_\e} \tq_\e (z) \left[ \dist_\M^2 (\Exp_{\cX_m} (\sqrt{\e} z), x_0) - \dist_\M^2 (\cX_m,x_0) \right] dz.
	\end{align*}
	For $z \in B_\e$, applying Lemma~\ref{Lipschitz:Levi} with $x = \cX_m$, $y = x_0$, $s=\sqrt{\e}$, $v = z$, $w=0$ yields
	$$ \dist_\M^2 (\Exp_{\cX_m} (\sqrt{\e} z),x_0) 
	- \dist_\M^2 (\cX_m,x_0) = - 2\sqrt{\e} \langle \Exp_{\cX_m}^{-1} (x_0), z \rangle_{\R^d} + \e \|z\|^2_{\R^d} +  E'(z)$$
	with
	$ |E'(z)| \leq C'_1 \e ( \dist_\M^2 (\cX_m,x_0) + \e \|z\|_{\R^d}^2 ) $ for some $\M$-dependent constant $C_1'$.
	This yields
	\begin{align*}
		\E [ \, \dist^2_\M (\cX_{m+1}, x_0) - \dist^2_\M (\cX_{m}, x_0) \mid \cF_m ]
		= \e \int_{B_\e} \tq_\e (z) \|z\|^2_{\R^d} dz + \int_{B_\e} \tq_\e (z) E'(z) dz.
	\end{align*}
	By applying Lemma~\ref{Lipschitz:normalapprox} and adjusting $r$ if necessary, we obtain
	\begin{equation*}
		\Big| \E [ \, \dist^2_\M (\cX_{m+1}, x_0) - \dist^2_\M (\cX_{m}, x_0) \mid \cF_m ] \Big|	\leq 3d \e,	
	\end{equation*}	
	and it gives $\E [ \cM'_{m+1} | \cF_m ] \leq \cM'_m$.
	Therefore, by combining the two cases, the claim follows.
	
	\mycase{5. Construction of the second supermartingale}
	We show the sequence of random variables $\{\cM''\}_{m\in\N_0}$, defined by
	\begin{equation*}
		\cM''_m = \dist_\M (\cX_{\cS \wedge m}, \cY_{\cS \wedge m}) - \dist^2_\M (\cX_{\cS \wedge m}, \cY_{\cS \wedge m}) + 6 \e  (\cS \wedge m)
	\end{equation*}
	where $C''$ is $\M$-dependent, forms a super-martingale with respect to  $\{\cF_m\}_{m\in\N_0}$.
	
	On the event $\{\cS \leq m\}$, it is clear that
	$ \E [\cM''_{m+1} \mid \cF_m] = \cM''_m. $ 
	
	On the event $\{\cS > m\}$, similar to \eqref{smooth:lipschitz:away-eq1}, it can be shown
	\begin{multline*}
		\E [ \dist^2_\M (\cX_{m+1}, \cY_{m+1}) - \dist^2_\M (\cX_m, \cY_m) \mid \cF_m ]
		\\= \int_{B_\e} \tq_\e (z)  \left[ \dist^2_\M ( \Exp_{\cX_m} (\sqrt{\e} z), \Exp_{\cY_m} (\sqrt{\e}
		\HPT_{\cX_m}^{\cY_m} (z))) - \dist^2_\M (\cX_m, \cY_m) \right] dz.
	\end{multline*}
	Applying Lemma~\ref{Lipschitz:Levi} with $x=\cX_m$, $y=\cY_m$, $s=\sqrt{\e}$, $v=z$ yields
	\begin{equation*}
		\dist_\M^2(\Exp_{\cX_m} (\sqrt{\e}z), \Exp_{\cY_m} (\sqrt{\e} \, \HPT_{\cX_m}^{\cY_m}(z))) - \dist^2_\M(\cX_m, \cY_m)
		= F''_1(z) + E''_1(z),
	\end{equation*}
	where
	\begin{align*}
		F''_1(z) &= -\frac{4 \sqrt{\e} \langle \Exp_{\cX_m}^{-1}(\cY_m), z \rangle_{\R^d}}{\dist_\M(\cX_m,\cY_m)} + \frac{4 \e \langle \Exp_{\cX_m}^{-1}(\cY_m), z \rangle^2_{\R^d}}{\dist^2_\M(\cX_m,\cY_m)},
		\\ |E_1''(z)| &\leq C_1'' \e (\dist_\M^2(\cX_m, \cY_m) + \e\|z\|^2_{\R^d} ),
	\end{align*}
	for some $\M$-dependent constant $C_1''$.
	By symmetry, we obtain
	\begin{multline*}
		\E [ \dist^2_\M (\cX_{m+1}, \cY_{m+1}) - \dist^2_\M (\cX_m, \cY_m) \mid \cF_m ]
		\\= 4\e \int_{B_\e} \tq_\e (z) \bigg\langle z, \frac{\Exp_{\cX_m}^{-1} (\cY_m)}{\|\Exp_{\cX_m}^{-1} (\cY_m) \|_{\R^d}} \bigg\rangle_{\R^d}^2 dz + \int_{B_\e} \tq_\e(z) E_1''(z) dz.
	\end{multline*}
	By applying Lemma~\ref{Lipschitz:normalapprox} and adjusting $r$ if necessary, we obtain
	\begin{align}
		\E [ \dist^2_\M (\cX_{m+1}, \cY_{m+1}) - \dist^2_\M (\cX_m, \cY_m) \,|\, \cF_m ] \ge 7\e.
		\label{smooth:lipschitz:away-eq3}
	\end{align}
	
	Similarly, we have
	\begin{multline*}
		\E [ \dist_\M (\cX_{m+1}, \cY_{m+1}) - \dist_\M (\cX_m, \cY_m) \mid \cF_m ]
		\\= \int_{B_\e} \tq_\e (z)  \left[ \dist_\M ( \Exp_{\cX_m} (\sqrt{\e} z),
		\Exp_{\cY_m} (\sqrt{\e}
		\HPT_{\cX_m}^{\cY_m} (z))) - \dist_\M (\cX_m, \cY_m) \right] dz.
	\end{multline*}
	For $z\in B_\e$, since $\frac{3 \sqrt{\e} \langle \Exp_{\cX_m}^{-1}(\cY_m), z \rangle_{T_x \M}}{\dist_\M(\cX,\cY)} \leq 3\sqrt{\e}\|z\|_{\R^d} \leq 3\sqrt{\e} \nu_\e < \dist_\M(\cX_m, \cY_m)$,
	applying Lemma~\ref{Lipschitz:Levi} with $x=\cX_m$, $y=\cY_m$, $s=\sqrt{\e}$, $v=z$ gives
	\begin{equation*}
		\dist_\M(\Exp_{\cX_m} (\sqrt{\e}z), \Exp_{\cY_m} (\sqrt{\e} \, \HPT_{\cX_m}^{\cY_m}(z))) - \dist_\M(\cX_m, \cY_m)
		= -\frac{2\sqrt{\e} \langle \Exp_{\cX_m}^{-1}(\cY_m), z \rangle_{\R^d}}{\dist_\M(\cX_m,\cY_m)} + E''_2(z),
	\end{equation*}
	where $|E''_2 (z)| \leq C_2'' \e \Big( \dist_\M(\cX_m,\cY_m) + \frac{\e \|z\|^2_{\R^d}}{\dist_\M(\cX_m,\cY_m)} \Big) \leq \frac{9C_2''r\e}{10} $.
	By symmetry, we have
	\begin{align*}
		\E [ \dist_\M (\cX_{m+1}, \cY_{m+1}) - \dist_\M (\cX_m, \cY_m) \mid \cF_m ]
		= \int_{B_\e} E_2''(z) dz.
	\end{align*}
	By applying Lemma~\ref{Lipschitz:normalapprox} and adjusting $r$ if necessary, we obtain
	\begin{equation}
		\big| \E [ \dist_\M (\cX_{m+1}, \cY_{m+1}) - \dist_\M (\cX_m, \cY_m) \mid \cF_m ] \big| \leq \e.
		\label{smooth:lipschitz:away-eq4}
	\end{equation}
	Combining \eqref{smooth:lipschitz:away-eq3} and \eqref{smooth:lipschitz:away-eq4} yields
	\begin{align*}
		\E [ \cM'_{m+1} - \cM'_{m}| \cF_m ]
		&= \E [ \dist_\M (\cX_{m+1}, \cY_{m+1}) - \dist_\M (\cX_m, \cY_m) \mid \cF_m ]
		\\&\quad - \E [ \dist^2_\M (\cX_{m+1}, \cY_{m+1}) - \dist^2_\M (\cX_m, \cY_m) \mid \cF_m ]
		+6\e \leq 0,
	\end{align*}
	which confirms the claim.
	
	\mycase{6. Final bound via the supermartingales}
	Applying Doob’s optional stopping theorem to the supermartingales \\$\{\cM'_m\}_{m \in \N_0}$ and $\{\cM''_m\}_{m \in \N_0}$ yields
	\begin{align*}
		\E ( \dist^2_\M (\cX_{\cS}, x_0) ) - 3d \e  \E\cS
		&= \E \cM'_\cS \leq \E \cM'_0 = 0, \\
		\E ( \dist_\M (\cX_{\cS}, \cY_{\cS}) )
		- \E ( \dist_\M^2 (\cX_{\cS}, \cY_{\cS}) )
		+ 6 \e  \E \cS &= \E \cM''_\cS \leq \E \cM''_0 = \dist_\M (x_0, y_0) - \dist_\M^2 (x_0, y_0),
	\end{align*}
	and consequently,
	\begin{align}
		\E (\dist^2_\M (\cX_{\cS}, x_0) )
		&\leq 3d \e  \E \cS, 
		\label{smooth:lipschitz:away-eq5}
		\\ \E \left( \dist_\M (\cX_{\cS}, \cY_{\cS}) \right) - \E ( \dist_\M^2 (\cX_{\cS}, \cY_{\cS}) ) + 6 \e \E \cS
		&\leq \dist_\M (x_0, y_0).
		\label{smooth:lipschitz:away-eq6}
	\end{align}
	By construction, we assure that $\dist_\M (\cX_{\cS}, \cY_{\cS}) \leq 1/2$ by adjusting $r$ if necessary. This implies
	$ \dist_\M (\cX_{\cS}, \cY_{\cS}) - \dist^2_\M (\cX_{\cS}, \cY_{\cS})
	\ge \frac{1}{2} \dist_\M (\cX_{\cS}, \cY_{\cS}), $
	and applying this to \eqref{smooth:lipschitz:away-eq6} gives
	\begin{equation}
		\max\bigg\{ \frac{1}{2} \E \big( \dist_\M (\cX_{\cS}, \cY_{\cS}) \big) ,\, 6 \e \E \cS \bigg\}
		\leq \frac{1}{2} \E \left( \dist_\M (\cX_{\cS}, \cY_{\cS}) \right) + 6 \e \E \cS
		\leq \dist_\M (x_0, y_0).
		\label{smooth:lipschitz:away-eq7}
	\end{equation}	
	Combining \eqref{smooth:lipschitz:away-eq5} and \eqref{smooth:lipschitz:away-eq7} yields
	\begin{equation*}
		\E (\dist^2_\M (\cX_{\cS}, x_0) ) \leq \frac{d \cdot \dist_\M(x_0,y_0)}{2}, \quad 
		\E (\dist_\M (\cX_{\cS}, x_0) ) \leq 2 \dist_\M(x_0,y_0), \quad
		\e \E \cS \leq \frac{\dist_\M(x_0,y_0)}{6}.
	\end{equation*}
	By using Chebyshev's inequality, we conclude 
	\begin{equation*}
		\tmbQ(\dist_\M (\cX_\cS, x_0) \geq r), \,
		\tmbQ(\dist_\M (\cX_\cS, \cY_\cS) \geq r), \,
		\e \E \cS \leq C_1 \dist_\M (x_0,y_0),
	\end{equation*}
	for some $\M$-dependent constant $C_1$.
	
	Therefore, inserting them into \eqref{smooth:lipschitz:away-eq2} completes the proof.
\end{proof}

\subsubsection*{Lipschitz regularity: close case}

We establish Lipschitz regularity for the difference $|f(x)-f(y)|$ when $\dist_\M (x,y) \le 3 \sqrt{\e} \nu_\e$,
as stated in Proposition~\ref{Lipschitz:close}.
This proposition is adapted from Theorem~9.1 in \cite{calder2022lipschitz} to our setting, although the proof differs due to the different kernel.
In fact, our proof relies on semigroup arguments; see, e.g., Chapter~2 of \cite{engel2000one}. 

We begin by stating auxiliary lemmas.

\begin{lemma} \label{Lipschitz-aux1}
	There exist $\M$-dependent constants $r$, and $C$ such that if $\e < r$, $x,y \in \M$ and $k \in \N$, we have
	\begin{equation*}
		\int_\M |\fw_\e^{*k} (x,z) -\fw_\e^{*k} (y,z)| \dP(z) 
		\leq C a(k,\e) \dist_\M(x,y)  + 2 \Big[ \{ 1+C \e\log^2(1/\e) \}^k - 1 \Big],
	\end{equation*}
	where $\fw_\e^{*k}$ denotes the $k$-fold convolution of $\fw_\e$, i.e.,
	\begin{equation*}
		\fw_\e^{*k} (x,z) = \int_\M \dots \int_\M \fw_\e(x,z_1) \fw_\e(z_1,z_2) \cdots \fw_\e (z_{k-1},z) d \P(z_1) \cdots d \P(z_{k-1}),
	\end{equation*}
	and $a(k,\e) = \frac{1}{\sqrt{k\e}}$ if $k\e < 1$ and $a(k,\e)=1$ otherwise.
\end{lemma}

\begin{proof}[Proof of Lemma~\ref{Lipschitz-aux1}]
	Without loss of generality, we may assume $x$ and $y$ lie in the same connected component of $\M$, so that $\dist_\M(x,y)$ is finite.
	
	We proceed with the proof in several steps
	
	\mycase{1. Approximating $\fw_\e^{*k} (x,z)$}
	We use a similar iterative argument as in \cite{shamarova2017uniform}.
	Using $\|\rho_\e -1 \|_\infty = O_\M(\e)$ in Lemma~\ref{smooth:KDE} and the Gaussian approximation of the heat kernel formulated in \eqref{smooth:rel1-gauss}, we have
	\begin{equation*}
		\fw_\e (x,y) = \fh_\e (x,y) + A_\e^{(1)} (x, y),
		\quad |A_\e^{(1)} (x, y)| \leq C \e \log^2(1/\e) (\fh_\e (x,y) + \e).
	\end{equation*}
	This gives
	\begin{align*}
		\fw_\e^{*2} (x,y)
		&= \int_\M \fw_\e (x,z) \fw_\e (z,y) \dP(z)
		\\&= \int_\M (\fh_\e (x,z) + A_\e^{(1)} (x, z)) (\fh_\e (z,y) + A_\e^{(1)} (z, y)) \dP(z)
		\\&= \int_\M \fh_\e (x,z) \fh_\e (z,y) \dP(z) + \int_\M \fh_\e (x,z) A_\e^{(1)} (z, y) \dP(z)
		\\&\quad + \int_\M A_\e^{(1)} (x, z) \fh_\e (z,y) \dP(z) + \int_\M A_\e^{(1)} (x, z) A_\e^{(1)} (z, y) \dP(z).
	\end{align*}
	Note that for any $t,s >0$ and $x,y \in \M$, the heat kernel $\fh_\e$ satisfies
	\begin{equation*}
		\int_\M \fh_t (x,z) \fh_s (z,y) \dP(z) = \fh_{t+s} (x,y), \quad
		\int_\M \fh_s (x,z) \dP(z) = 1.
	\end{equation*}
	These give
	\begin{align*}
		\int_\M \fh_\e (x,z) \fh_\e (z,y) \dP(z) &= \fh_{2\e} (x,y),
		\\
		\bigg| \int_\M \fh_\e (x,z) A_\e^{(1)} (z, y) \dP(z) \bigg|
		&\leq C \e \log^2(1/\e) \int_\M \fh_\e (x,z) (\fh_\e (z, y) + \e) \dP(z)
		\\&= C \e \log^2(1/\e) (\fh_{2\e} (x,y) + \e),
		\\
		\bigg| \int_\M A_\e^{(1)} (x, z) A_\e^{(1)} (z, y) \dP(z) \bigg|
		&\leq (C\e \log^2(1/\e))^2 \int_\M (\fh_\e (x, z) + \e) (\fh_\e (z, y) + \e) \dP(z)
		\\&= (C\e \log^2(1/\e))^2 (\fh_{2\e} (x,y) + 2\e + \e^2).
	\end{align*}
	Consequently, we obtain
	\begin{equation*}
		\fw^{*2}_\e (x,y) = \fh_{2\e} (x,y) + A_\e^{(2)} (x, y),
	\end{equation*}
	where
	\begin{equation*}
		|A_\e^{(2)} (x, y)| \leq 
		2C \e \log^2(1/\e) (\fh_{2\e} (x,y) + \e) + (C\e \log^2(1/\e))^2 (\fh_{2\e} (x,y) + 2\e + \e^2).
	\end{equation*}
	Iterating this argument $k$ times yields
	\begin{equation}
		\fw^{*k}_\e (x,y) = \fh_{k\e} (x,y) + A_\e^{(k)} (x, y),
		\label{Lipschitz-aux1-eq1}
	\end{equation}
	where
	\begin{equation*}
		| A_\e^{(k)} (x,y) | \leq \sum_{\ell=1}^k \binom{k}{\ell} (C \e \log^2(1/\e))^\ell 
		\bigg( \fh_{k\e}(x,y) + \sum_{m=1}^\ell \binom{\ell}{m} \e^m \bigg) .
	\end{equation*}
	
	\mycase{2. Concluding the proof}
	Using \eqref{Lipschitz-aux1-eq1}, we obtain
	\begin{multline}
		\int_\M |\fw_\e^{*k} (x,z) -\fw_\e^{*k} (y,z)| \dP(z) 
		\leq \int_\M |\fh_{k\e} (x,z) - \fh_{k\e} (y,z)| \dP (z)
		\\ \quad + 2 \sum_{\ell=1}^k \binom{k}{\ell} (C \e \log^2(1/\e))^\ell 
		\bigg(1 + \sum_{m=1}^\ell \binom{\ell}{m} \e^m \bigg).
		\label{Lipschitz-aux1-eq2}
	\end{multline}
	
	For the second term, simple algebra gives
	\begin{align}
		2 \sum_{\ell=1}^k \binom{k}{\ell} (C \e \log^2(1/\e))^\ell 
		\bigg(1 + \sum_{m=1}^\ell \binom{\ell}{m} \e^m \bigg) 
		&= 2 \Big[ \{ 1+C \e (1+\e)\log^2(1/\e) \}^k - 1 \Big]
		\notag\\&\leq 2 \Big[ \{ 1+2C \e \log^2(1/\e) \}^k - 1 \Big].
		\label{Lipschitz-aux1-eq3}
	\end{align}
	
	For the first term, let $\eta: [0,1] \to \M$ be the minimal geodesic from $y$ to $x$, parameterized with constant speed, i.e., $\eta(0)=y$, $\eta(1)=x$ and $\|\eta'(t)\|_{T_{\eta(t)}\M} = \dist_\M(x,y)$.
	Using $\eta$, we have
	\begin{align*}
		\fh_{k\e} (x,z) - \fh_{k\e} (y,z) 
		= \int_0^1 \frac{d}{dt} \fh_{k\e} (\eta(t),z) dt
		= \int_0^1 \langle \nabla_{\M,1} \fh_{k\e} (\eta(t),z), \eta'(t) \rangle_{T_{\eta(t)}\M} dt,
	\end{align*}
	where $\nabla_{\M,1}$ denotes the covariant derivative with respect to the first argument.
	This yields
	\begin{equation}
		\int_\M |\fh_{k\e} (x,z) - \fh_{k\e} (y,z)| d \P(z)
		\leq \dist_\M(x,y) \int_0^1  \int_\M \big\| \nabla_{\M,1} \fh_{k\e} (\eta(t),z) \big\|_{T_{\eta(t)} \M } d \P(z) dt.
		\label{Lipschitz-aux1-eq4}
	\end{equation}
	
	According to the Gaussian estimate of the derivative of heat kernels (see Theorem~6 of \cite{davies1989pointwise} and Theorem~5.5.3 of \cite{hsu2002stochastic}), it holds that
	\begin{equation}
		\big\| \nabla_{\M,1} \fh_{k\e} (\eta(t),z) \big\|_{T_{\eta(t)} \M } 
		\leq 
		\begin{cases}
			C \bigg( \frac{\dist_\M(\eta(t),z)}{k\e} +\frac{1}{\sqrt{k\e}} \bigg) \fh_{k\e} (\eta(t),z), & \text{if } k\e < 1,\\[1em]
			C \exp \Big( -\frac{\dist_\M^2(\eta(t),z)}{5k\e} \Big)\leq C, & \text{if } k\e \ge 1.
		\end{cases}
		\label{Lipschitz-aux1-eq6}
	\end{equation}
	By slightly modifying Lemma~5.4.2 of \cite{hsu2002stochastic}, it can be shown
	\begin{align*}
		\int_\M \dist_\M(\eta(t),z) \, h_{k\e} (\eta(t),z) \dP(z) \leq C \sqrt{k\e},
	\end{align*}
	and together with \eqref{Lipschitz-aux1-eq6} and $\int_\M \fh_{k\e} (\eta(t),z) \dP(z) = 1$, it follows that
	\begin{equation*}
		\int_\M \big\| \nabla_{\M,1} \fh_{k\e} (\eta(t),z) \big\|_{T_{\eta(t)} \M } d\P(z)
		\leq C a(k,\e).
	\end{equation*}
	Thus, \eqref{Lipschitz-aux1-eq4} gives
	\begin{equation}
		\int_\M |\fh_{k\e} (x,z) - \fh_{k\e} (y,z)| d \P(z)
		\leq C a(k,\e) \dist_\M(x,y).
		\label{Lipschitz-aux1-eq5}
	\end{equation}
	Finally, by combining \eqref{Lipschitz-aux1-eq2}, \eqref{Lipschitz-aux1-eq3} and \eqref{Lipschitz-aux1-eq5}, we complete the proof.
\end{proof}

\begin{lemma} \label{Lipschitz-aux2}
	For $t>0$, we have
	\begin{equation}
		\sum_{k=1}^\infty \frac{t^k}{k! \sqrt{k}} \leq \frac{4e^t}{\sqrt{t}}, \quad
		\label{Lipschitz-aux2-res1}
	\end{equation}
	and for $0< t \e < 1$,
	\begin{equation}
		\sum_{k:k\e>1} \frac{t^k}{k!} \leq (3t\e)^{1/\e}.
		\label{Lipschitz-aux2-res2}
	\end{equation}
\end{lemma}

\begin{proof}[Proof of Lemma~\ref{Lipschitz-aux2}]
	By letting $Y$ be a Poisson random variable with intensity $t$, we have
	\begin{align}
		\sum_{k=1}^\infty \frac{t^k}{k! \sqrt{k}} 
		= e^t \E \big[ Y^{-1/2} \chi_{\{Y>0\}} \big]
		= e^t \Big( \E \big[Y^{-1/2}\chi_{\{Y>t/2\}}\big] + \E \big[Y^{-1/2}\chi_{\{0<Y\le t/2\}}\big]\Big),
		\label{Lipschitz-aux2-eq1}
	\end{align}
	where $\chi$ denotes the indicator function.
	From $Y^{-1/2}\chi_{\{Y > t/2\}}\leq 2/\sqrt{t}$, we obtain 
	\begin{equation}
		\E \big[Y^{-1/2}\chi_{\{Y > t/2\}}\big] \leq \frac{2}{\sqrt{t}}.
		\label{Lipschitz-aux2-eq2}
	\end{equation}	
	Next, from $Y\geq1$ for $Y>0$, we have $\E \big[Y^{-1/2}\chi_{\{0<Y\le t/2\}}\big] \leq \P( Y \leq t/2)$.
	Applying the Chernoff bound yields
	\begin{equation}
		\E \big[Y^{-1/2}\chi_{\{0<Y\le t/2\}}\big]
		\leq \P( Y \leq t/2) \leq \bigg( \frac{2}{e} \bigg)^{t/2} \leq \frac{2}{\sqrt{t}}.
		\label{Lipschitz-aux2-eq3}
	\end{equation}	
	Hence, combining \eqref{Lipschitz-aux2-eq1}-\eqref{Lipschitz-aux2-eq3} yields the first claim \eqref{Lipschitz-aux2-res1}.
	
	To show the second claim \eqref{Lipschitz-aux2-res2}, using $Y$ and the Chernoff bound, we have
	\begin{equation*}
		\sum_{k:k\e>1} \frac{t^k}{k!} 
		= e^t \P ( Y > 1/\e)
		\leq e^t \exp\bigg( \frac{1}{\e} -t -\frac{1}{\e} \log \Big(\frac{1}{t\e}\Big) \bigg)
		= (e t\e)^{1/\e} \leq (3 t \e)^{1/\e}.
	\end{equation*}
	This completes the proof.
\end{proof}

Now, we establish the Lipschitz regularity for the close case as follows.

\begin{proposition} \label{Lipschitz:close}
	There exist $\M$-dependent constants $r$ and $C$ such that if $\e <r$ and $x,y \in \M$ with $\dist_\M (x,y) \leq 3 \sqrt{\e} \nu_\e$, then any bounded measurable function $f: \M \to \R$ satisfies
	\begin{equation*}
		|f(x)-f(y)| \leq C \big(\|f\|_\infty + \|\bfL_\e f \|_\infty \big) \e^{1/3} \log(1/\e).
	\end{equation*}
\end{proposition}

\begin{proof}[Proof of Proposition~\ref{Lipschitz:close}]
	We divide the proof into several steps.
	
	\mycase{1. Initial bound via semigroup associated with $\bfW_\e$}
	Let $\{ \bfQ_{t,\e} = e^{t \bfW_\e}\}_{t\ge0}$ denote the semigroup generated by $\bfW_\e$.
	Note that $\bfQ_{t,\e}:\LP \to \LP$ is well-defined and $t$ is independent of $\e$.
	Let $\fq_{t,\e}(x,y)$ be the kernel corresponding to $\bfQ_{t,\e}$, that is,
	\begin{equation*}
		\bfQ_{t,\e} f (x) = \int_\M \fq_{t,\e} (x,z) f(z) d \P (z).
	\end{equation*}
	From  $ \frac{d}{dt} \bfQ_{t,\e} = \bfQ_{t,\e} \bfW_\e $ and $\bfI_\LP - \bfW_\e = \e \bfL_\e$, 
	we have
	\begin{align*}
		\frac{d}{dt} (e^{-t} \bfQ_{t,\e}f)
		= - e^{-t} \bfQ_{t,\e} f + e^{-t} \bfQ_{t,\e} \bfW_\e f
		= -e^{-t} \bfQ_{t,\e} (f - \bfW_\e f)
		= -\e e^{-t} \bfQ_{t,\e} \bfL_\e f
	\end{align*}
	Integrating from $0$ to $t$ gives
	$$ e^{-t}\bfQ_{t,\e} f - f
	= - \e \int_0^t e^{-s} \bfQ_{s,\e} \bfL_\e f \, ds,$$
	and consequently,
	\begin{equation*}
		f = e^{-t}  \bfQ_{t,\e} f + \e \int_0^t e^{-s} \bfQ_{s,\e} \bfL_\e f \, ds.
	\end{equation*}
	Using $\fq_{t,\e}$, we rewrite the above expression as
	\begin{align*}
		f(x) = e^{-t} \int_\M \fq_{t,\e} (x,z) f(z) \dP(z)
		+ \e \int_0^t e^{-s} \int_\M \fq_{s,\e} (x,z) \bfL_\e f (z) \dP(z) ds,
	\end{align*}
	and it follows that 		
	\begin{equation}
		|f(x)-f(y)| \leq \|f\|_\infty g_t(x,y) + \e \|\bfL_\e f \|_\infty \int_0^t g_s(x,y) ds,
		\label{Lipscthiz:close-eq1}
	\end{equation}
	where $g_t (x,y) = e^{-t} \int_\M |\fq_{t,\e} (x,z) - \fq_{t,\e} (y,z)| \dP(z) $.
	Therefore, to complete the proof, it suffices to control $g_t (x,y)$ and $\int_0^t g_s (x,y) ds$ properly.
	
	\mycase{2. Bounding $g_t(x,y)$ and $\int_0^t g_s (x,y) ds$}
	From the expansion
	\begin{equation*}
		\bfQ_{t,\e} 
		= e^{t \bfW_\e} 
		= \sum_{k=0}^\infty \frac{(t \bfW_\e)^k}{k!}
		= \bfI_\LP + \sum_{k=1}^\infty \frac{t^k}{k!} \bfW_\e^k,
	\end{equation*}
	we see that $\fq_{t,\e}$ can be expressed by
	\begin{equation*}
		\fq_{t,\e}(x,y) = \mathrm{Dirac}_x(y) + \sum_{k=1}^{\infty} \frac{t^k}{k!} \fw_\e^{*k}(x,y),
	\end{equation*}
	where $\mathrm{Dirac}_x$ denotes the Dirac delta function at $x$. 
	Using this and applying Lemma~\ref{Lipschitz-aux1}, we obtain
	\begin{multline}
		\int_\M |\fq_{t,\e} (x,z) - \fq_{t,\e} (y,z)| d \P(z)
		\leq 2 + \sum_{k=1}^{\infty} \frac{t^k}{k!} 
		\int_\M |\fw_\e^{*k}(x,z) - \fw_\e^{*k}(x,z)| \dP(z)
		\\\leq 2 + C \dist_\M(x,y) \sum_{k: k\e <1} \frac{t^k}{k!} \cdot \frac{1}{\sqrt{k\e}}
		+ C \dist_\M(x,y) \sum_{k: k \e \ge 1 } \frac{t^k}{k!}
		+ 2 \sum_{k=1}^\infty \frac{t^k}{k!} \cdot 2 \Big[ \{ 1+2C \e \log^2(1/\e) \}^k - 1 \Big],
		\label{Lipscthiz:close-eq2}
	\end{multline}
	provided $\e <r$, where $r$ and $C$ are $\M$-dependent.
	
	The last term in \eqref{Lipscthiz:close-eq2} is simplified to
	\begin{equation}
		\sum_{k=1}^{\infty} \frac{t^k}{k!} 2 \Big[ \{ 1+2C \e\log^2(1/\e) \}^k - 1 \Big]
		= 2 (e^{t \{ 1+2C \e \log^2(1/\e) \} }-e^{t})
		= 2 e^{t} (e^{2 C t\e \log^2(1/\e) }-1).
		\label{Lipscthiz:close-eq3}
	\end{equation}
	For the first two terms \eqref{Lipscthiz:close-eq2}, applying Lemma~\ref{Lipschitz-aux2} yields
	\begin{equation}
		\sum_{k:k\e <1}^{\infty} \frac{t^k}{k!} \cdot \frac{1}{\sqrt{k\e}}
		\leq \sum_{k=1}^{\infty} \frac{t^k}{k!} \cdot \frac{1}{\sqrt{k\e}}
		\leq \frac{4e^t}{\sqrt{t\e}}, \quad\quad
		\sum_{k: k\e \ge 1} \frac{t^k}{k!} \leq (3 t\e)^{1/\e},
		\label{Lipscthiz:close-eq5}
	\end{equation}
	provided $ 0<t \e <1$.
	
	Combining \eqref{Lipscthiz:close-eq2}-\eqref{Lipscthiz:close-eq5}, we bound $g_t(x,y)$ for $0<t \e<1$ as
	\begin{equation}
		g_t(x,y) 
		\leq 2 e^{-t} + \frac{4C \dist_\M (x,y)}{\sqrt{t\e}} + e^{-t} (3t\e)^{1/\e}
		+ 2 (e^{2C t \e\log^2(1/\e) }-1).
		\label{Lipscthiz:close-eq6}
	\end{equation}
	and consequently,
	\begin{equation}
		\int_0^t g_s (x,y) ds
		\leq 2(1-e^{-t}) + \frac{8C\sqrt{t} \, \dist_\M (x,y)}{\sqrt{\e}} + t e^{-t} (3t\e)^{1/\e} + 2 t(e^{2C t\e \log^2(1/\e) }-1).
		\label{Lipscthiz:close-eq7}
	\end{equation}
	
	\mycase{3. Concluding the proof}
	By combining \eqref{Lipscthiz:close-eq1}, \eqref{Lipscthiz:close-eq6} and \eqref{Lipscthiz:close-eq7}, we have for $0<t\e<1$ and $\e <r$,
	\begin{multline*}
		|f(x)-f(y)| \leq \|f\|_\infty \Big( 2 e^{-t} + \frac{4C \dist_\M (x,y)}{\sqrt{t\e}} 
		+ e^{-t} (3t\e)^{1/\e} + 2 (e^{2C t \e\log^2(1/\e) }-1). \Big)
		\\+ \|\bfL_\e f \|_\infty \Big( 2\e (1-e^{-t}) + 8C\sqrt{t\e} \, \dist_\M (x,y) + t\e e^{-t} (3t\e)^{1/\e} + 2 t\e(e^{2C t\e \log^2(1/\e) }-1) \Big).
	\end{multline*}
	Under $\dist_\M(x,y) \leq 3\sqrt{\e}\nu_\e = 3 \sqrt{4d+8} \sqrt{\e \log(1/\e)}$,
	we complete the proof by choosing the optimal $t = \e^{-2/3} \log^{-1}(1/\e) $	and adjusting $r$ and $C$ appropriately.
\end{proof}

\subsubsection*{Lipschitz regularity and concluding the proof}

We finalize the proof of Lipschitz regularity, and we use this to complete the proof of Theorem~\ref{smooth:eigfun2}.

\begin{proposition} \label{Lipschitz:main}
	There exist $\M$-dependent constants $r$ and $C$ such that if $\e <r$,
	then any bounded measurable function $f: \M \to \R$ and any $x,y \in \M$ satisfies
	\begin{equation} 
		|f(x)-f(y)| \leq C \big(\|f\|_\infty + \|\bfL_\e f \|_\infty \big) \big( \dist_\M(x,y) 
		+ \e^{1/3} \log(1/\e) \big).
		\label{Lipschitz:main-res}
	\end{equation}
\end{proposition}

\begin{proof}[Proof of Proposition~\ref{Lipschitz:main}]
	Combining Propositions~\ref{Lipschitz:away} and \ref{Lipschitz:close}, we obtain \eqref{Lipschitz:main-res} for $\e < r$ and $\dist_\M(x,y) <r $. 
	Extending to arbitrary $x,y \in \M$ follows directly from the compactness of $\M$.
	This completes the proof.
\end{proof}

We present and prove two auxiliary lemmas below.

\begin{lemma} \label{smooth:eigfun2-aux1}
	There exist $\M$-dependent constants $C$ and $r$ such that for any $t >0$ and any bounded measurable function $f:\M \to \R$,
	if $\e < \frac{r}{t+1}$ and $\frac{\|\bfL_\e f \|_\infty}{\|f\|_\infty} \le t$, we have
	\begin{equation*}
		\|f\|_\infty \leq C \left( t + 1 \right)^d \|f\|_\LP.
	\end{equation*}
\end{lemma}

\begin{proof}[Proof of Lemma~\ref{smooth:eigfun2-aux1}]
	We use similar arguments to those in the proof of Corollary~2.5 in \cite{calder2022lipschitz}.
	
	Suppose $|f|$ attains its maximum at $y$, i.e., $\|f\|_\infty = |f(y)|$.
	Let $x \in \M$ with $\dist_\M (x,y) < b$.
	By letting $r_1$ and $C_1$ satisfy Proposition~\ref{Lipschitz:main}, we have for $\e<r_1$,
	\begin{align*}
		\|f\|_\infty - |f(x)| &\le |f(y)-f(x)| \leq C_1 (\|f\|_\infty + \|\bfL_\e f\|_\infty) (\dist_\M(x,y) + \e^{1/3}\log(1/\e))
		\\&\leq C_1 \|f\|_\infty (1+ t) \big(b+\e^{1/3}\log(1/\e)\big).
	\end{align*}
	This implies $\|f\|_\infty \leq 2 |f(x)|$, provided
	\begin{equation}
		C_1 (t + 1)(b+\e^{1/3}\log(1/\e)) \leq \frac{1}{2}, \quad\text{equivalently,}\quad 
		\e^{1/3}\log(1/\e) +b \leq \frac{1}{2C_1(t + 1)}.
		\label{smooth:lipschitz2-eq1}
	\end{equation}
	Under this, we obtain
	\begin{align*}
		b^d \|f\|_\infty \leq 2 \int_{\{x\in\M: \dist_\M(x,y)<b\}} |f(x)| \dP(x) 
		\leq 2 \|f\|_{\LP}.
	\end{align*}
	By adjusting $r$ if necessary, we ensure $\e^{1/3} \log(1/\e) \leq \frac{1}{4C_1(t+1)}$ if $\e < \frac{r}{t+1}$.
	Therefore, by letting $ b = \frac{1}{4 C_1(t +1)}, $ the condition \eqref{smooth:lipschitz2-eq1} is satisfied, and hence,
	\begin{equation*}
		\|f\|_\infty \leq 2 (4C_1)^d (t+1)^d \|f\|_\LP.
	\end{equation*}
	This concludes the proof.
\end{proof}

\begin{lemma} \label{smooth:eigfun2-aux2}
	There exist $(\M,\lambda_j)$-dependent constants $r$ and $C$ such that $\e < r$ implies
	\begin{equation*}
		\max_{k \in \cJ_j} \| \bfL_\e \xi_{k,\e} - \lambda_{k,\e} \xi_{k,\e}\|_\infty  \leq C \e.
	\end{equation*}
\end{lemma}

\begin{proof}[Proof of Lemma~\ref{smooth:eigfun2-aux2}]
	Note that
	$\bfL_\e \phi_{k,\e} = \lambda_{k,\e} \phi_{k,\e}
	= \lambda_{k,\e} \xi_{k,\e} + \lambda_{k,\e} \sum_{\ell \in \cJ_j} v_{k\ell,\e} \phi_\ell $
	and
	\begin{align*}
		\sum_{\ell \in \cJ_j} v_{k\ell,\e} \bfL_\e \phi_\ell
		= \sum_{\ell \in \cJ_j} v_{k\ell,\e} (\lambda_j  \phi_\ell + (\bfL_\e \phi_\ell-\lambda_j\phi_\ell))
		= \lambda_j \sum_{\ell \in \cJ_j} v_{k\ell,\e} \phi_\ell
		+ \sum_{\ell \in \cJ_j} v_{k\ell,\e} (\bfL_\e \phi_\ell-\lambda_j\phi_\ell).
	\end{align*}
	These yield
	\begin{align*}
		\bfL_\e \xi_{k,\e} 
		=\bfL_\e \phi_{k,\e} - \sum_{\ell \in \cJ_j} v_{k\ell,\e} \bfL_\e \phi_\ell
		= \lambda_{k,\e} \xi_{k,\e} + (\lambda_{k,\e} - \lambda_j) \sum_{\ell \in \cJ_j} v_{k\ell,\e} \phi_\ell + \sum_{\ell \in \cJ_j} v_{k\ell,\e} (\bfL_\e \phi_\ell-\lambda_j\phi_\ell),
	\end{align*}
	and consequently,
	$$\| \bfL_\e \xi_{k,\e} - \lambda_{k,\e} \xi_{k,\e}\|_\infty 
	\leq |\lambda_{k,\e}-\lambda_j| \sum_{\ell \in \cJ_j} |v_{k\ell,\e}| \cdot \|\phi_\ell\|_\infty + \sum_{\ell \in \cJ_j} |v_{k\ell,\e}|\cdot \|\bfL_\e \phi_\ell - \lambda_j \phi_\ell\|_\infty. $$
	Therefore, the claim follows by applying Theorems~\ref{smooth:eigval} and \ref{smooth:eigfun1} and Lemma~\ref{smooth:eigval:aux1}.
\end{proof}

Finally, we are ready to prove Theorem~\ref{smooth:eigfun2}.

\begin{proof}[Proof of Theorem~\ref{smooth:eigfun2}]
	We use similar arguments to those in the proof of Theorem~2.6 in \cite{calder2022lipschitz}.
	
	\eqref{smooth:eigfun2-res2} follows directly from \eqref{smooth:eigfun2-res1}, Lemma~\ref{smooth:eigval:aux1} and Theorem~\ref{smooth:eigfun1}.
	Hence, it suffices to prove \eqref{smooth:eigfun2-res1}.
	
	To this end, let $r$ satisfy Theorems~\ref{smooth:eigval}, \ref{smooth:eigfun1}, Proposition~\ref{Lipschitz:main} and Lemmas~\ref{smooth:eigfun2-aux1}, \ref{smooth:eigfun2-aux2}, and let $\e <r$.
	For $k \in \cJ_j$, we consider two possible cases:
	
	\mycase{Case 1. $(\lambda_{k,\e} + 1 ) \|\xi_{k,\e}\|_\infty \leq \| \bfL_\e \xi_{k,\e} \|_\infty$}
	Applying Lemma~\ref{smooth:eigfun2-aux2}, we have $\| \bfL_\e \xi_{k,\e} \|_\infty \leq \lambda_{k,\e} \|\xi_{k,\e}\|_\infty + C_1 \e$ for some $\M$-dependent constant $C_1$.
	This gives
	\begin{equation*}
		(\lambda_{k,\e} + 1 ) \|\xi_{k,\e}\|_\infty \leq \| \bfL_\e \xi_{k,\e} \|_\infty 
		\leq \lambda_{k,\e} \|\xi_{k,\e}\|_\infty  + C_1 \e,
	\end{equation*}
	and it follows 
	\begin{equation}
		\|\xi_{k,\e}\|_\infty \leq C_1 \e, \quad \| \bfL_\e \xi_{k,\e} \|_\infty \leq (C_1+\lambda_{k,\e})\e.
		\label{smooth:eigfun2-eq1}
	\end{equation}
	Inserting \eqref{smooth:eigfun2-eq1} to Proposition~\ref{Lipschitz:main}, we also have
	\begin{equation}
		\sup_{x,y \in \M} \frac{|\xi_{k,\e}(x)-\xi_{k,\e}(y)|}{\dist_\M(x,y) + \e^{1/3} \log(1/\e)} 
		\leq (2C_1 +\lambda_{k,\e})\e.
		\label{smooth:eigfun2-eq2}
	\end{equation}
	Since $\lambda_{k,\e} \leq \lambda_j + |\lambda_{k,\e}-\lambda_j |$, by applying Theorem~\ref{smooth:eigval},
	we see that the claim follows by \eqref{smooth:eigfun2-eq1} and \eqref{smooth:eigfun2-eq2}.
	
	\mycase{Case 2. $(\lambda_{k,\e} + 1 ) \|\xi_{k,\e}\|_\infty > \| \bfL_\e \xi_{k,\e} \|_\infty$}
	By Theorem~\ref{smooth:eigval}, the condition implies 
	\begin{equation*}
		\frac{\| \bfL_\e \xi_{k,\e} \|_\infty}{\|\xi_{k,\e}\|_\infty} \le C_2 (\lambda_j + 1 ),
	\end{equation*}
	for some $(\M,\lambda_j)$-dependent constant $C_2$.
	By Lemma~\ref{smooth:eigfun2-aux1}, we have 
	\begin{equation*}
		\|\xi_{k,\e}\|_\infty \leq C_3 (C_2 ( \lambda_j +2 ))^{d} \| \xi_{k,\e} \|_\LP.
	\end{equation*}
	for some $\M$-dependent constant $C_3$.
	Thus, applying Theorem~\ref{smooth:eigfun1} yields
	\begin{equation}
		\|\xi_{k,\e}\|_\infty \leq C_5 \e, \quad 
		\| \bfL_\e \xi_{k,\e} \|_\infty \leq C_5 \e.
		\label{smooth:eigfun2-eq3}
	\end{equation}
	for some $(\M,\lambda_j)$-dependent constant $C_5$.
	Arguing as in the derivation of \eqref{smooth:eigfun2-eq2}, one obtains
	\begin{equation}
		\sup_{x,y \in \M} \frac{|\xi_{k,\e}(x)-\xi_{k,\e}(y)|}{\dist_\M(x,y) + \e^{1/3} \log(1/\e)} 
		\leq C_6\e.
		\label{smooth:eigfun2-eq4}
	\end{equation}
	for some $(\M,\lambda_j)$-dependent constant $C_6$.
	
	Therefore, the proof is complete by \eqref{smooth:eigfun2-eq3} and \eqref{smooth:eigfun2-eq4}.
\end{proof}

\subsection{Proof of Lemma~\ref{smooth:extra2}} \label{smooth:extra2-proof}

We begin by stating and prove some auxiliary lemmas.

\begin{lemma} \label{smooth:eigfun:aux1}
	For any $k \in \N$, we have
	\begin{equation*}
		\lambda_{k,\e} \sqrt{\rho_\e} \phi_{k,\e} = \bfL_\e^\unnor \left(\frac{\phi_{k,\e}}{\sqrt{\rho_\e}}\right).
	\end{equation*}
\end{lemma}

\begin{proof}[Proof of Lemma \ref{smooth:eigfun:aux1}]
	By the definition of $\phi_{k,\e}$, we have for $x \in \M$,
	\begin{align*}
		\lambda_{k,\e} \phi_{k,\e} (x)
		&= \bfL_\e \phi_{k,\e} (x)
		= \frac{\phi_{k,\e}(x) - \int_\M \frac{\fg_\e (x,y)}{\sqrt{\rho_\e(x)}\sqrt{\rho_\e(y)}} \phi_{k,\e} (y) \dP(y)}{\e}.
	\end{align*}
	The numerator on the right-hand side can be rewritten as
	\begin{align*}
		\phi_{k,\e}(x) - \int_\M \frac{\fg_\e (x,y)}{\sqrt{\rho_\e(x)}\sqrt{\rho_\e(y)}} \phi_{k,\e} (y) \dP(y)
		= \frac{1}{\sqrt{\rho_\e(x)}} \left( \rho_\e (x) \frac{\phi_{k,\e}(x)}{\sqrt{\rho_\e(x)}} - \int_\M \fg_\e (x,y)  \frac{\phi_{k,\e}(y)}{\sqrt{\rho_\e(y)}} \dP(y) \right).
	\end{align*}
	By the definitions of $\rho_\e$ and $\bfL_\e$, it follows that
	\begin{align*}
		\rho_\e (x) \frac{\phi_{k,\e}(x)}{\sqrt{\rho_\e(x)}} - \int_\M \fg_\e (x,y)  \frac{\phi_{k,\e}(y)}{\sqrt{\rho_\e(y)}} \dP(y) 
		&= \int_\M \fg_\e (x,y) \left( \frac{\phi_{k,\e}(x)}{\sqrt{\rho_\e(x)}} - \frac{\phi_{k,\e}(y)}{\sqrt{\rho_\e(y)}} \right) \dP(y)
		\\&= \e \cdot \bfL_\e^\unnor \left( \frac{\phi_{k,\e}}{\sqrt{\rho_\e}} \right) (x).
	\end{align*}
	Therefore, combining the above yields
	\begin{equation*}
		\lambda_{k,\e} \phi_{k,\e} (x) = \frac{\bfL_\e^\unnor \left( \frac{\phi_{k,\e}}{\sqrt{\rho_\e}} \right) (x)}{\sqrt{\rho_\e(x)}}.
	\end{equation*}
	Multiplying both sides by $\sqrt{\rho_\e (x)}$ completes the proof.
\end{proof}

\begin{lemma} \label{smooth:eigfun:aux2}
	For any bounded and measurable functions $f,g : \M \rightarrow \R$, we have
	\begin{equation*}
		\bfL^\unnor_\e (f g) = f \, \bfL_\e^\unnor g + g \, \bfL_\e^\unnor f - \bfK_\e (f,g),
	\end{equation*}
	where $\bfK_\e$ defined as
	\begin{equation*}
		\bfK_\e (f,g) (x) = \int_\M \fg_\e (x,y) \frac{(f(x)-f(y))(g(x)-g(y))}{\e} \dP(y).
	\end{equation*}
\end{lemma}

\begin{proof}[Proof of Lemma \ref{smooth:eigfun:aux2}]
	From the identity
	\begin{equation*}
		f(x)g(x)-f(y)g(y) = f(x) (g(x)-g(y)) + g(x) (f(x)-f(y)) - (f(x)-f(y))(g(x)-g(y)),
	\end{equation*}
	we obtain
	\begin{align*}
		\bfL^\unnor_\e (f g) (x)
		&= \int_\M \fg_\e (x,y) \frac{f(x)g(x)-f(y)g(y)}{\e} \dP(y)
		\\&= f(x) \int_\M \fg_\e (x,y) \frac{g(x)-g(y)}{\e} \dP(y) + g(x) \int_\M \fg_\e (x,y) \frac{f(x)-f(y)}{\e} \dP(y)
		\\&\quad- \int_\M \fg_\e (x,y) \frac{(f(x)-f(y))(g(x)-g(y))}{\e} \dP(y)
		\\&= f (x)\, \bfL_\e^\unnor g (x) + g (x) \, \bfL_\e^\unnor f (x)- \bfK_\e (f,g) (x).
	\end{align*}
	This concludes the proof.
\end{proof}


To simplify the proofs below, we assume $\e < r$ with $r$ satisfying Lemmas~\ref{smooth:KDE} and \ref{gauss:tech:integral}.
Moreover, we let $r$ be sufficiently small that $1/\rho_\e$ is uniformly bounded by 1/2.
We also write $ a = O_\M (b)$ if $|a| \leq C |b|$ for some $\M$-dependent constant $C$.
Similarly, we write $a = O_{\M,\lambda_j} (b)$ when such $C$ is $(\M,\lambda_j)$-dependent.

\begin{lemma} \label{smooth:eigfun:aux3}
	There exist $\M$-dependent constants $r$ and $C$ such that if $\e <r$, we have
	\begin{equation*}
		\| \bfK_\e ( \rho_\e, \rho_\e) \|_\infty
		\leq  C\e^2.
	\end{equation*}
\end{lemma}

\begin{proof}[Proof of Lemma \ref{smooth:eigfun:aux3}]
	Noting $\sup_{x,y \in \M} |\rho(x) - \rho(y)| = O_\M(\e)$,
	plugging $f (\cdot) = (\rho(x) - \rho(\cdot))^2/\e$ into Lemma~\ref{gauss:tech:integral} yields
	\begin{align*}
		\bfK_\e ( \rho_\e, \rho_\e) (x) = \int_{B_\e} \tfg (z) \frac{\big(\rho_\e (x) -\rho_\e (\Exp_x(\sqrt{\e}z)\big)^2}{\e} dz + O_\M (\e^2).
	\end{align*}
	Similar to the proof of Lemma~\ref{smooth:KDE}, we obtain the Taylor expansion of $(\rho_\e (x) -\rho_\e (\Exp_x(\sqrt{\e}z))^2/\e$ around $z=0$ as 
	\begin{align*}
		\frac{(\rho_\e (x) -\rho_\e (\Exp_x(\sqrt{\e}z))^2}{\e}
		= \frac{( \e^{3/2} O_\M(\|z\|_{\R^d}) + \e^2 O_\M(\|z\|_{\R^d}^2))^2}{\e}
		= \e^2 O_\M(\|z\|^4_{\R^d}).
	\end{align*}
	Since $\tfg$ is a Gaussian kernel, the claim follows.
\end{proof}

\begin{lemma} \label{smooth:eigfun:aux4}
	There exist $\M$-dependent constants $r$ and $C$ such that if $\e <r$, then
	for any bounded measurable function $f :\M \to \R$, we have
	\begin{equation*}
		\left\| \bfK_\e \left( \frac{f}{\sqrt{\rho_\e}}, \frac{f}{\sqrt{\rho_\e}} \right) - \bfK_\e (f,f) \right\|_\infty
		\leq C \Big( \e \bfK_\e (f,f) + \e \|f\|_\infty \sqrt{\bfK_\e (f,f)} + \e^2 \|f\|_\infty^2\Big).
	\end{equation*}
\end{lemma}

\begin{proof}[Proof of Lemma \ref{smooth:eigfun:aux4}]
	From the identity
	\begin{align*}
		\frac{f(x)}{\sqrt{\rho_\e(x)}}- \frac{f(y)}{\sqrt{\rho_\e(y)}}
		&= f(x) - f(y) + \frac{1-\rho_\e(x)}{\sqrt{\rho_\e(x)}(\sqrt{\rho_\e(x)}+1)} (f(x) - f(y))
		\\&\quad+ \frac{f(y)}{\sqrt{\rho_\e(x)} \sqrt{\rho_\e(y)}(\sqrt{\rho_\e(x)}+\sqrt{\rho_\e(y)})} ( \rho_\e(y) - \rho_\e (x)),
	\end{align*}
	by squaring both sides and applying Lemma~\ref{smooth:KDE}, we obtain
	\begin{align*}
		&\quad \bigg(\frac{f(x)}{\sqrt{\rho_\e(x)}}- \frac{f(y)}{\sqrt{\rho_\e(y)}}\bigg)^2 - (f(x)-f(y))^2
		\\&= O_\M \Big(\e (f(x)-f(y))^2 + \|f\|_\infty |f(x)-f(y)||\rho_\e(x)-\rho_\e(y)|
		+ \|f\|_\infty (\rho_\e (x) - \rho_\e (y))^2 \Big).
	\end{align*}
	Noting 
	\begin{align*}
		\int_\M \fg_\e (x,y) \frac{|f(x)-f(y)||\rho_\e(x)-\rho_\e(y)|}{\e} d \P(y)
		\leq \sqrt{\bfK_\e (f,f) (x)} \sqrt{\bfK_\e (\rho_\e,\rho_\e) (x)},
	\end{align*}
	the definition of $\bfK_\e$ yields
	\begin{multline*}
		\bfK_\e \left( \frac{f}{\sqrt{\rho_\e}}, \frac{f}{\sqrt{\rho_\e}} \right) (x) -
		\bfK_\e \left( f, f \right) (x)
		\\= O_\M  \left( \e \bfK_\e (f,f) (x) + \|f\|_\infty \sqrt{\bfK_\e (f,f) (x)} \sqrt{\bfK_\e (\rho_\e,\rho_\e) (x)}+ \|f\|_\infty \bfK_\e(\rho_\e, \rho_\e) (x) \right).
	\end{multline*}
	Therefore, by applying Lemma~\ref{smooth:eigfun:aux3}, we complete the proof.
\end{proof}

Below, we additionally assume $\e<r$ with $r$ satisfying Theorems~\ref{smooth:eigval}, \ref{smooth:eigfun1} and \ref{smooth:eigfun2}, as well as Lemma~\ref{smooth:eigfun:aux4}.

\begin{lemma} \label{smooth:eigfun:part1}
	There exist $(\M,\lambda_j)$-dependent constants $r$ and $C$ such that if $\e <r$,
	\begin{equation*}
		\max_{k,\ell \in \cJ_j} \left\| \bfK_\e \left( \frac{\phi_{k}}{\sqrt{\rho_\e}}, \frac{\phi_{\ell}}{\sqrt{\rho_\e}} \right)
		- 2 \lambda_j \phi_{k} \phi_{\ell} + \Delta_\M (\phi_{k} \phi_{\ell})\right\|_\infty
		\leq C\e.
	\end{equation*}
\end{lemma}

\begin{proof}[Proof of Lemma~\ref{smooth:eigfun:part1}]
	For $k,\ell \in \cJ_j$,
	applying $f (\cdot) = (\phi_k(x) - \phi_k(\cdot))(\phi_\ell(x) - \phi_\ell(\cdot))/\e$ to Lemma~\ref{gauss:tech:integral} gives
	\begin{align*}
		\bfK_\e \left( \phi_k, \phi_\ell \right) 
		= \int_{B_\e} \tfg (z) \tf_{x} (z) \big(1+ O_\M(\e\|z\|^2_{\R^d})\big) dz + O_{\M,\lambda_j} (\e),
	\end{align*}
	where 
	\begin{equation*}
		\tf_{x} (z) = \frac{\big(\phi_k (x) - \tphi_{k,x}(z)\big)\big(\phi_\ell (x) -\tphi_{\ell,x} (z)\big)}{\e},
		\quad \tphi_{k,x} (z) = \phi_k (\Exp_x(\sqrt{\e}z).
	\end{equation*}	
	Using the Taylor expansion of $\tphi_{k,x}$ around $z=0$, we expand $\tf_x(z)$ as
	\begin{multline*}
		\tf_x(z) = z^\top \big[\nabla_{\R^d} \tphi_{k,x} (0)\big] \big[ \nabla_{\R^d} \tphi_{k,x} (0) \big]^\top z
		+ \frac{\sqrt{\e}}{2} z^\top \big[\nabla_{\R^d} \tphi_{k,x} (0)\big] z \big[\nabla_{\R^d} \tphi_{\ell,x} (0)\big]^\top z
		\\+ \frac{\sqrt{\e}}{2} z^\top \big[\nabla_{\R^d} \tphi_{\ell,x} (0)\big] z \big[\nabla_{\R^d} \tphi_{k,x} (0)\big]^\top z
		+ O_{\M,\lambda_j} (\e \|z\|_{\R^d}^6),
	\end{multline*}
	and consequently,
	\begin{align*}
		\bfK_\e \left( \phi_k, \phi_\ell \right) 
		&= \int_{B_\e} \tfg(z) z^\top \big[\nabla_{\R^d} \tphi_{k,x} (0)\big] \big[ \nabla_{\R^d} \tphi_{k,x} (0) \big]^\top z
		+ O_{\M,\lambda_j}(\e)
		\\&= \int_{\R^d} \tfg(z) z^\top \big[\nabla_{\R^d} \tphi_{k,x} (0)\big] \big[ \nabla_{\R^d} \tphi_{k,x} (0) \big]^\top z
		+ O_{\M,\lambda_j}(\e). 
	\end{align*}
	where in the last line we use the same integral approximation as in the proof of Lemma~\ref{Lipschitz:normalapprox}.
	Since $\tfg$ is the PDF of the $d$-dimensional mean-zero normal distribution with covariance $2 I_d$, we evaluate the integral as
	\begin{align*}
		\int_{\R^d} \tfg(z) \, z^\top [\nabla_{\R^d} \tphi_{k,x} (0)] [\nabla_{\R^d} \tphi_{\ell,x} (0)]^\top z dz 
		&= 2 \langle \nabla_{\R^d} \tphi_{k,x} (0), \nabla_{\R^d} \tphi_{\ell,x} (0) \rangle_{\R^d}
		\\&= 2 \langle \nabla_{\M} \phi_k (x), \nabla_{\M} \phi_\ell (x) \rangle_{T_x \M}.
	\end{align*}
	Since
	\begin{align*}
		\Delta_\M (\phi_k  \phi_\ell) 
		&= \phi_\ell \Delta_\M \phi_k + \phi_k \Delta_\M \phi_\ell -2 \langle \nabla_{\M} \phi_k, \nabla_{\M} \phi_\ell \rangle_{T\M}
		\\&= 2 \lambda_j \phi_k \phi_\ell -2 \langle \nabla_{\M} \phi_k, \nabla_{\M} \phi_\ell \rangle_{T\M},
	\end{align*}
	we obtain
	\begin{equation*}
		\bfK_\e \left( \phi_k, \phi_\ell \right) = 2\lambda_j \phi_k \phi_\ell - \Delta_\M (\phi_k\phi_\ell) + O_{\M,\lambda_j}(\e).
	\end{equation*}	
	Therefore, applying Lemma~\ref{smooth:eigfun:aux4} completes the proof.
\end{proof}

\begin{lemma} \label{smooth:eigfun:part2}
	There exist $(\M,\lambda_j)$-dependent constants $r$ and $C$ such that $\e < r$ implies
	\begin{equation*}
		\max_{k \in \cJ_j} \left\| \bfK_\e \left( \frac{\xi_{k,\e}}{\sqrt{\rho_\e}}, \frac{\xi_{k,\e}}{\sqrt{\rho_\e}} \right) \right\|_\infty \leq C \e^{\frac{5}{3}} \log^2(1/\e).
	\end{equation*}
\end{lemma}

\begin{proof}[Proof of Lemma \ref{smooth:eigfun:part2}]
	For $k \in \cJ_j$, plugging $f (\cdot) = (\xi_{k,\e}(x) - \xi_{k,\e}(\cdot))^2/\e$ into Lemma~\ref{gauss:tech:integral} gives
	\begin{align*}
		\bfK_\e \left( \xi_{k,\e}, \xi_{k,\e} \right)
		= \int_{B_\e} \tfg (z) \frac{\big(\xi_{k,\e}(x)-\xi_{k,\e}(\Exp_x (\sqrt{\e}z))\big)^2}{\e} dz + O_{\M,\lambda_j} (\e^2),
	\end{align*}
	Applying Theorem \ref{smooth:eigfun2}, we bound
	\begin{align*}
		\frac{ |\xi_{k,\e}(x)-\xi_{k,\e}(\Exp_x (\sqrt{\e}z))|^2}{\e}
		&= O_{\M,\lambda_j} \bigg(  \e \Big( \dist_\M(x,\Exp_x (\sqrt{\e}z)) + \e^{1/3} \log(1/\e) \Big)^2 \bigg)
		\\&= O_{\M,\lambda_j} \bigg(\e^{\frac{5}{3}} \log^2(1/\e) \Big( \e^{\frac{1}{3}} \log^{-2} \log(1/\e) \|z|\|_{\R^d}^2+1\Big) \bigg).
	\end{align*}
	Plugging this to the first equation yields
	\begin{equation*}
		\bfK_\e \left( \xi_{k,\e}, \xi_{k,\e} \right) = O_{\M,\lambda_j} (\e^{\frac{5}{3}} \log^2(1/\e)).
	\end{equation*}
	Therefore, applying Lemma~\ref{smooth:eigfun:aux4} completes the proof.
\end{proof}

\begin{lemma} \label{smooth:eigfun:part3}
	There exist $(\M,\lambda_j)$-dependent constants $r$ and $C$ such that $\e < r$ implies
	\begin{equation*}
		\max_{k,\ell \in \cJ_j} \left\| \bfK_\e \left( \frac{\phi_{k}}{\sqrt{\rho_\e}}, \frac{\xi_{\ell,\e}}{\sqrt{\rho_\e}} \right) \right\|_\infty \leq C \e^{5/6} \log(1/\e).
	\end{equation*}
\end{lemma}

\begin{proof}[Proof of Lemma~\ref{smooth:eigfun:part3}]
	Applying the Cauchy–Schwarz inequality gives
	\begin{equation*}
		\max_{k, \ell \in \cJ_j} \left\| \bfK_\e \left( \frac{\phi_{k}}{\sqrt{\rho_\e}}, \frac{\xi_{\ell,\e}}{\sqrt{\rho_\e}} \right) \right\|_\infty
		\leq \sqrt{ \max_{k \in \cJ_j}\left\| \bfK_\e \left( \frac{\phi_{k}}{\sqrt{\rho_\e}}, \frac{\phi_{k}}{\sqrt{\rho_\e}} \right) \right\|_\infty }
		\sqrt{ \max_{\ell \in \cJ_j} \left\| \bfK_\e \left( \frac{\xi_{\ell,\e}}{\sqrt{\rho_\e}}, \frac{\xi_{\ell,\e}}{\sqrt{\rho_\e}} \right) \right\|_\infty}.
	\end{equation*}
	Hence, the claim follows by Lemmas~\ref{smooth:eigfun:part1} and \ref{smooth:eigfun:part2}.
\end{proof}

Now, we are ready to prove Lemma~\ref{smooth:extra2}.

\begin{proof}[Proof of Lemma~\ref{smooth:extra2}]
	Once we show 
	\begin{equation}
		\Bigg\| \sum_{k \in \cJ_j} \mu_{k,\e}\bfL_\e^\unnor \bigg( \frac{\phi^2_{k,\e} }{\rho_\e} \bigg)
		- \sum_{k \in \cJ_j} \Delta_\M (\phi^2_k) \Bigg\|_\infty
		= O_{\M,\lambda_j} \big(\e^{5/6} \log(1/\e)\big)
		\label{smooth:eigfun:sum2-eq1}
	\end{equation}
	the claim follows from Hoeffding's inequality.
	
	Let $k \in \cJ_j$. By Lemmas \ref{smooth:eigfun:aux1} and \ref{smooth:eigfun:aux2}, we have
	\begin{equation}
		\bfL_\e^\unnor \left( \frac{\phi^2_{k,\e} }{\rho_\e} \right)
		= \frac{2\phi_{k,\e}}{\sqrt{\rho_\e}} \bfL_\e^\unnor \left( \frac{\phi_{k,\e} }{\sqrt{\rho_\e}} \right)
		- \bfK_\e \left( \frac{\phi_{k,\e}}{\sqrt{\rho_\e}},\frac{\phi_{k,\e}}{\sqrt{\rho_\e}} \right)
		= 2 \lambda_{k,\e} \phi_{k,\e}^2 - \bfK_\e \left( \frac{\phi_{k,\e}}{\sqrt{\rho_\e}},\frac{\phi_{k,\e}}{\sqrt{\rho_\e}} \right).
		\label{smooth:eigfun:sum2-eq2}
	\end{equation}
	We analyze the two terms on the right-hand side separately.
	
	First, from	
	\begin{align*}
		2 \lambda_{k,\e} \phi_{k,\e}^2 
		&= 2 \lambda_j \phi_{k,\e}^2 + 2 (\lambda_{k,\e}-\lambda_j) \phi_{k,\e}^2 
		= 2 \lambda_j \phi_{k,\e}^2 + O_{\M,\lambda_j}(\e ), \\
		\phi_{k,\e} &= \sum_{m \in \cJ_j} v_{km,\e} \phi_m + \xi_{k,\e}
		= \sum_{m \in \cJ_j} v_{km,\e} \phi_m + O_{\M,\lambda_j}(\e ),
	\end{align*}
	we obtain
	\begin{equation}
		2 \lambda_{k,\e} \phi_{k,\e}^2 
		= 2 \lambda_j \sum_{m_1 \in \cJ_j} \sum_{m_2 \in \cJ_j} v_{k m_1,\e} v_{k m_2,\e} \phi_{m_1} \phi_{m_2} + O_{\M,\lambda_j} (\e ).
		\label{smooth:eigfun:sum2-eq3}
	\end{equation}
	Next, again using $ \phi_{k,\e} = \sum_{m \in \cJ_j} v_{km,\e} \phi_m + \xi_{k,\e},$ we have for $x,y \in \M$,
	\begin{align*}
		\left( \frac{\phi_{k,\e}(x)}{\sqrt{\rho_\e(x)}}-\frac{\phi_{k,\e}(y)}{\sqrt{\rho_\e(y)}} \right)^2
		&= \sum_{m_1 \in \cJ_j} \sum_{m_2 \in \cJ_j} v_{km_1,\e} v_{km_2,\e} 
		\left( \frac{\phi_{m_1}(x)}{\sqrt{\rho_\e(x)}} - \frac{\phi_{m_1}(y)}{\sqrt{\rho_\e(y)}}\right) \left( \frac{\phi_{m_2}(x)}{\sqrt{\rho_\e(x)}} - \frac{\phi_{m_2}(y)}{\sqrt{\rho_\e(y)}}\right) 
		\\&\quad + 2 \sum_{m \in \cJ_j} v_{km,\e} \left( \frac{\phi_{m}(x)}{\sqrt{\rho_\e(x)}} - \frac{\phi_{m}(y)}{\sqrt{\rho_\e(y)}}\right) \left( \frac{\xi_{k,\e}(x)}{\sqrt{\rho_\e(x)}} - \frac{\xi_{k,\e}(y)}{\sqrt{\rho_\e(y)}}\right)
		\\&\quad + \left( \frac{\xi_{k,\e}(x)}{\sqrt{\rho_\e(x)}} - \frac{\xi_{k,\e}(y)}{\sqrt{\rho_\e(y)}}\right)^2,
	\end{align*}
	and consequently,
	\begin{align*}
		\bfK_\e \left( \frac{\phi_{k,\e}}{\sqrt{\rho_\e}},\frac{\phi_{k,\e}}{\sqrt{\rho_\e}} \right)
		&= \sum_{m_1 \in \cJ_j} \sum_{m_2 \in \cJ_j} v_{km_1,\e} v_{km_2,\e} 
		\bfK_\e \left( \frac{\phi_{m_1}}{\sqrt{\rho_\e}}, \frac{\phi_{m_2}}{\sqrt{\rho_\e}} \right) 
		\\&\quad + 2 \sum_{m \in \cJ_j} v_{km,\e} 
		\bfK_\e \left( \frac{\phi_{m}}{\sqrt{\rho_\e}}, \frac{\xi_{k,\e}}{\sqrt{\rho_\e}} \right)
		+ \bfK_\e \left( \frac{\xi_{k,\e}}{\sqrt{\rho_\e}}, \frac{\xi_{k,\e}}{\sqrt{\rho_\e}} \right).
	\end{align*}
	Hence, by applying Lemmas \ref{smooth:eigfun:part1}-\ref{smooth:eigfun:part3}, we obtain
	\begin{equation}
		\bfK_\e \left( \frac{\phi_{k,\e}}{\sqrt{\rho_\e}},\frac{\phi_{k,\e}}{\sqrt{\rho_\e}} \right)
		= \sum_{m_1 \in \cJ_j} \sum_{m_2 \in \cJ_j} v_{km_1,\e} v_{km_2,\e} 
		(2 \lambda_j \phi_{m_1} \phi_{m_2} - \Delta_\M (\phi_{m_1} \phi_{m_2}) ) + O_{\M,\lambda_j} \big(\e^{5/6} \log(1/\e)\big).
		\label{smooth:eigfun:sum2-eq4}
	\end{equation}
	
	Now, combining \eqref{smooth:eigfun:sum2-eq2}-\eqref{smooth:eigfun:sum2-eq4} yields
	\begin{align*}
		\bfL_\e^\unnor \bigg( \frac{\phi^2_{k,\e} }{\rho_\e} \bigg)
		&= \sum_{m_1 \in \cJ_j} \sum_{m_2 \in \cJ_j} v_{k m_1,\e} v_{k m_2,\e} \Delta_\M (\phi_{m_1} \phi_{m_2})
		+ O_{\M,\lambda_j} \big(\e^{5/6} \log(1/\e)\big).
	\end{align*}
	Moreover, from $\mu_{k,\e} = 1- \e\lambda_{k,\e} = 1 + O_{\M,\lambda_j} (\e)$, we have
	\begin{align*}
		\mu_{k,\e} \bfL_\e^\unnor \bigg( \frac{\phi^2_{k,\e} }{\rho_\e} \bigg)
		= \sum_{m_1 \in \cJ_j} \sum_{m_2 \in \cJ_j} v_{k m_1,\e} v_{k m_2,\e} \Delta_\M (\phi_{m_1} \phi_{m_2})
		+ O_{\M,\lambda_j} \big(\e^{5/6} \log(1/\e)\big).
	\end{align*}
	Therefore, by letting $A=(\Delta_\M(\phi_{m_1,m_2}))_{m_1,m_2 \in \cJ_j}$, we obtain 
	\begin{equation*}
		\sum_{k \in \cJ_j} \mu_{k,\e} \bfL_\e^\unnor \bigg( \frac{\phi^2_{k,\e} }{\rho_\e} \bigg) - 
		\sum_{k \in \cJ_j} \Delta_\M (\phi^2_k)
		= \tr(V_{j,\e} A V_{j,\e}^\top ) - \tr(A)
		+ O_{\M,\lambda_j} \big(\e^{5/6} \log(1/\e)\big).
	\end{equation*}
	Finally, by applying Theorem~\ref{smooth:eigfun1} and using $\|A\|_\F = O_{\M,\lambda_j}(1)$, we obtain \eqref{smooth:eigfun:sum2-eq1}.
	This completes the proof.
\end{proof}

\section{Spectral convergence of compact self-adjoint operators} \label{gauss:spec}

In this section, we briefly review the spectral convergence results for compact self-adjoint operators from Section~3 of \cite{bae2026asymptotic} and present them in a form suitable for our analysis.
Unless otherwise stated, we follow the notation introduced in the main text. 
To simplify the exposition, some symbols are reused with meanings specific to this appendix.
For instance, $\lambda$ denotes an eigenvalue of the operator considered here, rather than one in the main text.
The same convention applies to the subsequent appendices.

Let $\mbH$ be a real separable Hilbert space, and let $\hcH$ and $\cH$ be compact self-adjoint operators on $\mbH$.
This setting includes symmetric matrices as a special case, as well as symmetrically normalized smoothed graph Laplacians.
While the normalized smooth graph Laplacian itself is not compact, it can be written as the difference between the identity operator and a compact self-adjoint operator, and hence, it is covered by this setting.

We further let $\{\hlambda_k\}_{k \in \N}$ and $\{\lambda_k\}_{k \in \N}$ denote the eigenvalues of $\hcH$ and $\cH$,
respectively, counted with their multiplicities and arranged in non-increasing order.
We denote the corresponding $\mbH$-orthonormal eigenvectors by $\{\hpsi_k\}_{k\in\N}$ and $\{\psi_k\}_{k\in\N}$, respectively.

Let $\cJ \subset \N$ be a finite index set satisfying $\gammaJ = \min_{k \in \cJ, \ell \notin \cJ} |\lambda_k - \lambda_\ell|>0.$
We let $\theta_\tmax = \theta_1 > \dots > \theta_D = \theta_\tmin $ denote the distinct eigenvalues of $\cH$ in $\{\lambda_k\}_{k \in \cJ}$. When $D=1$, we simply write $\theta_\tmax = \theta_\tmin = \theta$. Note that $\{\theta_d\}_{d \in [D]}$ induce a partition $\{\cJ_d\}_{d \in [D]}$ of $\cJ$, where $\cJ_d = \{k \in \cJ: \lambda_k = \theta_d \}$.

The goal of this section is to quantify the convergence of the eigenvalues and eigenvectors of $\hcH$ indexed by $\cJ$ to the corresponding eigenvalues and eigenvectors of $\cH$.

\subsection{Weighted perturbation of compact self-adjoint operators}

The spectral convergence of $\hcH$ to $\cH$ can be characterized by two weighted perturbation quantities $\halpha_{\!\cJ}$ and $\hbeta_{\!\cJ}$, defined by
\begin{equation*}
	\halpha_{\!\cJ} = \big\|\cT_{\!\cJ}^{1/2} (\hcH - \cH) \cT_{\!\cJ}^{1/2} \big\|_{\op,\mbH}, \quad 
	\hbeta_{\!\cJ} = \gammaJ^{-1/2} \big\||\cR|_{\cJ}^{1/2} (\hcH-\cH) \cP_{\!\cJ} \big\|_{\HS,\mbH},
\end{equation*}
where
\begin{equation*}
	\cT_{\!\cJ}^{1/2} = \frac{1}{\sqrt{\gammaJ}} \cP_{\!\cJ} + |\cR|_{\!\cJ}^{1/2}, \quad
	|\cR|_{\cJ}^{1/2} = \sum_{k \notin \cJ} \frac{\cQ_k}{\sqrt{\min_{m \in \cJ} |\lambda_k - \lambda_m|}}, \quad
	\cQ_k = \psi_k \otimes_\mbH \psi_k, \quad
	\cP_{\!\cJ} = \sum_{k \in \cJ} \cQ_k.
\end{equation*}

We present how $\halpha_{\!\cJ}$ and $\hbeta_{\!\cJ}$ can be controlled under two situations that we use in the main text.

\begin{lemma}[Relative form boundedness] \label{gauss:spec:perturb1}
	Suppose there exist $\epsilon_1, \epsilon_2 >0$ such that
	\begin{equation*}
		|\langle (\hcH-\cH) v , v \rangle_\mbH| \leq \epsilon_1 \langle \cH v, v \rangle_\mbH + \epsilon_2 \|v\|_\mbH \quad \text{for all $v \in \mbH$},
	\end{equation*}
	and $\gammaJ \leq \theta_\tmax$.
	Then we have 
	\begin{equation*}
		\halpha_{\!\cJ} \leq \gammaJ^{-1} (2 \epsilon_1 \theta_\tmax + \epsilon_2), \quad
		\hbeta_{\!\cJ} \leq \gammaJ^{-1} \sqrt{\sum_{k \in \cJ} \|(\hcH-\cH)\psi_k\|_\mbH^2 }.
	\end{equation*}
\end{lemma}

\begin{proof}[Proof of Lemma~\ref{gauss:spec:perturb1}]
	For $\halpha_{\!\cJ}$, see Example~3.3 in \cite{bae2026asymptotic}.
	
	For $\hbeta_{\!\cJ}$, using $\|\cP_{\!\cJ}\|_{\op,\mbH}=1$ and $\||\cR|_{\!\cJ}^{1/2}\|_{\op,\mbH}\leq \gammaJ^{-1/2}$, it follows that
	\begin{equation*}
		\hbeta_{\!\cJ} 
		= \gammaJ^{-1/2} \big\| |\cR|_{\!\cJ}^{1/2} (\hcH-\cH) \cP_{\!\cJ} \big\|_{\HS,\mbH}
		\leq \gammaJ^{-1} \big\| (\hcH-\cH) \cP_{\!\cJ} \big\|_{\HS,\mbH}
		= \gammaJ^{-1} \sqrt{\sum_{k \in \cJ} \|(\hcH-\cH)\psi_k\|_\mbH^2 }.
	\end{equation*}
	This concludes the proof.
\end{proof}

\begin{remark}
	The condition $\gammaJ \leq \theta_\tmax$ is imposed just to simplify expressions.
	In fact, it always holds whenever there is an eigenvalue of $\cH$ smaller than $\theta_\tmin$.
\end{remark}

\begin{lemma}[Classical perturbation] \label{gauss:spec:perturb2}
	It holds that
	\begin{equation*}
		\halpha_{\!\cJ} \leq \gammaJ^{-1} \|\hcH-\cH\|_{\op,\mbH}, \quad
		\hbeta_{\!\cJ} \leq |\cJ| \gammaJ^{-1} \|\hcH-\cH\|_{\op,\mbH}.
	\end{equation*}
\end{lemma}

\begin{proof}[Proof of Lemma~\ref{gauss:spec:perturb2}]
	The claim follows immediately by using $\|\cP_{\!\cJ}\|_{\op,\mbH}=1$, $\|\cP_{\!\cJ}\|_{\HS,\mbH}=|\cJ|$, and $\||\cR|_{\!\cJ}^{1/2}\|_{\op,\mbH}\leq \gammaJ^{-1/2}$.
\end{proof}

\subsection{Eigenvalues of compact self-adjoint operators}

In this part and the following one, we present spectral convergence results for eigenvalues and eigenvectors, respectively.

\begin{proposition} \label{gauss:spec:eigval1}
	Under the same setting as in Lemma~\ref{gauss:spec:perturb1}, if $D=1$ and $\gammaJ^{-1} (2 \epsilon_1 \theta + \epsilon_2) \leq 1/4$, we have
	\begin{equation*}
		\big\| (\hlambda_k - \theta)_{k \in \cJ} \big\|_2
		\leq \frac{36+7\sqrt{2}}{18} \sqrt{\sum_{k \in \cJ} \| (\hcH - \cH) \psi_k \|_\mbH^2 }.
	\end{equation*}
\end{proposition}

\begin{proof}[Proof of Proposition~\ref{gauss:spec:eigval1}]
	According to Theorem~3.9 in \cite{bae2026asymptotic}, under $\halpha_{\!\cJ} \le 1/4$ and $D=1$, it holds that
	\begin{equation*}
		\Bigg\| (\hlambda_k - \theta)_{k \in \cJ} - \specd \Big( \langle (\hcH-\cH) \psi_k, \psi_\ell \rangle_\mbH \Big)_{k,\ell \in \cJ} \Bigg\|_2
		\leq \frac{36+14\sqrt{2}}{9}\gammaJ\halpha_{\!\cJ}\hbeta_{\!\cJ} 
		\leq \frac{18+7\sqrt{2}}{18}\gammaJ\hbeta_{\!\cJ},
	\end{equation*}
	where $\specd(\cdot)$ denotes the vector of eigenvalues arranged in non-increasing order.
	
	By the properties of the Frobenius norm and the Hilbert-Schmidt norm, it follows that 
	\begin{align*}
		\bigg\| \specd \Big( \langle (\hcH-\cH) \psi_k, \psi_\ell \rangle_\mbH \Big)_{k,\ell \in \cJ_d} \bigg\|_2
		= \bigg\| \Big( \langle (\hcH-\cH) \psi_k, \psi_\ell \rangle_\mbH \Big)_{k,\ell \in \cJ_d} \bigg\|_\F
		= \|\cP_\cJ (\hcH-\cH) \cP_\cJ\|_{\HS,\mbH}.
	\end{align*}
	Using $\|\cP_{\!\cJ}\|_{\op,\mbH}=1$, we bound 
	\begin{equation*}
		\|\cP_\cJ (\hcH-\cH) \cP_\cJ\|_{\HS,\mbH} 
		\leq \|(\hcH-\cH) \cP_\cJ\|_{\HS,\mbH} 
		= \sqrt{\sum_{k \in \cJ} \|(\hcH-\cH)\psi_k\|_\mbH^2 },
	\end{equation*}
	Therefore, the claim follows by combining the above and applying Lemma~\ref{gauss:spec:perturb1}.
\end{proof}

\begin{proposition} \label{gauss:spec:eigval2}
	Under the same setting as in Lemma~\ref{gauss:spec:perturb1}, if $D=1$ and $\gammaJ^{-1} (2 \epsilon_1 \theta + \epsilon_2) \leq 1/4$, we have
	\begin{equation*}
		\min_{k \in \cJ, \ell \notin \cJ} |\hlambda_k - \hlambda_\ell| \ge \frac{\gammaJ}{2}.
	\end{equation*}
\end{proposition}

\begin{proof}[Proof of Lemma~\ref{gauss:spec:eigval2}]
	According to Lemma~B.2 of \cite{bae2026asymptotic}, under $\halpha_\cJ \le 1/4$, it holds that
	\begin{equation}
		\max_{k \in \cJ} |\hlambda_k - \lambda_k| \leq \gammaJ \halpha_{\!\cJ}, \quad
		\min_{k \in \cJ, \ell \notin \cJ} |\lambda_k - \hlambda_\ell| \ge \gammaJ(1-\halpha_{\!\cJ}).
		\label{gauss:spec:eigval3-eq1}
	\end{equation}
	These yield
	\begin{align}
		\min_{k \in \cJ, \ell \notin \cJ} |\hlambda_k - \hlambda_\ell|
		\ge \min_{k \in \cJ, \ell \notin \cJ} |\lambda_k - \hlambda_\ell| - \max_{k \in \cJ} |\hlambda_k - \lambda_k|
		\ge \gammaJ(1-2\halpha_\cJ) \ge \frac{\gammaJ}{2}.
		\label{gauss:spec:eigval3-eq2}
	\end{align}
	Hence, the claim follows by Lemma~\ref{gauss:spec:perturb1}.
\end{proof}

\begin{proposition} \label{gauss:spec:eigval3}
	If $\gammaJ^{-1} \|\hcH-\cH\|_{\op,\mbH} \leq 1/4$, we have
	\begin{equation*}
		\bigg| \sum_{k\in\cJ} \hlambda_k - \sum_{k\in\cJ} \lambda_k
		- \sum_{k \in \cJ} \langle (\hcH-\cH) \psi_k, \psi_k \rangle_\mbH \bigg|
		\leq \bigg(\frac{5}{2\gammaJ} + \frac{4\theta_\tmax-4\theta_\tmin}{\gammaJ^2}\bigg) \|\hcH-\cH\|_{\op,\mbH}^2.
	\end{equation*}
\end{proposition}

\begin{proof}[Proof of Proposition~\ref{gauss:spec:eigval3}]
	According to Theorem~3.11 in \cite{bae2026asymptotic},
	under $\halpha_{\cJ} \le 1/4$, it holds that
	\begin{equation*}
		\bigg| \sum_{k\in\cJ} \hlambda_k - \sum_{k\in\cJ} \lambda_k
		- \sum_{k \in \cJ} \langle (\hcH-\cH) \psi_k, \psi_k \rangle_\mbH \bigg|
		\\ \leq \bigg(\frac{5}{2}\gammaJ + 4\theta_\tmax-4\theta_\tmin\bigg) \hbeta_{\!\cJ}^2.
	\end{equation*}
	Thus, the claim follows by Lemma~\ref{gauss:spec:perturb2}.
\end{proof}

\subsection{Eigenvectors of compact self-adjoint operators}


\begin{proposition} \label{gauss:spec:eigvec1}
	Under the same setting as in Lemma~\ref{gauss:spec:perturb1}, if $D=1$ and $\gammaJ^{-1} (2 \epsilon_1 \theta + \epsilon_2) \leq 1/4$, we have
	\begin{equation*}
		\sqrt{ \sum_{k \in \cJ} \Big\| \hpsi_k - \sum_{\ell \in \cJ} \langle \hpsi_k, \psi_\ell \rangle_\mbH \psi_\ell
			\Big\|_\mbH^2 }
		\leq \frac{91}{18 \gammaJ} \sqrt{\sum_{k \in \cJ} \| (\hcH - \cH) \psi_k \|_\mbH^2 }.
	\end{equation*}
\end{proposition}

\begin{proof}[Proof of Proposition~\ref{gauss:spec:eigvec1}]
	According to Theorem~3.7 of \cite{bae2026asymptotic}, under $\halpha_{\!\cJ} \leq 1/4$, it holds that
	\begin{multline*}
		\sqrt{ \sum_{k \in \cJ} \bigg\| \hpsi_k - \sum_{\ell \in \cJ} \langle \hpsi_k, \psi_\ell \rangle_\mbH 
			\bigg( \psi_\ell + \sum_{m \notin \cJ} \frac{ \langle (\hcH-\cH) \psi_\ell, \psi_m \rangle_\mbH}{\lambda_\ell - \lambda_m}  \psi_m \bigg) \bigg\|_\mbH^2 }
		\\ \leq \bigg( \frac{74\sqrt{D}}{9} + 8 \bigg) \halpha_{\!\cJ} \hbeta_{\!\cJ}
		\leq \bigg( \frac{37\sqrt{D}}{18} + 2 \bigg) \hbeta_{\!\cJ}.
	\end{multline*}
	A simple calculation gives
	\begin{equation*}
		\sqrt{ \sum_{k \in \cJ} \bigg\| \sum_{\ell \in \cJ} \langle \hpsi_k, \psi_\ell \rangle_\mbH 
			\sum_{m \notin \cJ} \frac{ \langle (\hcH-\cH) \psi_\ell, \psi_m \rangle_\mbH}{\lambda_\ell - \lambda_m}  \psi_m \bigg\|^2_\mbH }
		= \bigg\| \sum_{d=1}^D \cR_{d,\cJ} (\hcH-\cH) \cP_{\!\cJ_d} \hcP_{\!\cJ} \bigg\|_{\HS,\mbH},
	\end{equation*}
	where $\hcP_{\!\cJ} = \sum_{k\in\cJ} \hpsi_k \otimes_\mbH \hpsi_k$ and $\cR_{d,\cJ} = \sum_{m \notin \cJ} (\theta_d-\lambda_m)^{-1} \cQ_m$. 
	Using $\|\hcP_{\!\cJ}\|_{\op,\mbH}=1$ and the definitions of $\gammaJ$ and $\hbeta_\cJ$, we have
	\begin{multline*}
		\bigg\| \sum_{d=1}^D \cR_{d,\cJ} (\hcH-\cH) \cP_{\!\cJ_d} \hcP_{\!\cJ} \bigg\|_{\HS,\mbH}^2
		\leq \bigg\| \sum_{d=1}^D \cR_{d,\cJ} (\hcH-\cH) \cP_{\!\cJ_d} \bigg\|_{\HS,\mbH}^2
		\\ = \sum_{d=1}^D \sum_{\ell \in \cJ_d} \sum_{m \notin \cJ} \frac{\|\cQ_m (\hcH-\cH) \cQ_\ell \|_{\HS,\mbH}^2}{(\lambda_\ell-\lambda_m)^2}
		\leq \frac{1}{\gammaJ} \sum_{\ell \in \cJ} \sum_{m \notin \cJ} \frac{\|\cQ_m (\hcH-\cH) \cQ_\ell \|_{\HS,\mbH}^2}{\min_{k \in \cJ} |\lambda_\ell-\lambda_k|}
		= \hbeta_\cJ^2.
	\end{multline*}
	Combining the above and using the triangle inequality, we obtain
	\begin{equation}
		\sqrt{ \sum_{k \in \cJ} \Big\| \hpsi_k - \sum_{\ell \in \cJ} \langle \hpsi_k, \psi_\ell \rangle_\mbH \psi_\ell
			\Big\|_\mbH^2 }
		\leq \bigg( \frac{37\sqrt{D}}{18} + 3 \bigg) \hbeta_{\!\cJ}.
		\label{gauss:spec:eigvec1-eq1}
	\end{equation}	
	Hence, the claim follows by Lemma~\ref{gauss:spec:perturb1} and the assumption $D=1$.
\end{proof}

\begin{proposition} \label{gauss:spec:eigvec2}
	Let $V_{\!\cJ} = (\langle \hpsi_k, \psi_\ell \rangle_\mbH)_{k,\ell \in \cJ}$.
	Under the same setting as in Lemma~\ref{gauss:spec:perturb1}, if $D=1$ and $\gammaJ^{-1} (2 \epsilon_1 \theta + \epsilon_2) \leq 1/4$, we have
	\begin{equation*}
		\| \hV_{\!\cJ}^\top \hV_{\!\cJ} - I_{|\!\cJ|} \|_\F \leq 4 \gammaJ^{-2} \sum_{k \in \cJ} \|(\hcH-\cH)\psi_k\|^2_\mbH.
	\end{equation*}
\end{proposition}

\begin{proof}[Proof of Proposition~\ref{gauss:spec:eigvec2}]
	According to Lemma~5.5 in \cite{bae2026asymptotic}, under $\halpha_{\!\cJ} \leq 1/4$, it holds that
	\begin{equation}
		\| \hV_{\!\cJ}^\top \hV_{\!\cJ} - I_{|\!\cJ|} \|_\F \leq 4 \hbeta_\cJ^2.
		\label{gauss:spec:eigvec2-eq1}
	\end{equation}
	Hence, the claim follows by Lemma~\ref{gauss:spec:perturb1}.
\end{proof}

\begin{proposition} \label{gauss:spec:eigvec3}
	For any $\{u_\ell\}_{\ell \in \cJ} \subset \mbH$, it holds that
	\begin{multline*}
		\bigg| \sum_{k \in \cJ} \sum_{\ell \in \cJ} \langle \hpsi_k, u_\ell \rangle_\mbH^2	- \sum_{k \in \cJ} \sum_{\ell \in \cJ} \langle \psi_k, u_\ell \rangle_\mbH^2 \bigg|
		\\ \leq \bigg| \sum_{\ell \in \cJ} \langle \hcS_{\!\cJ} u_\ell, u_\ell \rangle_\mbH \bigg|
		+ \bigg( \sum_{\ell \in \cJ} \|u_\ell\|_\mbH^2 \bigg) (8D +32) \gammaJ^{-2} \|\hcH-\cH\|_{\op,\mbH}^2.
	\end{multline*}
	where $\hcS_{\!\cJ}$ is defined by
	\begin{equation*}
		\langle \hcS_{\!\cJ} \psi_k, \psi_\ell \rangle_\mbH
		= \begin{cases}
			\displaystyle\frac{ \langle (\hcH-\cH) \psi_k, \psi_\ell \rangle_\mbH }{\lambda_k - \lambda_\ell},
			& \text{if } k \in \cJ, \ell \notin \cJ, \\
			\displaystyle\frac{ \langle (\hcH-\cH) \psi_k, \psi_\ell \rangle_\mbH }{\lambda_\ell - \lambda_k},
			& \text{if } k \notin \cJ, \ell \in \cJ, \\
			0, & \text{otherwise,}
		\end{cases}
	\end{equation*}
	and it satisfies
	\begin{equation*}
		\| \hcS_{\!\cJ} \|_{\op,\mbH} \leq 2 \sqrt{D} \gammaJ^{-1} \|\hcH-\cH\|_{\op,\mbH}.
	\end{equation*}
\end{proposition}

\begin{proof}[Proof of Proposition~\ref{gauss:spec:eigvec3}]
	According to Theorem~3.4 and Lemma~5.4 in \cite{bae2026asymptotic}, it holds that
	\begin{gather}
		\| \hcP_{\!\cJ} - \cP_{\!\cJ} - \hcS_{\!\cJ} \|_{\op,\mbH} \leq (8D +32) \halpha^2_\cJ,
		\quad \|\hcS_\cJ \|_{\op,\mbH} \leq 2\sqrt{D} \halpha_{\!\cJ},
		\label{gauss:spec:eigvec3-eq1}		
	\end{gather}
	where $\hcP_{\!\cJ} = \sum_{k\in\cJ} \hpsi_k \otimes_\mbH \hpsi_k$. From the identities
	\begin{align*}
		\langle \hcP_{\!\cJ} u_\ell, u_\ell \rangle_\mbH
		= \sum_{k \in \cJ} \langle \hpsi_k, u_\ell \rangle_\mbH^2, \quad
		\langle \cP_{\!\cJ} u_\ell, u_\ell \rangle_\mbH
		= \sum_{k \in \cJ} \langle \psi_k, u_\ell \rangle_\mbH^2,
	\end{align*}
	we obtain
	\begin{equation}
		\bigg| \sum_{k \in \cJ} \sum_{\ell \in \cJ} \langle \hpsi_k, u_\ell \rangle_\mbH^2	- \sum_{k \in \cJ} \sum_{\ell \in \cJ} \langle \psi_k, u_\ell \rangle_\mbH^2 \bigg|
		\leq \bigg| \sum_{\ell \in \cJ} \langle \hcS_{\!\cJ} u_\ell, u_\ell \rangle_\mbH \bigg|
		+ \bigg| \sum_{\ell \in \cJ} \langle (\hcP_{\!\cJ} - \cP_{\!\cJ} - \hcS_{\!\cJ}) u_\ell, u_\ell \rangle_\mbH \bigg|.
		\label{gauss:spec:eigvec3-eq2}
	\end{equation}
	Since $\langle (\hcP_{\!\cJ} - \cP_{\!\cJ} - \hcS_{\!\cJ}) u_\ell, u_\ell \rangle_\mbH \leq \| \hcP_{\!\cJ} - \cP_{\!\cJ} - \hcS_{\!\cJ} \|_{\op,\mbH} \|u_\ell\|_\mbH^2$,
	the claim follows by combining \eqref{gauss:spec:eigvec3-eq1} and \eqref{gauss:spec:eigvec3-eq2}, and applying Lemma~\ref{gauss:spec:perturb2}.
\end{proof}

\section{Spectral convergence of kernel Gram matrices} \label{gauss:mercer}

In this section, we briefly summarize the spectral convergence results for kernel Gram matrices established in Section~4 of \cite{bae2026asymptotic} and reformulate them suited to our analysis.

Given an i.i.d. sample $X_1,\ldots,X_n$ from a compact metric space $\M$, equipped with a probability measure $\P$, and a symmetric, continuous, and positive semi-definite kernel $h: \M \times \M \to \R$, we consider the corresponding kernel Gram matrix $\hbfH_n \in \R^{n \times n}$ and kernel integral operator $\bfH:\LP \to \LP$, defined as
\begin{equation*}
	\hbfH_n = \bigg( \frac{h (X_i, X_{i'})}{n} \bigg)_{i,i' \in [n]} \quad\text{and}\quad
	\bfH f (x) = \int_\M h(x,y) f(y) d \P(y), \hquad\text{respectively.}
\end{equation*}

In the above setting, Mercer's theorem guarantees the existence of the associated reproducing kernel Hilbert space $\mbH$ equipped with the inner product
\begin{equation*}
	\langle f, g \rangle_\mbH = \sum_{k=1}^\infty \frac{\langle f, \phi_k \rangle_\LP \langle g, \phi_k \rangle_\LP}{\lambda_k}.
\end{equation*}
$\mbH$ is a subset of $\LP$ and contains $\{\phi_k\}_{k \in \N}$ and $\{h(x, \cdot)\}_{x \in \M}$.
In particular, it has the reproducing property, that is, $\langle h (x,\cdot), f \rangle_\mbH = f(x)$ for all $f \in \mbH, x \in \M.$

The discrepancy between the finite-dimensional matrix $\bfH_n$ and the (possibly infinite dimensional) operator $\bfH$ can be bridged through the (uncentered) covariance operators $\cH_n$ and $\hcH_n$ on $\mbH$, defined by
\begin{equation*}
	\hcH_n = \frac{1}{n} \sum_{i=1}^n h (X_i, \cdot) \otimes_\mbH h (X_i, \cdot), \quad
	\cH = \E h (X, \cdot) \otimes_\mbH h (X, \cdot).	
\end{equation*}
These covariance operators are closely related to $\hbfH_n$ and $\bfH$, respectively, via the reproducing property.
Specifically, for all $f,g \in \mbH$, it holds that 
\begin{gather*}
	\hbfH_n \sn{f} = \sn{\hH_n f}, \quad 
	\bfH f = \cH f, \\
	\langle \sn{f}, \sn{g} \rangle_\LP = \langle \hcH_n f, g \rangle_\mbH, \quad 
	\langle f, g \rangle_\LP = \langle \cH f, g \rangle_\mbH, \\
	\| \sn{f} \|_\LPn \leq \|\hcH_n\|_{\op,\mbH} \|f\|_\mbH^2, \quad
	\| f \|_\LP \leq \|\cH\|_{\op,\mbH} \|f\|_\mbH^2.
\end{gather*}
Moreover, this relationship extends to their spectral properties.

We let $\{\hlambda_{k,n}, \hvarphi_{k,n}\}_{k \in [n]}$ and $\{\lambda_k, \phi_k\}_{k \in \N}$ denote the eigenpairs of $\hbfH_n$ and $\bfH$, respectively. The eigenvalues are arranged in non-increasing order, and may be repeated according to their multiplicity.
The eigenvectors $\{\hvarphi_{k,n}\}_{k \in [n]}$ and eigenfunctions $\{\phi_k \}_{k \in \N}$ are $\LPn$- and $\LP$-orthonormal, that is, $\langle \hvarphi_{k,n}, \hvarphi_{\ell,n} \rangle_\LPn = \langle \phi_k, \phi_\ell \rangle_\LP = \delta_{k\ell}$.

Then, $\hbfH_n$ (resp. $\bfH$) shares the same eigenvalues with $\hcH_n$ (resp. $\cH$).
Moreover, there exist $\mbH$-orthonormal eigenfunctions $\{\hpsi_{k,n}\}_{k\in[n]}$ (resp. $\{\psi_k\}_{k\in\N}$) of $\hcH_n$ (resp. $\cH$) satisfying
\begin{equation*}
	\sn{\hpsi_{k,n}} =  \sqrt{\hlambda_{k,n}} \hvarphi_{k,n}
	\quad
	( \text{resp. } \psi_k = \sqrt{\lambda_k} \phi_k).
\end{equation*}
Hence, studying the spectral convergence of $\hcH_n$ to $\cH$ directly implies that of $\hbfH_n$ to $\bfH$.
This in particular allows us to directly apply several results stated in Appendix~\ref{gauss:spec}.

For $\cJ \subset [n]$, the goal of this section is to investigate the spectral quantities of $\hbfH$ and $\bfH$ corresponding to $\cJ$.
We use the same assumptions on $\cJ$ and the same notation $\hcH_n$ and $\cH$ as in Appendix~\ref{gauss:spec}, e.g., $\gammaJ >0$, $\cP_\cJ$ and $\theta_\tmax$.
For simplicity of exposition, we assume $\gammaJ\le \theta_\tmin$, which typically holds when $\LP$ is infinite-dimensional and $n$ is large enough.

\subsection{Weighted perturbation based quantities}

As in Appendix~\ref{gauss:spec}, we consider the weighted perturbation quantities $\halpha_{\cJ,n}$ and $\hbeta_{\cJ,n}$, defined in the same way using $\cJ$, $\hcH_n$ and $\cH$.
The randomness of $\halpha_{\cJ,n}$ and $\hbeta_{\cJ,n}$ can be quantified in terms of the following:
\begin{align*}
	\bfkappa_{\!\cJ}
	&= \sup_{x \in \M} \bigg( \frac{1}{\gammaJ}\sum_{k \in \cJ} \lambda_k \phi_k^2(x) + \sum_{k \notin \cJ} \frac{ \lambda_k \phi_k^2 (x) }{\min_{m \in \cJ} |\lambda_k - \lambda_m|} \bigg), \\
	\bfsigma_{\!\cJ}
	&= \frac{1}{\gammaJ}\sum_{k \in \cJ} \lambda_k + \sum_{k \notin \cJ} \frac{ \lambda_k }{\min_{m \in \cJ} |\lambda_k - \lambda_m|}, \\
	\bfvarkappa_{\!\cJ}
	&= \sup_{x \in \M} \sqrt{ \frac{1}{\gammaJ}\sum_{k \in \cJ} \lambda_k \phi_k^2(x) \sum_{\ell \notin \cJ} \frac{ \lambda_\ell \phi_\ell^2 (x) }{\min_{m \in \cJ} |\lambda_\ell - \lambda_m|} }, \\
	\bfvarsigma_{\!\cJ}
	&= \E \Bigg[ \frac{1}{\gammaJ}\sum_{k \in \cJ} \lambda_k \phi_k^2(X) \sum_{\ell \notin \cJ} \frac{ \lambda_\ell \phi_\ell^2 (X) }{\min_{m \in \cJ} |\lambda_\ell - \lambda_m|} \Bigg].
\end{align*}
To control these quantities, we use $\varrho := \sup_{x\in\M} h(x,x)$. From the spectral decomposition of the kernel $h(x,y) = \sum_{k=1}^\infty \lambda_k \phi_k (x) \phi_\ell (x)$, it can be easily shown that
\begin{equation}
	\bfsigma_{\!\!\cJ}, \bfkappa_{\!\cJ}\leq \frac{\varrho}{\gammaJ}, \quad
	\bfvarkappa_{\!\cJ} \leq \frac{\varrho}{2\gammaJ}, \quad
	\bfvarsigma_{\!\cJ} \leq \frac{\varrho}{\gammaJ^2}\sum_{k \in \cJ} \lambda_k.
	\label{gauss:mercer:condition}
\end{equation}
With $\varrho$, we formulate the lemmas below.

\begin{lemma} \label{gauss:mercer:perturb}
	For $\,0<\tau<1$, we have
	\begin{gather} 
		\P \bigg( \halpha_{\!\cJ,n} \le \sqrt{\frac{4\theta_\tmax \varrho \log ( 2 \varrho /(\gammaJ \tau) ) }{n \gammaJ^2}} + \frac{2 \varrho \log ( 2 \varrho /(\gammaJ\tau) ) }{3n\gammaJ} \bigg) \geq 1- \tau,
		\label{gauss:mercer:perturb-res1} \\
		\P \bigg( \hbeta_{\!\cJ,n} \le \sqrt{ \frac{2 \varrho \sum_{k \in \cJ} \lambda_k \log(2/\tau)}{n\gammaJ^2}} + \frac{\varrho \log(2/\tau)}{6n\gammaJ} \bigg) \geq 1 -\tau,
		\label{gauss:mercer:perturb-res2} \\
		\P ( \hbeta_{\!\cJ,n} > 1 ) 
		\leq 2 \exp \bigg(- \min \bigg\{ \frac{n \gammaJ^2}{8 \varrho \sum_{k \in \cJ} \lambda_k }, \frac{3n \gammaJ}{\varrho} \bigg\} \bigg).
		\label{gauss:mercer:perturb-res3}
	\end{gather}	
	In particular, we have for $n \ge \frac{38 \theta_\tmax \varrho}{\gammaJ^2}$,
	\begin{equation}
		\P \Big( \halpha_{\!\cJ,n} > \frac{1}{4}\,\Big) 
		\leq \frac{2 \varrho}{\gammaJ} \exp \bigg( -\frac{n \gammaJ^2}{70 \theta_\tmax \varrho} \bigg).
		\label{gauss:mercer:perturb-res4}
	\end{equation}	
\end{lemma}

\begin{proof}[Proof of Lemma~\ref{gauss:mercer:perturb}]
	According to Proposition~4.4 in \cite{bae2026asymptotic}, the following hold:
	
	For $\,0<\tau<1$, it holds that
	\begin{gather*} 
		\P \bigg( \halpha_{\!\cJ,n} \le \sqrt{\frac{4\theta_\tmax \bfkappa_{\!\cJ} \log ( 2\bfsigma_{\!\!\cJ} /\tau ) }{\gammaJ \,n}} + \frac{2\bfkappa_{\!\cJ} \log ( 2\bfsigma_{\!\!\cJ} /\tau ) }{3n} \bigg) \geq 1- \tau,
		\\
		\P \bigg( \hbeta_{\!\cJ,n} \le \sqrt{\frac{2\bfvarsigma_\cJ \log(2/\tau)}{n}} +  \frac{\bfvarkappa_{\!\cJ}\log(2/\tau)}{3n} \bigg) \geq 1 -\tau. 
	\end{gather*}	
	and for $n \ge \frac{38 \theta_\tmax \bfkappa_{\!\cJ}}{\gammaJ}$,
	\begin{equation*}
		\P \Big( \halpha_{\!\cJ,n} > 1/4 \Big) 
		\leq 2 \bfsigma_{\!\!\cJ} \exp \bigg( -\frac{n \gammaJ}{70 \theta_\tmax \bfkappa_{\!\cJ}} \bigg).
	\end{equation*}
	
	Therefore, \eqref{gauss:mercer:perturb-res1}, \eqref{gauss:mercer:perturb-res2} and \eqref{gauss:mercer:perturb-res4} follows by applying \eqref{gauss:mercer:condition}. Lastly, \eqref{gauss:mercer:perturb-res3} is obtained by a standard inversion argument applied to \eqref{gauss:mercer:perturb-res2}.	
\end{proof}

\subsection{Eigenvalues of kernel Gram matrices}

In this part and the following one, we present spectral convergence results for eigenvalues and eigenvectors, respectively.

\begin{proposition} \label{gauss:mercer:eigval1}
	If $n \ge \frac{38 \varrho\,\theta_\tmax }{\gammaJ^2}$, we have with probability at least 
	$1-\tau-\frac{2\varrho}{\gammaJ} \exp \Big( -\frac{n\gammaJ^2}{70 \varrho\,\theta_\tmax} \Big) $,
	\begin{multline*}
		\bigg| \sum_{k\in\cJ} \hlambda_{k,n} - \sum_{k\in\cJ} \lambda_k
		- (\Pn-\P) \Big( \sum_{k \in \cJ} \lambda_k \phi_k^2 \Big) \bigg| 
		\\ \leq \bigg(\frac{5}{2\gammaJ} + \frac{4\theta_\tmax-4\theta_\tmin}{\gammaJ^2}\bigg) 
		\bigg( \sqrt{ \frac{2 \varrho \sum_{k \in \cJ} \lambda_k \log(2/\tau)}{n}} +  \frac{\varrho \log(2/\tau)}{6n} \bigg)^2
	\end{multline*}
\end{proposition}

\begin{proof}[Proof of Proposition~\ref{gauss:mercer:eigval1}]	
	According to Theorem~4.13 of \cite{bae2026asymptotic}, if $n \ge \frac{38 \theta_\tmax \bfkappa_{\cJ}}{\gammaJ}$, we have
	with probability at least $1-\tau- 2 \bfsigma_{\!\!\cJ} \exp \Big( -\frac{n\gammaJ}{70 \theta_\tmax \bfkappa_{\!\cJ_d}} \Big),$
	\begin{multline*}
		\bigg| \sum_{k\in\cJ} \hlambda_{k,n} - \sum_{k\in\cJ} \lambda_k
		- (\Pn-\P) \Big( \sum_{k \in \cJ} \lambda_k \phi_k^2 \Big) \bigg| 
		\\ \leq \bigg( \frac{5}{2}\gammaJ + 4\theta_\tmax-4\theta_\tmin\bigg)
		\bigg( \sqrt{\frac{2\bfvarsigma_\cJ \log(2/\tau)}{n}} +  \frac{\bfvarkappa_{\!\cJ}\log(2/\tau)}{3n} \bigg)^2.
	\end{multline*}
	Hence, the claim follows by applying \eqref{gauss:mercer:condition}.
\end{proof}

\begin{proposition} \label{gauss:mercer:eigval2}
	Suppose $n \ge \frac{38 \varrho\,\theta_\tmax }{\gammaJ^2}$. 
	We have with probability at least
	$1-\tau-\frac{2\varrho}{\gammaJ} \exp \Big( -\frac{n\gammaJ^2}{70 \varrho\,\theta_\tmax} \Big) $,
	\begin{gather*}
		\min_{k \in \cJ, \ell \notin \cJ}|\hlambda_k - \hlambda_\ell|  \ge \frac{\gammaJ}{2}, \\
		\max_{k \in \cJ}|\hlambda_{k,n,\e} - \lambda_k| \leq \sqrt{\frac{4\theta_\tmax \varrho \log(2\varrho/(\gammaJ\tau))}{n}} + \frac{2\varrho\log(2\varrho/(\gammaJ \tau))}{3n}.
	\end{gather*}
\end{proposition}

\begin{proof}[Proof of Proposition~\ref{gauss:mercer:eigval2}]
	Using \eqref{gauss:spec:eigval3-eq1} and \eqref{gauss:spec:eigval3-eq2} in Lemma~\ref{gauss:spec:eigval2}, we have 
	\begin{equation*}
		\min_{k \in \cJ, \ell \notin \cJ}|\hlambda_k - \hlambda_\ell| \ge \frac{\gammaJ}{2}, \quad\text{and}\quad
		\max_{k \in \cJ}|\hlambda_{k,n,\e} - \lambda_k| \leq \gammaJ \halpha_{\!\cJ,n} 
	\end{equation*}	
	provided $\halpha_{\!\cJ,n} \leq 1/4$.
	Hence, the claim follows by Lemma~\ref{gauss:mercer:perturb}.
\end{proof}

\subsection{Eigenvectors of kernel Gram matrices}

\begin{proposition} \label{gauss:mercer:eigvec1}
	If $n \ge \frac{38 \varrho \theta_\tmax }{\gammaJ^2}$, we have with probability at least 
	$1-\tau-\frac{2\varrho}{\gammaJ} \exp \Big( -\frac{n\gammaJ^2}{70 \varrho\,\theta_\tmax} \Big) $,
	\begin{equation*}
		\Bigg| \sum_{k \in \cJ} \sum_{\ell \in \cJ} \langle \hvarphi_{k,n}, \sn{\phi_\ell} \rangle_\LP^2 - |\cJ| - (\P_n-\P) \Big( \sum_{k \in \cJ} \phi_k^2 \Big) \Bigg|
		\\ \leq \frac{11}{2\gammaJ^2} \bigg( \sqrt{ \frac{2 \varrho \sum_{k \in \cJ} \lambda_k \log(2/\tau)}{n}} +  \frac{\varrho \log(2/\tau)}{6n} \bigg)^2.
	\end{equation*}
\end{proposition}

\begin{proof}[Proof of Proposition~\ref{gauss:mercer:eigvec1}]
	According to Theorem~4.10 of \cite{bae2026asymptotic},
	if $n \ge \frac{38 \theta_\tmax \bfkappa_{\!\cJ}}{\gammaJ}$, we have with probability at least 
	$ 1-\tau - 2 \bfsigma_{\!\!\cJ} \exp \Big( -\frac{n \gammaJ}{70\theta_\tmax \bfkappa_{\!\cJ}} \Big),$
	\begin{equation*}
		\Bigg| \sum_{k \in \cJ} \sum_{\ell \in \cJ} \langle \hvarphi_{k,n}, \sn{\phi_\ell} \rangle_\LP^2 - |\cJ| - (\P_n-\P) \Big( \sum_{k \in \cJ} \phi_k^2 \Big) \Bigg|
		\\
		\leq \frac{11}{2}
		\bigg( \sqrt{\frac{2\bfvarsigma_\cJ \log(2/\tau)}{n}} +  \frac{\bfvarkappa_{\!\cJ}\log(2/\tau)}{3n} \bigg)^2.
	\end{equation*}
	Therefore, the claim follows by applying \eqref{gauss:mercer:condition}.
\end{proof}

\begin{proposition} \label{gauss:mercer:eigvec2}
	If $n \ge \frac{38 \varrho \theta_\tmax }{\gammaJ^2}$, with probability at least 
	\begin{equation*}
		1-\tau-\frac{2\varrho}{\gammaJ} \exp \Big( -\frac{n\gammaJ^2}{70 \varrho\,\theta_\tmax} \Big) 
		- 2 \exp \bigg(- \min \bigg\{ \frac{n \gammaJ^2}{8 \varrho \sum_{k \in \cJ} \lambda_k }, \frac{3n \gammaJ}{\varrho} \bigg\} \bigg)
	\end{equation*}
	we have for any symmetric matrix $A \in \R^{|\cJ|\times|\cJ|}$,
	\begin{multline}
		\bigg| \sum_{k \in \cJ} \hlambda_{k,n} \langle A \hvarphi_{k,n}, \hvarphi_{k,n} \rangle_\LPn
		- \sum_{k \in \cJ} \lambda_k \langle A \sn{\phi_k}, \sn{\phi_k} \rangle_\LPn \bigg|
		\\ \leq \frac{C_{|\cJ|} \|A\|_{\op,\LPn} \| \hbfH_n \|_{\op,\LPn}}{\gammaJ}
		\bigg( \sqrt{ \frac{2 \varrho \sum_{k \in \cJ} \lambda_k \log(2/\tau)}{n}} +  \frac{\varrho \log(2/\tau)}{6n} \bigg),
		\label{gauss:mercer:eigvec2-res1}
	\end{multline}
	where $C_{|\cJ|}$ is a constant depending only on $|\cJ|$.
\end{proposition}

\begin{proof}[Proof of Proposition~\ref{gauss:mercer:eigvec2}]
	Using $ \hlambda_{k,n}^{1/2} \hvarphi_{k,n} = \sn{\hpsi_{k,n}}$ and $\hlambda_{k,n}^{1/2} \phi_{k,n} = \psi_{k,n}$, we see that
	\begin{align*}
		\sqrt{\hlambda_{k,n} \hlambda_{\ell,n}} \langle A \hvarphi_{k,n}, \hvarphi_{\ell,n} \rangle_\LPn
		&= \langle A \sn{\hpsi_{k,n}}, \sn{\hpsi_{\ell,n}} \rangle_\LPn, \\
		\sqrt{\lambda_k \lambda_\ell} \langle A \sn{\phi_k}, \sn{\phi_\ell} \rangle_\LPn
		&= \langle A \sn{\psi_k}, \sn{\psi_\ell} \rangle_\LPn,
	\end{align*}
	Hence, by letting 
	\begin{align*}
		\hPsi_{\!\cJ\!,n} &= \big( \langle A \sn{\hpsi_{k,n}}, \sn{\hpsi_{\ell,n}} \rangle_\LPn \big)_{k,\ell\in\cJ}, \quad
		\Psi_{\!\cJ} = \big( \langle A \sn{\psi_{k}}, \sn{\psi_{\ell}} \rangle_\LPn \big)_{k,\ell \in \cJ},
	\end{align*}
	we have
	\begin{equation}
		\sum_{k \in \cJ} \hlambda_{k,n} \langle A \hvarphi_{k,n}, \hvarphi_{k,n} \rangle_\LPn
		- \sum_{k \in \cJ} \lambda_k \langle A \sn{\phi_k}, \sn{\phi_k} \rangle_\LPn
		= \tr(\hPsi_{\!\cJ,n} - \Psi_{\!\cJ}).
		\label{gauss:mercer:eigvec2-eq1}
	\end{equation}
	Define
	$\hV_{\!\cJ\!,n} = (\hv_{k\ell,n})_{k,\ell \in \cJ} = (\langle \hpsi_{k,n}, \psi_\ell \rangle_\mbH)_{k,\ell \in \cJ}$
	and
	$\hxi_{k,n} = \hpsi_{k,n} - \sum_{m \in \cJ} \hv_{km,n} \psi_m$.
	Then we have
	\begin{multline*}
		\langle A \sn{\hpsi_{k,n}}, \sn{\hpsi_{\ell,n}} \rangle_\LPn
		= \sum_{m_1 \in \cJ} \sum_{m_2 \in \cJ} \hv_{km_1,n} \hv_{\ell m_2,n} \langle A \sn{\psi_{m_1}}, \sn{\psi_{m_2}} \rangle_\LPn
		\\ + \sum_{m_1 \in \cJ} \hv_{km_1,n} \langle A \sn{\psi_{m_1}}, \sn{\hxi_{\ell,n}} \rangle_\LPn
		+ \sum_{m_2 \in \cJ} \hv_{\ell m_2,n} \langle A \sn{\hxi_{k,n}}, \sn{\psi_{m_2}} \rangle_\LPn
		+ \langle A \sn{\hxi_{k,n}}, \sn{\hxi_{\ell,n}} \rangle_\LPn.
	\end{multline*}
	By letting 
	$E_1 = \big( \langle A \sn{\hxi_{k,n}}, \sn{\psi_{\ell}} \rangle_\LPn \big)_{k,\ell \in \cJ}$
	and
	$E_2 = \big( \langle A \sn{\hxi_{k,n}}, \sn{\hxi_{\ell,n}} \rangle_\LPn \big)_{k,\ell \in \cJ}$,
	we obtain
	\begin{align*}
		\hPsi_{\!\cJ\!,n} - \Psi_{\!\cJ} 
		= \hV_{\!\cJ,n} \Psi_{\!\cJ} \hV_{\!\cJ,n}^\top - \Psi_{\!\cJ}
		+ \hV_{\!\cJ,n} E_1 + E_1^\top \hV_{\!\cJ,n}^\top + E_2,
	\end{align*}
	and consequently,
	\begin{align}
		|\tr(\hPsi_{\!\cJ\!,n} - \Psi_{\!\cJ})|
		&\leq |\tr(\hV_{\!\cJ,n} \Psi_{\!\cJ} \hV_{\!\cJ,n}^\top) - \tr(\Psi_{\!\cJ})| + 2 |\tr(\hV_{\!\cJ,n} E_1)| + |\tr(E_2)|
		\notag\\&\leq \|\Psi_{\!\cJ}\|_\F \|\hV_{\!\cJ,n}^\top\hV_{\!\cJ,n} - I_{|\cJ|}\|_\F + 2 \|\hV_{\!\cJ,n}\|_\F \|E_1\|_\F + |\tr(E_2)|.
		\label{gauss:mercer:eigvec2-eq2}
	\end{align}
	
	On the other hand, from \eqref{gauss:spec:eigvec1-eq1} and \eqref{gauss:spec:eigvec2-eq1}, $\hV_{\!\cJ,n}$ and $\hxi_{k,n}$ satisfy
	\begin{equation*}
		\|\hV_{\!\cJ\!,n}^\top \hV_{\!\cJ\!,n} - I_{|\cJ|} \|_\F \leq 4 \hbeta_{\cJ,n}^2, \quad
		\sqrt{\sum_{k\in\cJ} \|\hxi_{k,n} \|_\mbH^2 } \leq \bigg( \frac{65\sqrt{D}}{18} + 5\bigg) \hbeta_{\cJ,n},
	\end{equation*}
	provided $\halpha_{\!\cJ\!,n} \leq 1/4$.
	Together with Bessel's inequality, and $\langle \sn{f}, \sn{g} \rangle_\LPn = \langle \hcH_n f, g \rangle_\mbH$,
	we bound the terms on the right-hand side of \eqref{gauss:mercer:eigvec2-eq2} as
	\begin{align*}
		\|\Psi_{\!\cJ}\|_\F
		&= \sqrt{ \sum_{k\in\cJ} \sum_{\ell\in\cJ} \langle A \sn{\psi_{k}}, \sn{\psi_{\ell}} \rangle_\LPn^2 }
		\leq \|A\|_{\op,\LPn} \sum_{k\in\cJ} \|\sn{\psi_k}\|_\LPn^2
		\\&\leq \|A\|_{\op,\LPn} \|\hcH_n\|_{\op,\mbH} \sum_{k \in \cJ} \|\psi_k\|_\mbH^2
		= |\cJ| \cdot \|A\|_{\op,\LPn} \|\hcH_n\|_{\op,\mbH},
	\end{align*}	
	\begin{align*}
		\| \hV_{j,n} \|_\F
		&= \sqrt{\tr( \hV_{j,n}^\top \hV_{j,n} )}
		= \sqrt{\tr( \hV_{j,n}^\top \hV_{j,n} - I_{|\cJ|}) + |\cJ|}
		\leq \sqrt{ \sqrt{|\cJ|} \cdot \|\hV_{\!\cJ,n}^\top \hV_{\!\cJ,n} - I_{|\cJ|}\|_\F + |\cJ|}
		\\&\leq |\cJ|^{\frac{1}{4}} \sqrt{ \|\hV_{\!\cJ,n}^\top \hV_{\!\cJ,n} - I_{|\cJ|}\|_\F} + \sqrt{|\cJ|}
		\leq 2 |\cJ|^{\frac{1}{4}} \hbeta_{\cJ,n} + \sqrt{|\cJ|},
	\end{align*}
	\begin{multline*}
		\| E_1 \|_\F
		= \sqrt{\sum_{k \in \cJ} \sum_{\ell \in \cJ} \langle A \sn{\hxi_{k,n}}, \sn{\psi_{\ell}} \rangle_\LPn^2}
		\leq \|A\|_{\op,\LPn} \sqrt{\sum_{k \in \cJ} \| \sn{\hxi_{k,n}} \|_\LPn^2} \sqrt{ \sum_{\ell \in \cJ} \|\sn{\psi_{\ell}} \|_\LPn^2}
		\\ \leq \|A\|_{\op,\LPn} \|\hcH_n\|_{\op,\mbH} \sqrt{\sum_{k \in \cJ} \| \hxi_{k,n} \|_\mbH^2} \sqrt{ \sum_{\ell \in \cJ} \|\psi_{\ell} \|_\mbH^2}
		\leq \sqrt{|\cJ|} \cdot \|A\|_{\op,\LPn} \|\hcH_n\|_{\op,\mbH} \bigg( \frac{65\sqrt{D}}{18} + 5\bigg) \hbeta_{\cJ,n},
	\end{multline*}
	\begin{align*}
		\tr(E_2)
		&= \sum_{k\in\cJ} \langle A \sn{\hxi_{k,n}}, \sn{\hxi_{k,n}} \rangle_\LPn
		\leq \|A\|_{\op,\LPn} \sum_{k \in \cJ} \|\sn{\hxi_{k,n}}\|_\LPn^2
		\\&\leq \|A\|_{\op,\LPn} \|\hcH_n\|_{\op,\mbH} \sum_{k \in \cJ}  \|\hxi_{k,n}\|_\mbH^2
		\leq \|A\|_{\op,\LPn} \|\hcH_n\|_{\op,\mbH} \bigg( \frac{65\sqrt{D}}{18} + 5\bigg)^2 \hbeta_{\cJ,n}^2.
	\end{align*}
	Consequently, under $\halpha_{\!\cJ,n} \leq 1/4$ and $\hbeta_{\!\cJ,n} \leq 1$, using $D \leq |\cJ|$ yields
	\begin{multline}
		|\tr(\hPsi_{\!\cJ\!,n} - \Psi_{\!\cJ})| 
		\leq |\cJ| \cdot \|A\|_{\op,\LPn} \|\hcH_n\|_{\op,\mbH} \cdot 4 \hbeta_{\cJ,n}^2
		+ 2 (2 |\cJ|^{\frac{1}{4}} \hbeta_{\cJ,n} + \sqrt{|\cJ|}) \\ \times \sqrt{|\cJ|} \cdot \|A\|_{\op,\LPn} \|\hcH_n\|_{\op,\mbH} \bigg( \frac{65\sqrt{D}}{18} + 5\bigg) \hbeta_{\!\cJ,n}
		+ \|A\|_{\op,\LPn} \|\hcH_n\|_{\op,\mbH} \bigg( \frac{65\sqrt{D}}{18} + 5\bigg)^2 \hbeta_{\cJ,n}^2
		\\ \leq C_{|\cJ|} \|A\|_{\op,\LPn} \|\hcH_n\|_{\op,\mbH} \hbeta_{\!\cJ,n}.
		\label{gauss:mercer:eigvec2-eq3}
	\end{multline}
	Since $\hcH_n$ and $\hbfH_n$ shares the same eigenvalues, it holds that $\|\hcH_n\|_{\op,\mbH}=\|\hbfH_n\|_{\op,\LPn}$.
	Therefore, the claim follows by \eqref{gauss:mercer:eigvec2-eq3} and Lemma~\ref{gauss:mercer:perturb}
\end{proof}

The following two lemmas are used to prove Proposition~\ref{gauss:mercer:eigvec3} below.

\begin{lemma} \label{gauss:mercer:eigvec3:aux1}
	If $\halpha_{\!\cJ,n} \leq 1/4$, we have
	\begin{equation*}
		\| (\hcP_{\!\cJ,n} - \cP_{\!\cJ}) \cP_{\!\cJ} \|_{\HS,\mbH} 
		\leq \bigg(\frac{16\sqrt{D}}{9}+3 \bigg) \hbeta_{\!\cJ,n},
	\end{equation*}
	where $\hcP_{\cJ,n} = \sum_{k \in \cJ} \hpsi_{k,n} \otimes_\mbH \hpsi_{k,n}$.
\end{lemma}

\begin{proof}[Proof of Lemma~\ref{gauss:mercer:eigvec3:aux1}]
	From equation (40) in \cite{bae2026asymptotic}, it holds that if $\halpha_{\!\cJ} \leq 1/4$,
	\begin{equation}
		\| (\hcP_{\!\cJ,n}-\cP_{\!\cJ}-\hcS_{\!\cJ,n}) \cP_{\!\cJ} \|_{\HS,\mbH} 
		\leq \bigg(\frac{56\sqrt{D}}{9}+8 \bigg) \halpha_{\!\cJ} \hbeta_{\!\cJ}
		\leq \bigg(\frac{16\sqrt{D}}{9}+2 \bigg) \hbeta_{\!\cJ},
		\label{gauss:mercer:eigvec3:aux1-eq1}
	\end{equation}
	where 
	\begin{equation*}
		\hcS_{\!\cJ,n}
		= \sum_{k \in \cJ} \sum_{\ell \notin \cJ} \frac{\cQ_k (\hcH_n - \cH) \cQ_\ell + \cQ_\ell (\hcH_n - \cH) \cQ_k}{\lambda_k - \lambda_\ell}, \quad
		\cQ_k = \psi_k \otimes_\mbH \psi_k.	
	\end{equation*}
	Then we bound $\|\hcS_{\!\cJ,n} \cP_{\!\cJ}\|_{\HS,\mbH}$ as
	\begin{align}
		\|\hcS_{\!\cJ,n} \cP_{\!\cJ}\|_{\HS,\mbH}^2
		= \sum_{k \in \cJ} \sum_{\ell \notin \cJ} \frac{\|\cQ_\ell (\hcH_n-\cH)\cQ_k\|_{\HS,\mbH}^2}{(\lambda_k-\lambda_\ell)^2}
		\leq \frac{1}{\gammaJ}\sum_{k \in \cJ} \sum_{\ell \notin \cJ} \frac{\|\cQ_\ell (\hcH_n-\cH)\cQ_k\|_{\HS,\mbH}^2}{\min_{m \in\cJ} |\lambda_m-\lambda_\ell|}
		= \hbeta_{\!\cJ,n}^2.
		\label{gauss:mercer:eigvec3:aux1-eq2}
	\end{align}
	Thus, the proof is complete by combining \eqref{gauss:mercer:eigvec3:aux1-eq1} and \eqref{gauss:mercer:eigvec3:aux1-eq2}.
\end{proof}

\begin{lemma} \label{gauss:mercer:eigvec3:aux2}
	For any $t>0$, 
	\begin{equation*}
		\P \Big( \big\| \big( \phi_k(X_i) \phi_\ell(X_i) - \delta_{k \ell} \big)_{k,\ell \in \cJ} \big\|_\F \ge t \Big)
		\leq 2 \exp \bigg(-\frac{nt^2}{2(|\cJ| \max_{k \in \cJ} \|\phi_k\|^2_\infty + \sqrt{|\cJ|})^2} \bigg).
	\end{equation*}
	Equivalently, with probability at least $1-\tau$,
	\begin{equation*}
		\| \big( \phi_k(X_i) \phi_\ell(X_i) - \delta_{k \ell} \big)_{k,\ell \in \cJ} \|_\F 
		\le \big( |\cJ| \max_{k \in \cJ} \|\phi_k\|^2_\infty + \sqrt{|\cJ|} \,\big) \sqrt{\frac{2\log(2/\tau)}{n}}.		
	\end{equation*}
\end{lemma}

\begin{proof}[Proof of Lemma~\ref{gauss:mercer:eigvec3:aux2}]
	Note that
	\begin{align*}
		\big( \langle \sn{\phi_k}, \sn{\phi_\ell} \rangle_\LPn - \delta_{k\ell} \big)_{k,\ell \in \cJ}
		= \frac{1}{n}\sum_{i=1}^n \big( \phi_k(X_i) \phi_\ell(X_i) - \delta_{k \ell} \big)_{k,\ell \in \cJ},
	\end{align*}
	which is a sum of mean-zero i.i.d. matrices $\big( \phi_k(X_i) \phi_\ell(X_i) - \delta_{k \ell} \big)_{k,\ell \in \cJ}$.
	Since 
	\begin{align*}
		\| \big( \phi_k(X_i) \phi_\ell(X_i) - \delta_{k \ell} \big)_{k,\ell \in \cJ} \|_\F
		\leq \sqrt{|\cJ|} \max_{k \in \cJ} \|\phi_k\|_\infty^2 + |\cJ|,
	\end{align*} 	
	we complete the proof by applying Hoeffding's inequality in Hilbert spaces (see, e.g., Section~3.2 of \cite{martinez-taboada2026empirical}).
\end{proof}

\begin{proposition} \label{gauss:mercer:eigvec3}
	If $n \ge \frac{38 \varrho \theta_\tmax }{\gammaJ^2}$, with probability at least 
	\begin{multline*}
		1-\tau-\frac{2\varrho}{\gammaJ} \exp \Big( -\frac{n\gammaJ^2}{70 \varrho\,\theta_\tmax} \Big) 
		- 2 \exp \bigg(- \min \bigg\{ \frac{n \gammaJ^2}{8 \varrho \sum_{k \in \cJ} \lambda_k }, \frac{3n \gammaJ}{\varrho} \bigg\} \bigg)
		\\ - 2 \exp \bigg(-\frac{n}{ C_{|\cJ|} (\max_{k \in \cJ} \|\phi_k\|^2_\infty + 1)^2} \bigg),
	\end{multline*}
	we have for any symmetric matrix $A \in \R^{n \times n}$,
	\begin{multline}
		\bigg| \sum_{k \in \cJ} \langle A \hvarphi_{k,n}, \hvarphi_{k,n} \rangle_\LPn
		- \sum_{k \in \cJ} \langle A \sn{\phi_k}, \sn{\phi_k} \rangle_\LPn \bigg|
		\\ \leq \frac{C_{|\cJ|} \|A\|_{\op,\LPn} \Big( \sqrt{\theta_\tmin^{-1} \|\hbfH_n\|_{\op,\LPn}} + \theta_\tmin^{-1} \|\hbfH_n\|_{\op,\LPn} \Big)}{\gammaJ} \hspace{8em}
		\\ \times \bigg( \sqrt{ \frac{2 \varrho \sum_{k \in \cJ} \lambda_k \log(4/\tau)}{n}} +  \frac{\varrho \log(4/\tau)}{6n} \bigg)
		+ C_{|\cJ|} (\max_{k \in \cJ} \|\phi_k\|_\infty^2 + 1) \sqrt{\frac{\log(4/\tau)}{n}},
		\label{gauss:mercer:eigvec3-res1}
	\end{multline}
	where $C_{|\cJ|}$ is a constant depending only on $|\cJ|$.
\end{proposition}

\begin{proof}[Proof of Proposition~\ref{gauss:mercer:eigvec3}]
	We divide the proof into several steps.
	Since the proof is similar to that of Proposition~\ref{gauss:mercer:eigvec2}, we omit some details
	
	\mycase{1. Decomposition}
	We let
	\begin{gather*}
		\hPhi_{\!\cJ\!,n} = \big( \langle A \hvarphi_{k,n}, \hvarphi_{\ell,n} \rangle_\LPn \big)_{k,\ell\in\cJ}, \quad
		\Phi_{\!\cJ} = \big( \langle A \sn{\phi_{k}}, \sn{\phi_{\ell}} \rangle_\LPn \big)_{k,\ell \in \cJ}, \\
		\hu_{k\ell,n} = \langle \sn{\phi_k}, \hvarphi_{\ell,n} \rangle_\LPn, \quad
		\heta_{k,n} = \sn{\phi_k} - \sum_{\ell \in \cJ} \hu_{k\ell,n} \hvarphi_{\ell,n}, \quad
		\hU_{\cJ,n} = (\hu_{k\ell,n})_{k,\ell \in \cJ}, \\
		E_1 = \big( \langle A \heta_{k,n}, \hvarphi_{\ell,n} \rangle_\LPn \big)_{k,\ell \in \cJ}, \quad 
		E_2 = \big( \langle A \heta_{k,n}, \heta_{\ell,n} \rangle_\LPn \big)_{k,\ell \in \cJ}.
	\end{gather*}
	By using the same argument in the proof of Proposition~\ref{gauss:mercer:eigvec2}, with $\Phi_{\!\cJ}$ instead of $\hOmega_{\!\cJ}$, we obtain
	\begin{multline}
		\bigg| \sum_{k \in \cJ} \langle A \hvarphi_{k,n}, \hvarphi_{k,n} \rangle_\LPn
		- \sum_{k \in \cJ} \langle A \sn{\phi_k}, \sn{\phi_k} \rangle_\LPn \bigg|
		= |\tr (\Phi_{\!\cJ}-\hPhi_{\!\cJ\!,n})|
		\\ \leq C_{|\cJ|} \Big( \| \hPhi_{\!\cJ} \|_\F \|\hU_{\!\cJ,n}^\top\hU_{\!\cJ,n} - I_{|\cJ} \|_\F
		+ \|E_1\|_\F \sqrt{1+\|\hU_{\!\cJ,n}^\top\hU_{\!\cJ,n} - I_{|\cJ|} \|_\F}  + |\tr(E_2)\|_\F \Big).
		\label{gauss:mercer:eigvec3-eq1}
	\end{multline}
	Below, we bound each term.
	
	\mycase{2. Bounding $\hPhi_{\!\cJ}$}
	Noting $\{\hvarphi_{k,n}\}_{k \in \cJ}$ are $\LPn$-orthonormal, Bessel's inequality gives
	\begin{align}
		\| \hPhi_{\!\cJ} \|_\F
		= \sqrt{\sum_{k \in \cJ} \sum_{\ell \in \cJ} \langle A \hvarphi_{k,n}, \hvarphi_{\ell,n} \rangle_\LPn^2}
		\leq \sqrt{|\cJ|} \, \cdot \|A\|_{\op,\LPn}.
		\label{gauss:mercer:eigvec3-eq2}
	\end{align}
	
	\mycase{3. Bounding $\|\hU_{\!\cJ,n}^\top\hU_{\!\cJ,n} - I_{|\cJ|}\|_\F$}
	By the reproducing property and the definitions of $\hcP_{\!\cJ,n}$ and $\cP_{\!\cJ}$, we have
	\begin{align}
		\sn{ \hcP_{\!\cJ,n} \phi_k } 
		&= \sum_{m \in \cJ} \langle \phi_k, \hpsi_{m,n} \rangle_\mbH \sn{\hpsi_{m,n}}
		= \sum_{m \in \cJ} \hlambda_{m,n}^{-1/2} \langle \phi_k, \hcH_n \hpsi_{m,n} \rangle_\mbH \hvarphi_{m,n}
		\notag\\&= \sum_{m \in \cJ} \hlambda_{m,n}^{-1/2} \langle \sn{\phi_k}, \sn{\hpsi_{m,n}} \rangle_\LPn \hvarphi_{m,n}
		= \sum_{m \in \cJ} \langle \sn{\phi_k}, \hvarphi_{m,n} \rangle_\LPn \hvarphi_{m,n},
		\label{gauss:mercer:eigvec3-eq3}
	\end{align}
	and
	\begin{align*}
		\sn{ \hcP_{\!\cJ,n} \phi_k } 
		= \sn{ \cP_{\!\cJ} \phi_k } + \sn{ (\hcP_{\!\cJ,n}-\cP_{\!\cJ}) \phi_k } 
		= \sn{\phi_k} + \sn{ (\hcP_{\!\cJ,n}-\cP_{\!\cJ}) \phi_k } .
	\end{align*}
	Consequently, we obtain
	\begin{align}
		\heta_{k,n} = \sn{\phi_k} - \sum_{m \in \cJ} \langle \sn{\phi_k}, \hvarphi_{m,n} \rangle_\LPn \hvarphi_{m,n} 
		= - \sn{ (\hcP_{\!\cJ,n}-\cP_{\!\cJ}) \phi_k }.
		\label{gauss:mercer:eigvec3-eq4}
	\end{align}
	
	On the other hand, from the definition of $\heta_{k,n}$, we have for $k,\ell \in \cJ$,
	\begin{align*}
		\langle \sn{\phi_k}, \sn{\phi_\ell} \rangle_\LPn = \sum_{m \in \cJ} \langle \sn{\phi_k}, \hvarphi_{m,n} \rangle_\LPn \langle \hvarphi_{m,n}, \sn{\phi_\ell} \rangle_\LPn  
		+ \langle \heta_{k,n}, \sn{\phi_\ell} \rangle_\LPn,
	\end{align*}
	and this gives
	\begin{align*}
		\big( \langle \sn{\phi_k}, \sn{\phi_\ell} \rangle_\LPn \big)_{k,\ell \in \cJ}
		= \hU_{\!\cJ,n} \hU_{\!\cJ,n}^\top + \big( \langle \heta_{k,n}, \sn{\phi_\ell} \rangle_\LPn \big)_{k,\ell \in \cJ}.
	\end{align*}
	Consequently, we obtain
	\begin{align*}
		\| \hU_{\!\cJ,n} \hU_{\!\cJ,n}^\top - I_{|\cJ|} \|_\F 
		\leq \big\| \big( \langle \sn{\phi_k}, \sn{\phi_\ell} \rangle_\LPn - \delta_{k\ell} \big)_{k,\ell \in \cJ} \big\|_\F
		+ \big\| \big( \langle \heta_{k,n}, \sn{\phi_\ell} \rangle_\LPn \big)_{k,\ell \in \cJ} \big\|_\F.
	\end{align*}	
	The first term on the right-hand side above can be controlled by Lemma~\ref{gauss:mercer:eigvec3:aux2}, which will be applied in the final step.
	
	For the second term, we bound
	\begin{multline*}
		\big\| \big( \langle \heta_{k,n}, \sn{\phi_\ell} \rangle_\LPn \big)_{k,\ell \in \cJ} \big\|_\F
		= \sqrt{ \sum_{k \in \cJ} \sum_{\ell \in \cJ} \langle \heta_{k,n}, \sn{\phi_\ell} \rangle_\LPn^2 }
		\\= \sqrt{ \sum_{k \in cJ} \sum_{\ell \in \cJ} \langle \sn{ (\hcP_{\!\cJ,n}-\cP_{\!\cJ}) \phi_k }, \sn{\phi_\ell} \rangle_\LPn^2 }
		= \sqrt{|\cJ|}\theta_\tmin^{-1} \|\hcH_n\|_{\op,\mbH} \| (\hcP_{\!\cJ,n}-\cP_{\!\cJ})\cP_{\!\cJ} \|_{\HS,\mbH}.
	\end{multline*}
	Hence, applying Lemma~\ref{gauss:mercer:eigvec3:aux1} yields
	\begin{equation}
		\| \hU_{\!\cJ,n} \hU_{\!\cJ,n}^\top - I_{|\cJ|} \|_\F 
		\leq \big\| \big( \langle \sn{\phi_k}, \sn{\phi_\ell} \rangle_\LPn - \delta_{k\ell} \big)_{k,\ell \in \cJ} \big\|_\F
		+ \sqrt{|\cJ|} \theta_\tmin^{-1} \|\hcH_n\|_{\op,\mbH} \bigg(\frac{16\sqrt{D}}{9}+3\bigg) \hbeta_{\!\cJ,n}
		\label{gauss:mercer:eigvec3-eq5}
	\end{equation}
	
	\mycase{4. Bounding $\|E_1\|_\F$}
	Using Bessel's inequality, we bound
	\begin{align*}
		\|E_1\|_\F 
		= \sqrt{\sum_{k \in \cJ} \sum_{\ell \in \cJ} \langle A \heta_{k,n}, \hvarphi_{\ell,n} \rangle_\LPn^2 }
		\leq \|A\|_{\op,\LPn} \sqrt{\sum_{k \in \cJ} \| \heta_{k,n} \|_\LPn^2 }.
	\end{align*}
	On the other hand, by \eqref{gauss:mercer:eigvec3-eq4}, we have
	\begin{align*}
		\sqrt{\sum_{k \in \cJ} \| \heta_{k,n} \|_\LPn^2}
		= \sqrt{\sum_{k \in \cJ} \| \sn{ (\hcP_{\!\cJ,n}-\cP_{\!\cJ}) \phi_k } \|_\LPn^2}
		\leq \sqrt{\theta_\tmin^{-1} \|\hcH_n\|_{\op,\mbH}} \| (\hcP_{\!\cJ,n}-\cP_{\!\cJ}) \cP_{\!\cJ} \|_{\HS,\mbH},
	\end{align*}
	and applying Lemma~\ref{gauss:mercer:eigvec3:aux1} gives
	\begin{align}
		\sqrt{\sum_{k \in \cJ} \| \heta_{k,n} \|_\LPn^2} 
		\leq \sqrt{\theta_\tmin^{-1} \|\hcH_n\|_{\op,\mbH}} \bigg(\frac{16\sqrt{D}}{9}+3 \bigg) \hbeta_{\cJ,n}.
		\label{gauss:mercer:eigvec3-eq6}
	\end{align}
	Consequently, we obtain
	\begin{equation}
		\|E_1\|_\F \leq  \|A\|_{\op,\LPn} \sqrt{\theta_\tmin^{-1}\|\hcH_n\|_{\op,\mbH}} \bigg(\frac{16\sqrt{D}}{9}+3 \bigg) \hbeta_{\cJ,n}.
		\label{gauss:mercer:eigvec3-eq7}
	\end{equation}
	
	\mycase{5. Bounding $|\tr(E_2)|$}
	Using \eqref{gauss:mercer:eigvec3-eq6}, we have
	\begin{equation}
		|\tr(E_2)| 
		= \bigg| \sum_{k \in \cJ} \langle A \heta_{k,\e}, \heta_{k,\e} \rangle_\LPn \bigg|
		\leq \|A\|_{\op,\LPn} \theta_\tmin^{-1}\|\hcH_n\|_{\op,\mbH} \bigg(\frac{16\sqrt{D}}{9}+3 \bigg)^2 \hbeta^2_{\cJ,n}
		\label{gauss:mercer:eigvec3-eq8}
	\end{equation}
	
	\mycase{6. Concluding the proof}
	Combining the results from the previous steps, under $\halpha_{\!\cJ,n} \leq 1/4$, $\hbeta_{\!\cJ,n} \leq 1$ and $\big\| \big( \langle \sn{\phi_k}, \sn{\phi_\ell} \rangle_\LPn - \delta_{k\ell} \big)_{k,\ell \in \cJ} \big\|_\F \leq 1$, we have
	\begin{multline*}
		\bigg| \sum_{k \in \cJ} \langle A \hvarphi_{k,n}, \hvarphi_{k,n} \rangle_\LPn
		- \sum_{k \in \cJ} \langle A \sn{\phi_k}, \sn{\phi_k} \rangle_\LPn \bigg|
		\\ 
		\leq C_{|\cJ|} \|A\|_{\op,\LPn} \bigg( \big\| \big( \langle \sn{\phi_k}, \sn{\phi_\ell} \rangle_\LPn - \delta_{k\ell} \big)_{k,\ell \in \cJ} \big\|_\F +\theta_\tmin^{-1} \|\hcH_n\|_{\op,\mbH} \hbeta_{\cJ,n}
		\\ + \sqrt{\theta_\tmin^{-1} \|\hcH_n\|_{\op,\mbH} } \hbeta_{\!\cJ,n} \bigg).
	\end{multline*}
	Hence, the proof is complete by applying Lemmas~\ref{gauss:mercer:perturb} and \ref{gauss:mercer:eigvec3:aux2}.
\end{proof}

\section{Auxiliary results}

Here, we provide some useful results for proving our main results.

The lemma below provides the Taylor expansion-type expression for the matrix 
$$(I_m+E)^{-1/2}A(I_m + E)^{-1/2},$$
provided that $\|E\|_\op$ is sufficiently small.
This result can be extended to compact self-adjoint operators on a separable Hilbert space.
For our purposes, we adapt the proof of Lemma~5.4 of \cite{koltchinskii2000random}, replacing Frobenius-norm bounds with operator-norm bounds.

\begin{lemma} \label{gauss:sup:taylor}
	Let $A,E \in \R^{m\times m}$ be symmetric matrices.
	If $\|E\|_\op <1/2$, we have
	\begin{gather*}
		\bigg\| (I_m+E)^{-1/2}A(I_m + E)^{-1/2} - A \bigg\|_\op
		\leq \frac{5 \|A\|_\op \|E\|_\op}{2},\\
		\bigg\| (I_m+E)^{-1/2}A(I_m + E)^{-1/2} - A + \frac{AE + EA}{2} \bigg\|_\op 
		\leq \frac{17 \|A\|_\op \|E\|_\op^2}{4}.
	\end{gather*}
\end{lemma}

\begin{proof}[Proof of Lemma~\ref{gauss:sup:taylor}]
	We only prove the second assertion. The first can be shown in the same way.
	
	By letting
	$$ R = (I_m + E)^{-1/2} - \left( I_m -\frac{E}{2}\right), $$
	we have the following expansion:
	\begin{align*}
		(I_m+E)^{-1/2} A (I_m+E)^{-1/2} - A
		= \left(I_m -\frac{E}{2} + R \right) A \left( I_m -\frac{E}{2} + R \right) - A
		= -\frac{AE + EA}{2} + S,
	\end{align*}
	where
	\begin{align*}
		S &= \frac{EAE}{4} + RA + AR -\frac{RAE+EAR}{2} + RAR.
	\end{align*}
	Hence, to complete the proof, it suffices to properly bound $\|S\|_\op$.
	
	An elementary calculation shows that 
	$$ \left|(1+x)^{-1/2} - \left(1-\frac{x}{2}\right) \right| \leq |x|^2 \quad\text{for}\hquad |x| < 1/2. $$
	Applying the matrix function argument (see, e.g., Theorem 6.2.8 of \cite{horn1991topics}) yields
	$$ \|R\|_\op = \left\| (I_m + E)^{-1/2} - \left( I_m -\frac{E}{2}\right) \right\|_\op \leq \|E\|_\op^2. $$
	Using this and the sub-multiplicity of operator norms, we have
	\begin{align*}
		\|S\|_\op 
		&\leq \left(\frac{1}{4}+1+1+\frac{2}{2}+1\right) \|A\|_\op \|E\|_\op^2
		= \frac{17\|A\|_\op \|E\|_\op^2}{4},
	\end{align*}
	which concludes the proof.
\end{proof}

Here we state concentration bounds for $U$-statistics stated below.
To obtain sharper bounds, we also consider the special case of $\P$-canonical functions.
A measurable function $f : \M^m := \underbrace{\M \times \dots \times \M}_{\text{$m$ times}} \to \R$ is called $\P$-canonical if
\begin{equation*}
	\E[ f(X_1, \dots, X_m) \mid X_1,\dots, X_{k-1}, X_{k+1}, \dots, X_m] =0 \quad\text{a.s. for each $k=1,\dots,m$.}
\end{equation*}
This degeneracy property allows the sums of $\P$-canonical functions to concentrate more tightly, leading to sharper high-probability bounds.

\begin{lemma}[Proposition~2.3 of \cite{arcones1993limit}] \label{gauss:sup:ustat}
	For i.i.d. random variables $X_1, \dots, X_n$ from a probability measure $\P$ on a measurable space $\M$,
	let $f : \M^m \to \R$ be a bounded measurable function with $\|f\|_\infty \leq B$, $\E f(X_1, \dots, X_m)=0$ and $\E f^2(X_1, \dots, X_m) \leq \sigma^2$. Then, we have with probability at least $1-\tau$,
	\begin{equation*}
		\left|\frac{1}{n^m} \sum_{(i_1,\dots,i_m)\in \mathsf{I}_m^n } f(X_{i_1}, \dots, X_{i_m}) \right|
		\le C_m \bigg( \sqrt{\frac{\sigma^2 \log(1/\tau)}{n}} + \frac{B \log(1/\tau)}{n} \bigg),
	\end{equation*}
	where $C_m>0$ is a constant depending only on $m$, and $\mathsf{I}_m^n$ is defined by
	$$ \mathsf{I}_m^n = \{ (i_1, \dots, i_m) \,:\, i_j \in \N,\, 1\leq i_j \leq n , i_j \ne i_k\text{ if }j \ne k \}. $$
	Moreover, if $f$ is $\P$-canonical, we have with probability at least $1-\tau$,
	\begin{equation*}
		\left|\frac{1}{n^m} \sum_{(i_1,\dots,i_m)\in \mathsf{I}_m^n } f(X_{i_1}, \dots, X_{i_m}) \right|
		\le C_m' \bigg( \frac{\sigma (\log(1/\tau))^{m/2}}{n^{m/2}} + \frac{B (\log(1/\tau))^{(m+1)/2}}{n^{(m+1)/2}} \bigg),
	\end{equation*}
	where $C_m'>0$ is a constant depending only on $m$.
\end{lemma}

\printbibliography[heading=bibintoc]

\end{document}